\documentclass{article}

\usepackage{amsthm}
\usepackage{amsfonts}
\usepackage{graphicx}
\usepackage{amsmath}
\usepackage{amssymb}
\usepackage{enumerate}

\usepackage[small,nohug,heads=vee]{diagrams}
\diagramstyle[labelstyle=\scriptstyle]

\DeclareGraphicsExtensions{.pdf,.png,.jpg} 
\usepackage[margin=3.5cm,nohead]{geometry}
\usepackage{epstopdf}

\makeatletter
\def\blfootnote{\gdef\@thefnmark{}\@footnotetext}
\makeatother

\newtheorem{teorema}{Theorem}[section]
\newtheorem{lema}[teorema]{Lemma}

\newtheorem{corolario}[teorema]{Corollary}
\newtheorem{proposicion}[teorema]{Proposition}

\newtheorem{asuncion}[teorema]{Assumption}
\newtheorem{notacion}[teorema]{Notation}

\newtheorem*{lema.sn}{Lemma}
\newtheorem*{def.sn}{Definition}
\newtheorem*{teoa}{Theorem A}
\newtheorem*{teob}{Theorem B}
\newtheorem*{teoc}{Theorem C}
\newtheorem*{teod}{Theorem D}
\newtheorem*{claim.sn}{Claim}

\theoremstyle{definition}
\newtheorem{remark}[teorema]{Remark}

\theoremstyle{definition}
\newtheorem{definicion}[teorema]{Definition}

\theoremstyle{definition}
\newtheorem{claim}[teorema]{Claim}

\newenvironment{prueba}[1][Proof.]{\begin{trivlist}
   \item[\hskip \labelsep {\itshape #1}]}{\rule{0.5em}{0.5em}\end{trivlist}}

\newcommand{\ol}{\overline}
\newcommand{\wt}{\widetilde}
\newcommand{\wh}{\widehat}
\newcommand{\ra}{\rightarrow}
\newcommand{\pr}{\partial}
\newcommand{\txt}{\textnormal}
\newcommand{\tl}{\tilde}
\newcommand{\empt}{\emptyset}
\newcommand{\e}{\epsilon}
\newcommand{\minus}{\setminus}

\newcommand{\R}{\mathbf R} 
\newcommand{\C}{\mathbf C}
\newcommand{\N}{\mathbf N}
\newcommand{\T}{\mathbf T}
\newcommand{\Z}{\mathbf Z}
\newcommand{\Q}{\mathbf Q}
\newcommand{\cl}{\mathcal}
\newcommand{\D}{\mathbf D}
\newcommand{\A}{\mathbf A}

\title{On annular maps of the torus and sublinear diffusion}

\author{Pablo D\'avalos}

\date{}

\begin{document}

\maketitle

\begin{abstract}
There is a classification by Misiurewicz and Ziemian \cite{mz} of elements in Homeo$_0(\T^2)$ by their rotation set $\rho$, according to wether $\rho$ is a point, a segment or a set with nonempty interior. A recent classification of nonwandering elements in Homeo$_0(\T^2)$ by Koropecki and Tal was given in \cite{st}, according to the itrinsic underlying ambient where the dynamics takes place: planar, annular and strictly toral maps. We study the link between these two classifications, showing that, even abroad the nonwandering setting, annular maps are characterized by rotation sets which are \textit{rational segments}. Also, we obtain information on the \textit{sublinear diffusion} of orbits in the -not very well understood- case that $\rho$ has nonempty interior.   
\end{abstract}

\tableofcontents

\blfootnote{This work was partially supported by FAPESP-Brasil and ITESO-M\'exico.}

\section{Introduction}

A very useful tool for understanding the dynamics of a homeomorphism $f$ of a surface $M$ in the isotopy class of the identity is the \textit{rotation set}; a topological invariant which is a subset of $H_1(M,\R)$ and roughly consists of asymptotic homological velocity vectors of points under iteration. 

It was first introduced by H. Poincar\'e \cite{po} as the \textit{rotation number} of circle homeomorphisms, and in this case he proved that, from the topological point of view, the rotation number leads to a complete classification of the dynamics. The concept was later extended for homeomorphisms of any manifold \cite{schw, mz, pol, lc2} or even for metric spaces \cite{mat}. For the case of surfaces there exist results relating the dynamics of a homeomorphism with its rotation set; for example showing that certain vectors with rational coordinates in the rotation set imply the existence of periodic orbits for $f$ \cite{f3,mat}, and showing that `large' rotation sets have positive topological entropy \cite{lm, mat}.  

In this article we deal with the case $M= \T^2$. For any lift $\wh{f}:\R^2\to\R^2$ of $f$, the rotation set of $\wh{f}$, denoted $\rho(\wh{f})$, is defined as the set of accumulation points of sequences of the form 
$$\left( \frac{\wh{f}^{n_i}(x_i)-x_i}{n_i} \right)_{i\in\N}$$
where $m_i\ra\infty$ and $x_i\in\R^2$, and it is known to be a compact, convex subset of $H_1(\T^2,\R) \simeq \R^2$ \cite{mz}. An interesting fact about the rotation set is that it prompts a classification of the set Homeo$_0(\T^2)$ of homeomorphisms of $\T^2$ homotopic to the identity into three (disjoint) cases, depending wether the rotation set is: 
\begin{enumerate}[(i)]
\item a set consisting of a single point;
\item a compact segment;
\item a set with nonempty interior. 
\end{enumerate} 

In case $(i)$ we say that $f$ is a \textit{pseudorotation}. Pseudorotations have been thoroughly studied. For example, in \cite{j1} it is given a Poincar\'e-like classification theorem for conservative pseudorotations, and in \cite{kt2} a classification theorem is given for \textit{rational} pseudorotations, that is, for the case that $\rho(\wh{f})$ is a single rational vector. Also, in \cite{kt1, fay1, bcjr} are constructed examples with exotic dynamical properties. 

In case $(iii)$ it is known that $f$ must have positive entropy \cite{lm}, any vector in $\rho(\wh{f})$ is realized by a minimal set \cite{mz}, and any rational vector in the interior of $\rho(\wh{f})$ is realized by a periodic orbit \cite{f3}. Thus, an interesting mechanism for creating entropy and (infinitely many) periodic orbits is given by the creation of (at least three) points with non-colinear asymptotic velocity vectors. Worth to mention, the only known examples in case $(iii)$ are polygons or ``infinite polygons'' with a countable set of extremal points \cite{kw2,kw}.

For case $(ii)$ it is conjectured in \cite{fm} that the only possible examples are \textit{rational segments}, meaning segments of rational slope containing rational vectors, and segments of irrational slope containing a rational endpoint. For the latter case, an example in \cite{handel} attributed to Katok gives essentially the unique mechanism known to the author to produce examples with such a rotation set. Examples with rational segments are an analog of twist maps for elements in Homeo$_0(\T^2)$, and they have been recently studied. For example, results concerning dynamical models and the existence of periodic orbits associated to any rational vector in $\rho(\wh{f})$ include \cite{kk, gkt, dav}.\\

Aside this classification given by the rotation set, a recent new classification for the nonwandering elements in Homeo$_0(\T^2)$ was given in \cite{st}. In that article, it is shown that any nonwandering element of Homeo$_0(\T^2)$ falls in one of the following categories:
\begin{enumerate}[(a)]
\item there is $k\in\N$ such that Fix$(f^k)$ is \textit{fully essential}, that is, $\T^2\minus\txt{Fix}(f^k)$ is a union of topological open discs. 
\item $f$ is \textit{annular}: there is $k\in\N$ and a lift $F:\R^2\to\R^2$ of $f^k$ such that the deviations in the direction of some nonzero $v\in\Z^2$ are uniformly bounded:
$$ -M \leq \left\langle  F^n(x)-x, v \right\rangle  \leq M \ \ \ \txt{for all $x\in\R^2$ and $n\in\Z$}.$$
\item $f$ is \textit{strictly toral}. This roughly means that the dynamics of $f$ cannot be embedded in the plane or the annulus, and also, the `irrotational' part of the dynamics is contained in `elliptic islands' (below we explain this concept precisely). 
\end{enumerate}     

In case $(a)$ the dynamics of $f$ takes place essentially in the plane, as the set $\T^2\minus\txt{Fix}(f^k)$ can be embedded in the plane. 

In case $(b)$ it is easy to see that there is a finite covering of $\T^2$ such that the lift of $f$ to this covering has an invariant annular set (see for example Remark 3.10 in \cite{j2}), so that in some sense the dynamics of $f$ in a finite covering is embedded in an annulus. The notion of annular homeomorphism is equivalent to saying that $f$ is \textit{rationally bounded} in the sense of \cite{j2}.

For case $(c)$, let us first explain the notion of \textit{essential points} for a nonwandering homeomorphism. An \textit{essential set} $K\subset\T^2$ is a set which is not contained in a topological open disc. A point $x\in\T^2$ is essential for $f$ if the orbit of every neighborhood of $f$ is an essential subset of $\T^2$. Roughly speaking, this means that $x$ has a weak form of rotational recurrence. The set of essential points for $f$ is denoted Ess$(f)$, and the set of inessential points is Ine$(f)=\T^2\minus\txt{Ess}(f)$. The precise statement of item $(c)$ is the following: the set Ess$(f)$ is nonempty, connected and fully essential, and Ine$(f)$ is the union of pairwise disjoint bounded open discs. Therefore, if $f$ is strictly toral, there is a decomposition of the dynamics into a union Ine$(f)$ (possibly empty) of periodic bounded discs which can be regarded as ``elliptic islands'', and a fully essential set Ess$(f)$ which carries the ``rotational'' part of the dynamics. 

Case $(c)$ is disjoint from cases $(a)$ and $(b)$. However, cases $(a)$ and $(b)$ intersect (the identity belongs to both). In order to distinguish ``planar from annular'' we make the following definition. We say that a homeomorphism $f$ is \textit{planar} if either $f$ belongs to case $(a)$ or the orbits of (any lift to $\R^2$ of) $f$ are uniformly bounded. Also, we say that $f$ is \textit{strong annular} if $f$ is annular and is not planar. 

We therefore have that the set of nonwandering elements of Homeo$_0(\T^2)$ is a disjoint union of planar, strong annular and strictly toral maps.\\

There exist some links of course between these two classifications in the nonwandering case. Let us see how cases $(i),(ii),(iii)$ fall into the categories of planar, strong annular and strictly toral. 

For case $(iii)$, namely if $\rho(\wh{f})$ has nonempty interior, it is clear that $f$ cannot be annular. It cannot either be planar, as in \cite{st} it is shown that planar maps are \textit{irrotational}: $\rho(\wh{f})=\{(0,0)\}$. Therefore in case $(iii)$ $f$ is strictly toral. 

As for case $(i)$, an \textit{irrational pseudorotation} (that is, when $\rho(\wh{f})$ is a totally irrational vector) is analogously seen to be strictly toral. If $f$ is a rational pseudorotation, then $f$ might be either planar (e.g. the identity), strictly toral (see $\S$1.2 in \cite{kt2} for an example\footnote{Such example is a non-wandering map with a unique invariant measure supported on a fixed point. In contrast, in \cite{tal1} it is shown that an irrational pseudorotation which preserves a measure with full suport cannot be strictly toral.}), and it is not known wether $f$ may be strong annular (see Question 3 in \cite{tal1}). Finally, if $\rho(\wh{f})$ is a vector which is neither rational nor totally irrational, then $f$ may be either strong annular (e.g. a rigid translation) or strictly toral (as is Furstenberg's example \cite{fur}). 

The options for case $(ii)$ are either strong annular or strictly toral. As we mentioned above, the only known examples are rational segments and segments with irrational slope containing a rational endpoint. Note that if any other example may exist, it cannot be annular and therefore it must be strictly toral. The case of a segment with irrational slope containing a rational endpoint is also strictly toral, and rational segments may of course be realized by annular maps (e.g. the twist-like map $(x,y)\mapsto (x,y+\sin(2\pi x)))$. The question wether rational segments may or may not be realized by strictly toral maps remained open.

\begin{figure}[h]        
\begin{center} 
\includegraphics{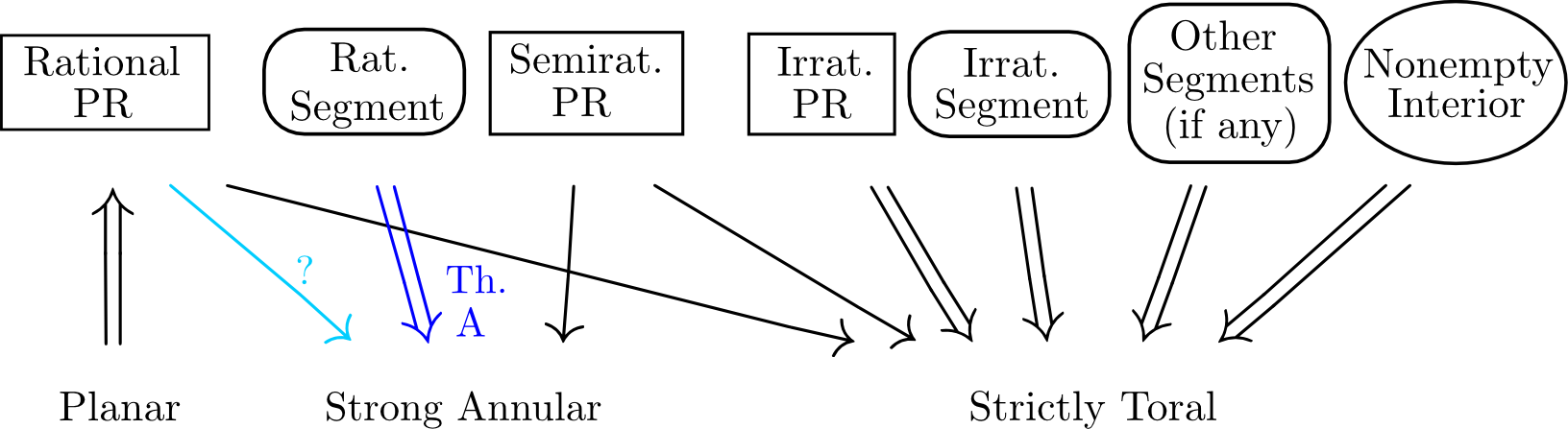}
\caption{Link between the classifications of Misiurewicz-Ziemian and Koropecki-Tal. Double arrows represent implications, and single arrows existence of examples. ``Semirat. PR'' stands for the pseudorotations which are neither rational nor irrational.}
\label{fig.intro}
\end{center}  
\end{figure}

Summarizing, the only two remaining open questions in this direction are the following:
\begin{enumerate}
\item If $f$ is a non-wandering pseudorotation, can $f$ be strong annular?
\item If $\rho(\wh{f})$ is a rational segment, can $f$ be strictly toral?
\end{enumerate}
A negative answer to Question 2 was given in \cite{gkt} in the area preserving case. \\

One of the objectives of this work is to contrubute to understanding the link between these two classifications by giving a negative answer to Question 2, \textit{even abroad the nonwandering setting}, namely, for any element of Homeo$_0(\T^2)$. This is given by the following. 

\begin{teoa}   
If $\rho(\wh{f})$ is a segment with rational slope containing rational points, then $f$ is annular. 
\end{teoa}

Another purpose of this article is to give a (very small) step in understanding case $(iii)$; the case when $\rho(\wh{f})$ has nonempty interior. As we said, the only known examples are polygons or ``infinite polygons'' with a countable set of extremal points. It is also not known if there are convex subsets of $\R^2$ which are \textit{not} realizable as rotation sets. The next theorem gives us information about the \textit{sublinear diffusion} of displacements away from the rotation set, in the case that $\pr \rho(\wh{f})\subset\R^2$ contains a rational segment $S$. Note that without loss of generality we may assume that the segment $S$ is contained in a line with rational slope passing through the origin in $\R^2$ (see Section \ref{sec.rotset} for the basic properties of the rotation set), and the next theorem then tells us that displacements must be uniformly bounded in the direction orthogonal to $S$ and outwards $\rho(\wh{f})$.

\begin{teob}    
If $\rho(\wh{f})$ is contained in the half-plane $\{ x\in\R^2 \ : \left\langle x,v \right\rangle \leq 0 \}$, for some $v\in\Q^2$, and if the line $v^{\perp}$ contains more than one point of $\rho(\wh{f})$, then there exists $M>0$ such that
$$  \langle  \wh{f}^n(x)-x , v \rangle  < M \ \ \ \txt{for all $x\in\R^2$ and $n\in\N$}. $$ 
\end{teob}

Finally, we adress the following question: \textit{to what extent does the rotation set capture the rotation information?} In \cite{kt1} is constructed a (smooth, Lebesgue-ergodic) example with $\rho(\wh{f})=\{(0,0)\}$ and such that almost every point rotates at a sublinear speed in \textit{every direction}, that is, for almost every $x\in\R^2$, the sequence
$$\left( \frac{\wh{f}^n(x)-x}{|\wh{f}^n(x)-x|} \right)_{n\in\N}$$ 
accumulates in the whole circle $\T^1$. In other words, the dynamics is far away to be rotationless, even though $\rho(\wh{f})$ is reduced to the point $\{(0,0)\}$. This phenomenon of sublinear diffusion may also occur for any pseudorotation \cite{kk2,j2}. 

For larger rotation sets the situation turns out to be quite different. For example, Theorem A tells us that if $\rho(\wh{f})$ is a rational segment, then \textit{there is no sublinear diffusion in the direction perpendicular to $\rho(\wh{f})$}. Moreover, if $\rho(\wh{f})$ is a rational polygon, that is, a non-degenerate polygon whose extremal points are rational vectors, we will show that \textit{there is no sublinear diffusion at all} (note that any rational polygon is realized as a rotation set \cite{kw2}). This is given by the following.

\begin{teoc}
If $\rho(\wh{f})$ is a non-degenerate polygon with rational endpoints, then there exists $M>0$ such that
$$ \sup_{x\in\R^2, n\in\N} \   d( \wh{f}^n(x)-x, n\rho(\wh{f}) )  < M. $$ 
\end{teoc}

We make the important remark that in the last three theorems, the bound $M = M(f) >0$ in them will be constructed \textit{explicitly}. Once the constant $M(f)$ is constructed, we will prove the following main theorem, which will have as corollaries Theorems A, B and C. It gives us our main dichotomy between \textit{bounded mean motion} and positive \textit{linear speed}. 

\begin{teod}   
Suppose that $\rho(\wh{f})$ contains two vectors $(0,a)$ and $(0,b)$, with $a<0<b$. 

Then, one of the following holds: 
\begin{itemize}
\item $\sup_{x\in\R^2, n\in\N} \ \txt{pr}_1(\wh{f}^n(x)- x)  < M(f)$;
\item $\exists x\in\R^2 \ \exists N>0$ s.t. $\txt{pr}_1(\wh{f}^{nN}(x)-x) > n \ \ \forall n\in\N$. 
\end{itemize}
\end{teod}

This article is organized as follows. In Section \ref{sec.prel} we introduce some basic facts about the rotation set, Atkinson's Lemma from ergodic theory, Brouwer theory, Poincar\'e-Bendixson theory and a theorem of Handel. In Section \ref{sec.bcd} we prove Theorems A, B and C from Theorem D, and Section \ref{sec.proof} is devoted to the proof of Theorem D.\\

\noindent \textbf{Acknowledgments}. I am very grateful to F. A. Tal for carefully listening to the proof of the main theorem, for improvements and for key suggestions. I am also grateful to F. Le Roux and P. Le Calvez for their hospitality during a short visit to Universit\'e Paris VII and for suggestions which led to simplifications in the proof of the main theorem. Finally, I also thank T. J\"ager for suggestions that led to improvements and A. Koropecki , A. Passeggi and J. Xavier for useful conversations.

\section{Preliminaries}     \label{sec.prel}

\subsection{Notations.}      \label{sec.notac}

We denote $\N$ the set of positive integers, and $\N_0=\N\cup \{0\}$. Also we denote the circle $\T^1=\R / \Z$ and the two-torus $\T^2=\R^2 / \Z$. If $S= \R\times\T^1$ or $\R^2$, we will denote the translation $T_1: S\to S$, $T_1(x,y)=(x+1,y)$, and $T_2:\R^2\to\R^2$ denotes $T_2(x,y)=(x,y+1)$. 

If $S$ is any surface, a map $\gamma:[a,b]\to S$ and its image will be both referred to as a \textbf{curve}. By an \textbf{arc} we mean a simple compact curve $\gamma:[0,1]\to S$. For the concatenation of two arcs $\gamma_1,\gamma_2:[0,1]\ra S$ we will use the `left-to-right' notation 
$$\gamma_1\cdot\gamma_2 = \left\{
\begin{array}{l l}
\gamma_1(2t) & t\in[0,1/2] \\
\gamma_2(2t-1) & t\in[1/2,1]
\end{array} \right.
$$
The arc $\gamma$ parametrized with the opposite orientation will be denoted by $\gamma^{-1}$. If $S=\R\times\T^1$ or $\T^2$, a \textbf{vertical curve} $\gamma\subset S$ is a loop which has homotopy class $(0,n)\in H_1(S,\Z)$, for some $n\in\Z$.

For an isotopy $(f_t)_t$ on a surface $S$, define its \textbf{canonical lift} to a covering space $\wt{S}$ as the lift $(\wt{f}_t)_t$ such that $\wt{f}_0=\txt{Id}$. For the homeomorphism $f=f_1$, we call the \textbf{canonical lift of $f$ with respect to the isotopy $(f_t)_t$} as the lift $\wt{f}=\wt{f}_1$. If the isotopy between Id and $f$ is implicitly known, we just refer to the \textbf{canonical lift of $f$}. If $(f_t)_{t\in[0,1]}$ is an isotopy, we will denote also by $(f_t)_t$ its extension for $t\in\R$, namely $f_t= f_{t \mod 1} \circ f^{\left\lfloor t \right\rfloor}$, where we use the notation 
$$ t \mod 1 = t - \left\lfloor t \right\rfloor \ \ \ \txt{and} \ \ \ 
\left\lfloor t \right\rfloor = \left\{ 
\begin{array}{l l}
\left\lfloor  t\right\rfloor & t\geq 0 \\
-\left\lfloor |t| \right\rfloor & t < 0
\end{array}  \right.$$

If $X\subset S$, two arcs $\gamma_1,\gamma_2\subset S\minus X$ are said to be \textbf{homotopic Rel$(X)$} if there is a homotopy on $S\minus X$ from $\gamma_1$ to $\gamma_2$. If $\gamma_1$ and $\gamma_2$ are homotopic with fixed endpoints, we will abreviate this by saying that $\gamma_1,\gamma_2$ are homotopic \textbf{wfe}. The \textbf{interior of an arc} $\gamma$, denoted int$(\gamma)$ is defined as $\gamma\minus \{\gamma(0),\gamma(1)\}$.

If  $S=\R\times\T^1$ or $\T^2$, a set $K\subset S$ is called \textbf{essential} if it is not contained in a topological open disc. If $K\subset S$ is not essential, we say that $K$ is \textbf{inessential}. If $K\subset\T^2$, a set is called \textbf{fully essential} if $K$ contains the complement of some disjoint union of topological open discs. If $K\subset S$ is compact and essential, we say that $K$ is \textbf{vertical} if $K$ is contained in a topological open annulus homotopic to the vertical annulus $\{(x,y)\,:\, 0< x < 1/2 \}$. Also, a set $K\subset S$ is called \textbf{annular} if it is a nested intersection of compact annuli $A_i$ such that the inclusion $A_{i+1} \hookrightarrow A_i$ is a homotopy equivalence. A \textbf{circloid} in an annular set which does not conain any proper annular subset. Let $K\subset S$ be a compact, connected, essential and vertical set, and let $\wh{K}\subset\R^2$ be a connected component of the preimage of $K$ by the (canonical) covering map $\R^2\to S$. As $K$ is vertical, the set $\R^2\minus \wh{K}$ has exactly one connected component $U_r$ which is unbounded to the right, and exactly one $U_l$ which is unbounded to the left. We denote
$$  R(\wh{K}) = U_r \ \ \ \and \ \ \ L(\wh{K})= U_l.$$

A \textbf{line} $\ell$ is a proper embedding $\ell:\R\ra\R^2$. Given a line $\ell$, by Shoenflies' Theorem (\cite{cairns}), there exists an orientation preserving homeomorphism $h$ of $\R^2$ such that $h\circ \ell(t)=(0,t)$, for all $t\in\R$. Then, the open half-plane $h^{-1}((0,\infty)\times\R)$ is independent of $h$, and we call it the \textbf{right} of $\ell$, and denote it by $R(\ell)$. Analogously, we define $L(\ell)=h^{-1}((-\infty,0)\times\R)$ the open half-plane to the \textbf{left} of $\ell$.  The sets $\ol{R}(\ell)$ and $\ol{L}(\ell)$ denote the closures of $R(\ell)$ and $L(\ell)$, respectively. If $\ell,\ell'$ are two lines in $\R^2$, we define $(\ell,\ell')=R(\ell)\cap L(\ell')$, and $[\ell,\ell']=\ol{R}(\ell)\cap \ol{L}(\ell)$. A \textbf{Brouwer curve} is a line $\ell$, such that $f(\ell)\subset R(\ell)$, and a \textbf{Brouwer $(0,1)$-curve} $\ell$ is a Brouwer curve whose image by the canonical projection $\R^2\to\R\times\T^1$ is a simple essential loop and such that $\ell$ is oriented upwards.  

Suppose that $K\subset S = \R^2,\R\times\T^1$ or $\T^2$ is an inessential set, and let $U\subset S$ be a topological open disc containing $K$. Let $\wh{U}\subset\R^2$ be a lift of $U$ by the canonical covering $\pi:\R^2\to S$, and let $\wh{K}= \pi^{-1}(K)\cap\wh{U}$. The \textbf{filling} of $K$, denoted Fill$(K)$ is defined as $\pi(\R^2\minus V)$, where $V$ is the unbounded connected component of $\R^2\minus \wh{K}$. A set $K\subset S$ such that $K= \txt{Fill}(K)$ is called a \textbf{filled} set. 

Consider a foliation $\cl{F}$ with singularities of $\R^2$, and denote sing$(\cl{F})$ the set of singularities of $\cl{F}$. If $\gamma\subset\R^2$ is a bounded leaf of $\cl{F}$, denote $\alpha(\gamma)$ and $\omega(\gamma)$ the alpha and omega-limit sets of $\gamma$, respectively. Define
$$\txt{sing}(\gamma) = \txt{sing}(\cl{F})\cap ( \txt{Fill}(\omega(\gamma)) \cup \txt{Fill}(\alpha(\gamma)) ).$$

For subsets of $\R$ or $\R^2$ we denote the metric $d(A,B)= \sup \{ |a-b| \,:\, a\in A, b\in B\}$, and when $A$ consists of a point $A=\{a\}$, we write $d(a,B)=d(A,B)$. We denote the \textbf{diameter} of a subset $A$ of $\R$ or $\R^2$ by diam$(A)$, and for a subset $A\subset\R^2$ the \textbf{horizontal diameter} diam$_1(A)=\txt{diam}(\txt{pr}_1(A))$. Similarly, if $A\subset\T^2$, we define diam$(A)= \txt{diam}(\pi^{-1}(A))$ and diam$_1(A)= \txt{diam}_1(\pi^{-1}(A))$, where $\pi:\R^2\to\T^2$ denotes the canonical projection. 

For two sets $A,B\subset\R^2$, we say that $A$ is \textbf{above} $B$ if $\inf \txt{pr}_2(A) > \sup \txt{pr}_2(B)$, and similarily, we say that $A$ is \textbf{below} $B$ if $B$ is above $A$.

\subsection{The rotation set.}    \label{sec.rotset}

Let $f:\T^2\to\T^2$ be a homeomorphism homotopic to the identity, and let $\wh{f}:\R^2\to\R^2$ be any lift of $f$. Consider the rotation set $\rho(\wh{f})$ of $\wh{f}$ as defined in the introduction. For a point $x\in\T^2$, we define its \textit{rotation set}, denoted $\rho(x,\wh{f})$, as the subset of $\R^2$ of accumulation points of sequences of the form
$$\left( \frac{\wh{f}^{n_i}(\hat{x})-\hat{x}}{n_i} \right)_{i\in\N}$$
where $m_i\to\infty$ and $\hat{x}\in\pi^{-1}(x)$. When $\rho(x,\wh{f})$ consists of a single vector $v$, we call it the \textit{rotation vector} of $x$. For any subset $K\subset\T^2$, we analogously define the \textit{rotation set} of $K$ by 
$$\rho(K,\wh{f}) = \cup_{x\in K} \rho(x,\wh{f}).$$
If $K\subset\R^2$, we denote $\rho(K,\wh{f})= \rho(\pi(K),\wh{f})$. 
  
It follows easily from the definition that 
\begin{equation}    \label{eq.pq}
\rho(T_1^p T_2^q \wh{f}^n) = T_1^p T_2^q( n \rho(\wh{f})).
\end{equation}
Under topological conjugacies, the rotation set behaves in the following way (see for example \cite{kk} for a proof).
\begin{lema}   \label{lema.rsconj}
Let $f:\T^2\to\T^2$ be a homeomorphism homotopic to the identity, let $A\in \txt{SL}(2,\Z)$ and let $h:\T^2\to\T^2$ be a homeomorphism isotopic to the map $\T^2\to\T^2$ induced by $A$. Let $\wh{f}$ and $\wh{h}$ be the respective lifts of $f$ and $h$. Then
$$ \rho(\hat{h}\hat{f}\hat{h}^{-1})= A\rho(\wh{f}).$$
In particular, $\rho(A\wh{f}A^{-1})= A\rho(\wh{f})$.
\end{lema}

Then, if $\rho(\wh{f})$ is a segment with rational slope, we may find $A\in \txt{SL}(2,\Z)$ such that $\rho(A\wh{f}A^{-1})$ is vertical. Indeed, if $\rho(\wh{f})$ has slope $p/q$, we may find integers $x,y$ such that $px+qy=1$, and then letting
$$ A= 
\begin{pmatrix}
p & -q \\
y & x
\end{pmatrix}
$$
we have det$(A)=1$, and as $A(q,p)=(0,1)$, $A\rho(\wh{f})$ is vertical.

By this, and using (\ref{eq.pq}) and Lemma \ref{lema.rsconj} one may easily show the following.

\begin{lema}    \label{lema.crlado}
Suppose that $\rho(\wh{f})$ is a polygon with rational endpoints, and let $S$ be a side of $\rho(\wh{f})$. Then, there exist $A\in\txt{SL}(2,\Z)$, $m,n\in\Z$ and $p\in\N$ such that the map $G=  T_1^{m}T_2^n A \wh{f}^p A^{-1}$ satisfies:
\begin{itemize}
\item $\rho(G) \subset \{ (x,y)\, : \, x\leq 0\}$,
\item $\rho(G)\cap\{(x,y)\, : \, x=0\} = D$, where $D$ is the side of the polygon $\rho(G)$ given by $D=T_1^mT_2^n A(pS)$.
\end{itemize}
\end{lema}

\noindent \textbf{The rotation set and periodic orbits.} A rational point $(p_1/q,p_2/q)\in\rho(\wh{f})$ (with gcd$(p_1,p_2,q)=1$) is \textit{realized by a periodic orbit} if there exists a periodic point for $f$ with rotation vector $(p_1/q,p_2/q)$, or equivalently, if there is $x\in\R^2$ such that 
$$\wh{f}^q(x)=x+(p_1,p_2).$$
We have the following realization results

\begin{teorema}[\cite{f2}]    \label{teo.frint}
Any rational point in the interior of $\rho(\wh{f})$ is realized by a periodic orbit.
\end{teorema}

\begin{teorema}[\cite{mz}]     \label{mz.cc}
Every extremal point of $\rho(\wh{f})$ is the rotation vector of some point.
\end{teorema}

\noindent \textbf{The rotation set and invariant measures.} For a compact $f$-invariant set $\Lambda\subset\T^2$, we denote by $\cl{M}_{f}(\Lambda)$ the family of $f$-invariant probability measures with support in $\Lambda$, and $\cl{M}_{f}=\cl{M}_{f}(\T^2)$. Define the \textit{displacement function} $\varphi:\T^2\ra\R^2$ by
$$\varphi(x)= \wh{f}(\wh{x})-\wh{x}, \ \ \ \txt{for $\wh{x}\in\pi^{-1}(x)$}.$$
This is well defined, as any two preimages of $x$ by the projection $\pi:\R^2\ra\T^2$ differ by an element of $\Z^2$, and $\wh{f}$ is $\Z^2$-periodic. Now, for $\mu\in\cl{M}_{f}$ we define the \textit{rotation vector} of $\mu$ as
$$\rho(\mu,\wh{f})=\int \varphi \, d\mu.$$
Then, we define the sets
$$\rho_{mes}(\Lambda,\wh{f})=\left\{ \rho(\mu,\wh{f})\, : \, \mu\in \cl{M}_{f}(\Lambda)\right\},$$
and
$$\rho_{erg}(\Lambda,\wh{f})=\left\{ \rho(\mu) \, : \, \mu \txt{ is ergodic for $f$ and $\txt{supp}(\mu)\subset\Lambda$}  \right\}.$$
When $\Lambda=\T^2$ we simply write $\rho_{mes}(\wh{f})$ and $\rho_{erg}(\wh{f})$.

\begin{proposicion}[\cite{mz}]    \label{prop.mz}
It holds the following:
$$\rho(\wh{f})=\rho_{mes}(\wh{f})=\txt{conv} (\rho_{erg}(\wh{f})).$$
\end{proposicion}

\subsection{Atkinson's Lemma.}    \label{sec.atk}

Let $(X,\mu)$ be a probability space, $T:X\ra X$ be an ergodic transformation with respect to $\mu$, and $\phi:X\ra\R$ a measurable map. We say that the pair $(T,\phi)$ is \textit{recurrent} if for any measurable set $A\subset X$ of positive measure, and every $\epsilon>0$ there is $n>0$ such that 
$$ \mu \left( A\cap T^{-n}(A) \cap \left\{ x\, : \, \sum_{i=0}^{n-1} |\phi(T^i(x))| < \epsilon \right\} \right) >0. $$

In \cite{atk} it is proved the following theorem. 

\begin{teorema}    
Let $(X,\mu)$ be a non-atomic probability space, $T:X\ra X$ an ergodic automorphism, and $\phi:X\ra \R$ an integrable function. Then, the pair $(T,\phi)$ is recurrent if and only if $\int \phi d\mu =0$.
\end{teorema}  

From this theorem, it is not difficult to obtain the following corollary, usually known as `Atkinson's Lemma'.


\begin{corolario}      \label{atk}
Let $X$ be a separable metric space and $\mu$ a probability measure in $X$ which is ergodic with respect to a measurable transformation $T:X\ra X$. Let $\phi:X\ra\R$ be an integrable function, with $\int \phi d\mu=0$. Then, there exists a full $\mu$-measure set $\wt{X}\subset X$ such that for any $x\in \wt{X}$, there is a sequence of positive integers $n_i$ with 
$$ T^{n_i}(x) \ra x      \    \txt{ and }     \    \left| \sum_{j=0}^{n_i-1} \phi(T^j(x)) \right| \ra 0 \ \ \ \ \ \txt{ as }  i\ra \infty.$$
\end{corolario}

\subsection{A theorem of Handel}

Denote by $\A$ the compact annulus $\A= [0,1]\times\T^1$, and by $\wt{\A}$ its universal cover $\wt{\A} = [0,1]\times\R$. Let $h:\A\ra\A$ be an orientation preserving, boundary component preserving homeomorphism, and let $\wt{h}:\wt{\A}\ra\wt{\A}$ be a lift of $h$. 

For a point $x\in\A$, and a lift $\tl{x}\in\wt{\A}$ of $x$, the limit
$$ \rho(x,\wt{h}) = \lim_{n\ra\infty} \frac{\txt{pr}_2( \wt{h}^n(\wt{x}) - \wt{x}) }{n},$$
whenever it exists is independent of the lift $\tl{x}$, and it is called the \textit{rotation number} of $x$. The \textit{pointwise rotation set} of $\wt{h}$ is defined as 
$$\rho_{point}(\wt{h})=\bigcup \rho(x,\wt{h}),$$
where the union is taken over all the $x$ in the domain of $\rho(\cdot,\wt{h})$.

The following theorem is part of Theorem 0.1 in \cite{han}.

\begin{teorema}   \label{teo.handel}
If $h:[0,1]\times \T^1 \ra [0,1]\times \T^1$ is an orientation preserving, boundary component preserving homeomorphism and $\wt{f}:[0,1]\times \R\ra [0,1]\times\R$ is any lift, then:
\begin{enumerate}
\item $\rho_{point}(\wt{f})$ is a closed set.
\item With the exception of at most a discrete set of values $r$ in $\rho_{point}(\wt{h})$, there is a compact invariant set $Q_r$ such that $\rho(x,\wt{h})=r$ for all $x\in Q_r$. If $r$ is rational, then $Q_r$ exists and is realized by a periodic orbit. 
\end{enumerate} 
\end{teorema}

\subsection{Stable and unstable sets, in the case there is a Brouwer curve.}        \label{sec.conjest}

Let $\wh{f}$ be any lift of a homeomorphism of $\T^2$ homotopic to the identity. Suppose there is a Brouwer $(0,1)$-curve $\ell$ for $\wh{f}$. For $i\in\N_0$, denote $\ell_i = T_1^i(\ell)$. 

\begin{figure}[h] 
\begin{center} 
\includegraphics{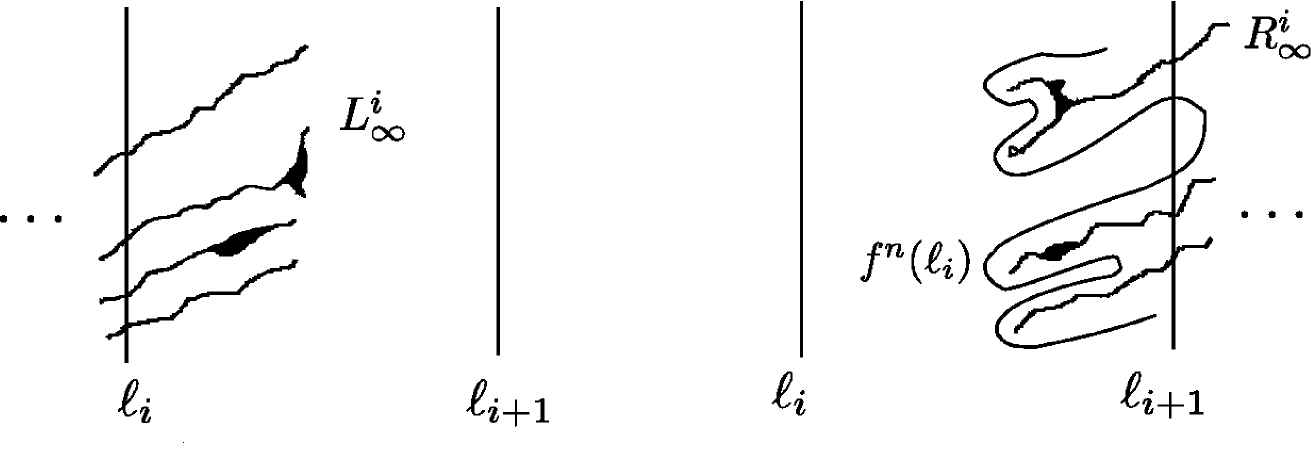}
\caption{Some examples of the sets $L_{\infty}^i$ and $R_{\infty}^i$.}
\label{fig.xi}  
\end{center}  
\end{figure}

Following \cite{dav}, for each $i\in\N$, we define the sets $L_{\infty}^i$ and $R_{\infty}^i$, which in some sense are the `stable' and `unstable' sets, respectively, of the maximal invariant set in $[\ell_i,\ell_{i+1}]$ for $\wh{f}$. Let
$$R_{\infty}^i= \bigcap_{n\in\Z} R \left(\wh{f}^n(\ell_i)\right), \ \ \ \txt{and} \ \ \    L_{\infty}^i=\bigcap_{n\in\Z} L (\wh{f}^{-n}(\ell_{i+1}))$$
(see Fig. \ref{fig.xi}.) 

By definition, the sets $R_{\infty}^i,$ and $L_{\infty}^i$ are $\wh{f}$-invariant. As $\ell$ is a Brouwer curve for $\wh{f}$, $\wh{f}(\ell_i) \subset R(\ell_i)$ for all $i$. Therefore,
$$R_{\infty}^i= \{x\in \R^2\ : \  \wh{f}^{-n}(x)\in R(\ell_i) \ \forall n\geq 0\},$$
and
$$L_{\infty}^i= \{x\in \R^2\ : \  \wh{f}^{n}(x)\in L(\ell_{i+1}) \ \forall n\geq 0\},$$
Therefore, for each $i$, the set $R_{\infty}^i\cap L_{\infty}^i$ is the maximal invariant set of $[\ell_i,\ell_{i+1}]$ for $\wh{f}$. Observe that, either the set $R_{\infty}^i\cap L_{\infty}^i$ is non-empty for all $i$, or $\wh{f}^n(\ell_0)\subset R(\ell_1)=R(T_1(\ell_0))$ for some $n$.

We have
$$R_{\infty}^i\cap[\ell_i,\ell_{i+1}]= \{x\in[\ell_i,\ell_{i+1}]\ : \ d(\wh{f}^{-n}(x),L_{\infty}^i\cap R_{\infty}^i)\ra 0 \txt{ as } n\ra\infty\},$$
and
$$L_{\infty}^i\cap[\ell_i,\ell_{i+1}]=\{x\in [\ell_i,\ell_{i+1}] \ : \ d(\wh{f}^{n}(x),L_{\infty}^i\cap R_{\infty}^i)\ra 0 \txt{ as } n\ra\infty\}.$$
That is, the set $L_{\infty}^i\cap[\ell_i,\ell_{i+1}]$ can be thought as the `local stable set' of $R_{\infty}^i\cap L_{\infty}^i$, and $R_{\infty}^i\cap[\ell_i,\ell_{i+1}]$ can be thought as the `local unstable set' of $R_{\infty}^i\cap L_{\infty}^i$.

The following lemma is proved in \cite{dav} (Lemma 6.8 in that article).

\begin{lema}\label{discos}
For every $i\geq 0$:
\begin{enumerate}
\item if $C$ is a connected component of $R_{\infty}^i$, then $\sup \txt{pr}_1(C)=+\infty$, and
\item if $C'$ is a connected component of $L_{\infty}^i$, then $\inf \txt{pr}_1(C')=-\infty$.
\end{enumerate}
\end{lema}

The following lemma is a consequence of the periodicity of $\wh{f}$. For a proof, see Lemma 6.10 in \cite{dav}.

\begin{lema}  \label{lema.limit}
Let $i\in\N$ and suppose that there is $n>0$ such that $\wh{f}^n(\ell)\cap R(\ell_i)\neq\empt$. There is $K>0$  such that if $C\subset(\ell,\ell_i)$ is forward invariant, then $\txt{diam}_2(C) < K$.
\end{lema}

As a consequence of the last two lemmas we have the following. 

\begin{lema}     \label{lema.Lcomp}
Let $i\in\N$ and suppose that there is $n>0$ such that $\wh{f}^n(\ell)\cap R(\ell_i)\neq\empt$. If $C$ is a connected component of $L_{\infty}^i\cap R(\ell)$, we have:
\begin{itemize}
\item $\ol{C}$ is compact, 
\item $\rho(\ol{C},\wh{f})$ consists of a point,
\item $\ol{C}\cap\ell\neq\empt$.
\end{itemize}
\end{lema}
\begin{prueba}
To prove that $\ol{C}$ is compact, first note that $\ol{C}\subset L_{\infty}^1\cap \ol{R}(\ell)$, $\ol{C}$ is forward invariant, and $\wh{f}^n(\ol{C})\subset L(\ell_2)$ for all $n$. If $\ol{C}$ was not compact, it would be unbounded in the vertical direction, and this would contradict Lemma \ref{lema.limit}. 

If $\rho(C,\wh{f})$ consisted of more than one point, then its forward iterates would be contained in $(\ell,\ell_2)$ and would have arbitrarily large vertical diameter, which would also contradict Lemma \ref{lema.limit}.

Finally, the fact that $\ol{C}\cap\ell\neq\empt$ is a consequence of Lemma \ref{discos}.
\end{prueba}

\subsection{Brouwer Theory}

\subsubsection{Brouwer foliations.}

Let $\cl{I}=(f_t)_{t\in[0,1]}$ be an isotopy on a surface $S$ from $f_0=\txt{Id}$ to a homeomorphism $f_1=f$. We say that a point $x\in S$ is \textit{contractible} for $f$ if the isotopy loop $f_t(x)$ is homotopically trivial in $S$. 

Given a topological oriented foliation $\cl{F}$ of $S$, we say that the isotopy $\cl{I}$ is \textit{tranverse} to $\cl{F}$ if for each $x\in S$, the isotopy path $f_t(x)$ is homotopic with fixed endpoints to an arc which is positively transverse to $\cl{F}$ in the usual sense. 

A clear obstruction to the existence of a foliation transverse to $\cl{I}$ is the existence of contractible points for $f$. The following theorem from Le Calvez tells us that this is the only obstruction.  

\begin{teorema}[\cite{lc2}]     \label{teo.lc}
If there are no contractible fixed points for $f$, then there exists a foliation $\cl{F}$ without singularities which is transverse to $\cl{I}$. 
\end{teorema}

In the case that there exist contractible fixed points for $f$, the previous theorem is still useful. This is thanks to the following result from Jaulent. 

\begin{teorema}[\cite{jau}]    \label{teo.jau}
Given an isotopy $\cl{I} = (f_t)_{t\in[0,1]}$ on $S$ from the identity to a homeomorphism $f$, there exists a closed set $X\subset\txt{Fix}(f)$ and an isotopy $\cl{I}'=(f'_t)_{t\in[0,1]}$ on $S\minus X$ from $\txt{Id}_{S\minus X}$ to $f|_{S\minus X}$ such that
\begin{itemize}
\item for each $z\in S\minus X$ the arc $(f'_t(z))_{t\in[0,1]}$ is homotopic with fixed endpoints (in $S$) to $(f_t(z))_{t\in[0,1]}$,
\item there are no contractible fixed points for $f|_{S\minus X}$ with respect to $\cl{I}'$.
\end{itemize}
\end{teorema}

\begin{remark}    \label{remark.jau1}
By Theorem \ref{teo.lc} it follows that there exists a foliation $\cl{F}_X$ of $S\minus X$ without singularities which is transverse to $\cl{I}'$. 
\end{remark}

\begin{remark}    \label{remark.jau2}
If the set $X$ is totally disconnected, the isotopy $\cl{I}'$ may be extended to an isotopy on $S$ that fixes $X$; that is, $f'_t(x)=x$ for each $x\in X$ and all $t\in[0,1]$. Similarly, the foliation $\cl{F}_X$ may be extended to a foliation with singularities $\cl{F}$ on $S$, where the set of singularities of $\cl{F}$ coincides with $X$. Also, if $\wh{\cl{I}}'=(\wh{f}'_t)_{t\in[0,1]}$ and $\wh{\cl{I}}=(\wh{f}_t)_{t\in[0,1]}$ are the respective canonical lifts of $\cl{I}'$ and $\cl{I}$ to the universal cover of $S$, then $\wh{f}'_1=\wh{f}_1$. This follows from the fact that, if $z\in S\minus X$, the paths $(f'_t(z))_t$ and $(f_t(z))_t$ are homotopic with fixed points in $S$. 
\end{remark}

\subsubsection{Fixing the isotopy $(f_t)_t$, the foliation $\cl{F}$ transversal to $(f_t)_t$, the set $X\subset\txt{Fix}(f)$ and their lifts $(\wh{f}_t)_t$, $\wh{\cl{F}}$, $\wh{X}$ to $\R^2$.}      \label{sec.isot}

Let $f$ and $\wh{f}$ be as in Theorem D. From now on we fix an isotopy $(f_t)_t$ on $\T^2$ from $f_0=\txt{Id}$ to $f_1=f$ such that the canonical lift $(\wh{f}_t)_t$ of $(f_t)_t$ to $\R^2$ ends in $\wh{f}_1= \wh{f}$. By Theorems \ref{teo.lc} and \ref{teo.jau} and by Remark \ref{remark.jau1}, there exists a set $X\subset \txt{Fix}(f)$, an isotopy $(f_t)_t$ on $\T^2\minus X$ from Id to $f$ and a foliation $\cl{F}_X$ of $\T^2\minus X$ such that $(f_t)_t$ is transverse to $\cl{F}_X$. Let $\wh{\cl{F}}$ and $\wh{X}$ be lifts of $\cl{F}$ and $X$ to $\R^2$, respectively. 

Note that if $X$ is totally disconnected, then by Remark \ref{remark.jau2} we have that:
\begin{itemize}
\item $(f_t)_t$ extends to an isotopy on $\T^2$, denoted also by $(f_t)_t$, that fixes $X$,
\item there exists a foliation $\cl{F}$ of $\T^2$ with sing$(\cl{F})=X$ and such that $\cl{F}\minus\txt{sing}(\cl{F})$ is transverse to the isotopy $(f_t)_t$ restricted to $\T^2\minus X$.
\end{itemize}

\begin{remark}
Suppose that the set $X$ is totally disconnected. Consider a closed set $F\subset\wh{X}$ and the universal covering $\D\to \R^2\minus F$. Whenever we refer to the canonical lift of $\wh{f}$ to $\D$, it will be implicit that it is the canonical lift with respect to the isotopy $(f_t)_t$. 
\end{remark}

\subsubsection{Some lemmas.}

Through all of this section we suppose that the set $X$ from $\S$\ref{sec.isot} is totally disconnected. The following lemma is trivial.

\begin{lema}   \label{lemabrouwer1}
Suppose that a leaf $\ell$ of $\wh{\cl{F}}$ is a Brouwer curve for $\wh{f}$. Let $\beta$ be a compact arc such that $\beta(1)\in\gamma$. We have:
\begin{enumerate}
\item if $\beta(0)\in \wh{X}\cap L(\ell)$ then $\wh{f}^n(\beta) \cap \ell\neq\empt$ for all $n\in\N$. 
\item if $\beta(0)\in \wh{X}\cap R(\ell)$ then $\wh{f}^{-n}(\beta) \cap \ell\neq\empt$ for all $n\in\N$. 
\end{enumerate}
\end{lema}

Recall that if $\gamma\in\wh{\cl{F}}$, sing$(\gamma)$ denotes the set of singularities of $\wh{\cl{F}}$ contained in $\txt{Fill}(\omega(\gamma))\cup\txt{Fill}(\alpha(\gamma))$.

\begin{definicion}    \label{def.btl}
Let $\gamma\subset\R^2$ be a leaf of $\wh{\cl{F}}$, and let $\beta:[0,1]\to\R^2$ be such that $\beta(t_0)\in\gamma$, for some $t_0\in[0,1]$. Consider a lift $\ol{\gamma}$ of $\gamma$ to the universal cover of $\R^2\minus\txt{sing}(\gamma)$, and let $\ol{\beta}$ be the lift of $\beta$ such that $\ol{\beta}(t_0)\in\ol{\gamma}$.
\begin{enumerate}
\item We say that $\beta$ \textit{arrives in $\gamma$ in $t_0$ by the left} if $t_0>0$ and if there is $t_1\in(0,t_0)$ such that $\ol{\beta}(t)\in L(\ol{\gamma})$ for all $t\in(t_1,t_0)$. Similarily, we say that $\beta$ \textit{arrives in $\gamma$ in $t_0$ by the right} if $t_0>0$ and there is $t_1\in(0,t_0)$ such that $\ol{\beta}(t)\in R(\ol{\gamma})$ for all $t\in(t_1,t_0)$.
\item We say that $\beta$ \textit{leaves $\gamma$ in $t_0$ by the right} if $t_0<1$ and if $\beta^{-1}$ arrives in $\gamma$ in $t_0$ by the right. Also, we say that $\beta$ \textit{leaves $\gamma$ in $t_0$ by the left} if $t_0<1$ and $\beta^{-1}$ arrives in $\gamma$ in $t_0$ by the left. 
\end{enumerate}
\end{definicion}

Next lemma describes how an arc might get `anchored' to a leaf of $\wh{\cl{F}}$. For a leaf $\gamma\in\wh{\cl{F}}$ recall our notation $\txt{sing}(\gamma)= \txt{sing}(\wh{\cl{F}})\cap (\txt{Fill}(\alpha(\gamma)) \cup \txt{Fill}(\omega(\gamma)))$. 

\begin{figure}[h]        
\begin{center} 
\includegraphics{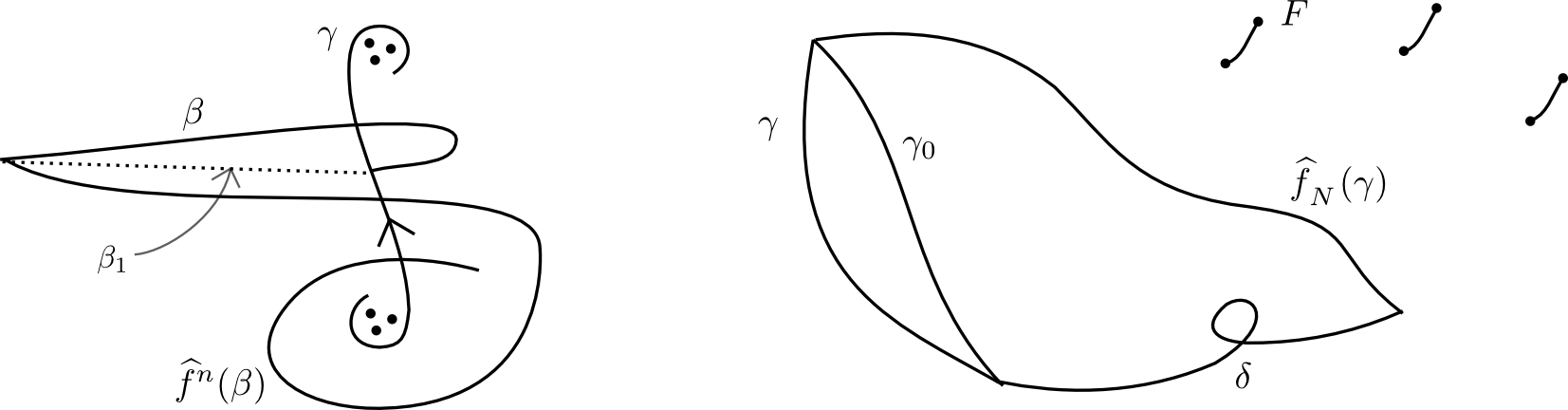}
\caption{Lemma \ref{lemabrouwer2} (left), and Lemma \ref{isot.lema2} (right).}
\label{fig.br}
\end{center}  
\end{figure}

\begin{lema}           \label{lemabrouwer2}
Let $\gamma$ be a non-compact leaf of $\wh{\cl{F}}$ with compact closure, and let $S = \txt{sing}(\gamma)$ (see Fig. \ref{fig.br}). Suppose that $\beta\subset \R^2\minus S$ is an arc such that $\beta(0)\in \wh{X}\minus S$, $\beta(1)\in\gamma$, and $\beta$ is homotopic wfe $\txt{Rel}(S)$ to an arc $\beta_1$ such that the only intersection of $\beta_1$ with $\gamma$ is the point $\beta_1(1)$.
\begin{enumerate}
\item If $\beta_1$ arrives in $\gamma$ in 1 by the left, then $\wh{f}^n \beta \cap \gamma\neq\empt$ for all $n\in\N$. 
\item If $\beta_1$ arrives in $\gamma$ in 1 by the right, then $\wh{f}^{-n} \beta \cap \gamma\neq\empt$ for all $n\in\N$. 
\end{enumerate}
\end{lema}

\begin{prueba}
We will prove item 1, item 2 being analogous. By the Poincar\'e-Bendixson Theorem we know that $S$ is non-empty. Consider the universal covering space of $\R^2\setminus S$, which is homeomorphic to the unit disc $\D$, and the canonical lift $\bar{f}:\D\ra\D$ of $\wh{f}|_{\R^2\setminus S}$. Any lift of $\gamma$ to $\D$ is a Brouwer curve for $\bar{f}$.

Consider lifts $\bar{\beta}$, $\bar{\beta}_1$ and $\bar{\gamma}$ of the arcs $\beta$, $\beta_1$ and the leaf $\gamma$, respectively, being these lifts such that $\bar{\beta}(1)=\bar{\beta}_1(1)\in\bar{\gamma}$. As $\beta$ is homotopic wfe Rel$(S)$ to $\beta_1$, we have that $\ol{\beta}(0)=\ol{\beta}_1(0)$. The arc $\bar{\beta}_1$ intersects $\bar{\gamma}$ only in $\bar{\beta}_1(1)$ and arrives in $\gamma$ in 1 by the left, and therefore $\bar{\beta}(0)=\bar{\beta}_1(0) \in L(\bar{\gamma})$. By Lemma \ref{lemabrouwer1} we have that $\bar{f}^n(\bar{\beta})$ intersects $\bar{\gamma}$ for all $n\in\N$, and therefore $\wh{f}^n \beta$ intersects $\gamma$ for all $n$.   
\end{prueba}

\begin{remark}
In last proof we see that, for all $n\in\N$, $\wh{f}^n \beta|_{[0,1)} \cap \neq\empt$. 
\end{remark}

For $t\in\R$, recall our notation 
$$\wh{f}_t := \wh{f}_{t\mod 1} \circ \wh{f}^{\left\lfloor t \right\rfloor}$$
(cf. Section \ref{sec.notac}), and note that $\wh{f}_n=\wh{f}^n$ for $n\in\Z$.

\begin{lema}        \label{isot.lema2}
Let $F\subset \R^2$ be a closed subset, and let $N\in\N$. Let $\gamma,\gamma_0:[0,N]\ra\R^2$ be arcs disjoint from $F$, and such that $\gamma$ and $\gamma_0$ are homotopic wfe $\txt{Rel}(F)$ (see Fig. \ref{fig.br}). Let $\delta$ be the arc $(\wh{f}_t(\gamma(1)))_{t\in[0,N]}$. Then, if $\{ \wh{f}_t(F) \, : \, t\in[0,N] \}$ is disjoint from $\gamma_0\cdot \delta$, the arc $\wh{f}_{N}(\gamma)$ and the curve $\gamma_0\cdot \delta$ are homotopic wfe $\txt{Rel}(\wh{f}^N(F))$.  
\end{lema}

\begin{prueba}
For $t\in[0,N]$ define $\delta_t:= \delta|_{[0,t]}$. For each $0\leq t \leq N$, the arc $\Gamma_t:=\wh{f}_t(\gamma)\cdot \delta_t^{-1} \cdot \gamma_0^{-1}$ is a (non necessarily simple) loop. As $\gamma$ and $\gamma_0$ are disjoint from $F$, and homotopic with fixed endpoints and $\txt{Rel}(F)$, we have that $F$ is contained in the unbounded component of $\R^2\setminus \Gamma_0$. By continuity, and as by hypothesis $\{ \wh{f}_t(F) \, : \, t\in[0,N] \}$ is disjoint from $\gamma_0\cdot \delta=\gamma_0\cdot\delta_N$, we have that for any $t$, $\wh{f}_t(F)$ is contained in the unbounded component of $\R^2\setminus \Gamma_t$, and in particular for $t=N$. This implies in turn that $\Gamma_N$ is homotopically trivial $\txt{Rel}(\wh{f}^N(F))$, or equivalently, the arcs $\wh{f}^N(\gamma)$ and $\gamma_0\cdot \delta_N=\gamma_0\cdot\delta$ are homotopic wfe $\txt{Rel}(\wh{f}^N(F))$.
\end{prueba}

\subsection{Poincar\'e-Bendixson theory.}

This section follows \cite{st}. Suppose that the set $X=\txt{sing}(\cl{F})$ is totally disconnected. We call a leaf of $\cl{F}$ a \textit{connection} if both its $\alpha$-limit and its $\omega$-limit are sets consisting of one element of sing$(\cl{F})$. Define a \textit{generalized cycle of connections} of $\cl{F}$ as a loop $\Gamma$ such that $\Gamma\minus\txt{sing}(\cl{F})$ is a (not necessarily finite) disjoint union of regular leaves of $\cl{F}$, with their orientation coinciding with that of $\Gamma$.

The following is a generalization of the Poincar\'e-Bendixson Theorem, and is a particular case of a theorem of Solntzev \cite{soln} (cf. also \cite{NemStep} $\S$1.78) and can be stated in terms of continuous flows due to a theorem of Guti\'errez \cite{Gut}.

\begin{teorema}     \label{poinc.ben}
Let $\phi = \{ \phi_t \}_{t\in\R}$ be a continuous flow on $\R^2$ with a totally disconnected set of singularities. If the forward orbit of a point $\{\phi_t(z) \}_{t\geq 0}$ is bounded, then its $\omega$-limit is one of the following:
\begin{itemize}
\item a singularity;
\item a closed orbit;
\item a generalized cycle of connections.
\end{itemize}
\end{teorema}
As $\cl{F}$ is an oriented foliation with singularities, it can be embeded in a flow (cf. \cite{whi33,whi41}) and we may apply last theorem to $\cl{F}$.

\subsection{Intersection number of certain lines.}     \label{sec.in}

Let $\Gamma_1,\Gamma_2\subset\R^2$ be two lines. We say that $\Gamma_1$ and $\Gamma_2$ have \textit{disjoint ends} if there is $T>0$ such that 
$$\inf_{t>T} \, d\left( \, \{\Gamma_1(t),\Gamma_1(-t) \} , \{\Gamma_2(t), \Gamma_2(-t)  \}  \, \right)   >0.$$
If $\Gamma_1$ and $\Gamma_2$ have disjoint ends, then the algebraic intersection number of $\Gamma_1$ and $\Gamma_2$, denoted IN$(\Gamma_1,\Gamma_2)$ is well defined (although it is not invariant under arbitrary isotopies). We use the usual convention that IN$(\Gamma_1,\Gamma_2)=+1$ if $\Gamma_1$ traverses $\Gamma_2$ from the left to the right. We will make use of the following basic fact.

\begin{lema}    \label{lema.in}
Suppose that $\Gamma_1$ and $\Gamma_2$ are lines with disjoint ends, and $\txt{IN}(\Gamma_1,\Gamma_2)=1$. If $t_1, t_2$ are such that
$$ t_1 < \min \{ t\in\R \, : \, \Gamma_1(t) \in \Gamma_2 \},$$
$$ t_2 > \max \{ t\in\T \, : \, \Gamma_1(t) \in \Gamma_2 \},$$
then $\Gamma_1(t_1)\in L(\Gamma_2)$ and $\Gamma_1(t_2)\in R(\Gamma_2)$.
\end{lema}

\section{Theorems A, B and C.}               \label{sec.bcd}

\noindent \textbf{Theorem B.} This theorem is an immediate consequence of Theorem D if the vector $v\in\Q^2$ from the statement is vertical. The general case $v\in\Q^2\minus\{(0,0)\}$, is reduced to the case of vertical $v$ by standard arguments, using (\ref{eq.pq}) and Lemma \ref{lema.rsconj}.\\

\noindent \textbf{Theorem A.} Also, if $\rho(\wh{f})$ is a vertical segment of the form $\{0\}\times I$, Theorem A follows immediately from Theorem D, and the case of a general rational interval is reduced to this case by the use of (\ref{eq.pq}) and Lemma \ref{lema.rsconj}.\\  

\noindent \textbf{Theorem C.} This theorem is an easy consequence of Theorem B. Let $S$ be a side of the polygon $\rho(\wh{f})$. Recall that by Lemma \ref{lema.crlado}, there exist $A\in\txt{SL}(2,\Z)$, $m,n\in\Z$ and $p\in\N$ such that the map $G=T_1^mT_2^nA \wh{f}^p A^{-1}$ satisfies
\begin{itemize}
\item $\rho(G) \subset \{ (x,y)\, : \, x\leq 0\}$,
\item $\rho(G)\cap\{(x,y)\, : \, x=0\} = D$, where $D$ is the side of the polygon $\rho(G)$ given by $D=T_1^mT_2^n A(pS)$.
\end{itemize}

Let $w\in\Q^2$ be the unit vector orthogonal to $S$ pointing outwards $\rho(\wh{f})$. From the definition of the rotation set, (\ref{eq.pq}) and Lemma \ref{lema.rsconj} one may also show the following:
\begin{itemize}
\item there is $m_1>0$ such that 
\begin{equation}   \label{eq.pp}
\sup \{ \ \txt{pr}_1(G^n(x) - x)  \ : \ x\in\R^2, n\in\N \ \} < m_1
\end{equation}
if and only if there is $m_2>0$ such that 
$$\sup \left\{ \langle \wh{f}^n(x) - x - nr, w \rangle \ : \ x\in\R^2, r\in\rho(\wh{f}), n\in\N \right\} < m_2.$$
\end{itemize}

Observe that, as $\rho(G)\subset \{(x,y)\,:\, x\leq 0\}$, by Theorem B there is $m_1>0$ such that (\ref{eq.pp}) holds. In this way, we have that for any side $S_i$ of $\rho(\wh{f})$, if $w_i$ is the unit vector orthogonal to $S_i$ pointing outwards $\rho(\wh{f})$, there is $m_1^i>0$ such that
$$\sup \left\{ \langle \wh{f}^n(x) - x - nr, w_i \rangle \ : \ x\in\R^2, r\in\rho(\wh{f}), n\in\N \right\} < m_1^i.$$
Letting $m= 2\max \{m_1^i\}$, it follows that 
$$ \sup \{ \ d(\wh{f}^n(x)-x, n\rho(\wh{f})) \ : \ x\in\R^2, n\in\N\ \} < m,$$
as desired.

\section{Proof of Theorem D.}     \label{sec.proof}

The organization of the proof of Theorem D is as follows. Recall that the isotopy $(f_t)_t$ from Id to $f$ and the foliation $\cl{F}$ with set of singularities $X\subset \txt{Fix}(f)$ were fixed in $\S$\ref{sec.isot}. In Section \ref{sec.totdisc} we study mainly under what conditions one may assume that sing$(\wh{\cl{F}})$ is inessential. In Section \ref{sec.bh} we show there exists an essential circloid $\cl{C}\subset\R\times\T^1$ which is formed by leaves and singularities of the foliation $\wt{\cl{F}}$, which is the lift of $\cl{F}$ to $\R\times\T^1$. Using such circloid, in Section \ref{sec.M} we construct the bound $M(f)>0$ which will be used to show that, if there are horizontal displacements larger than $M(f)$, then we may assume that sing$(\wh{\cl{F}})$ is totally disconnected (Section \ref{sec.col}) and there is linear horizontal speed (Sections \ref{sec.ss} and \ref{sec.cs}).

\subsection{About inessential sets and the set $\txt{sing}(\cl{F})$.}      \label{sec.totdisc}

We begin with the following. 

\begin{lema}   \label{lema.fijos1}
There exist essential vertical loops $c_1$ and $c_2$ in $\T^2$ which are positively transverse to $\cl{F}$, disjoint from $\txt{sing}(\cl{F})$, $c_1$ is oriented upwards and $c_2$ is oriented downwards.
\end{lema}
\begin{prueba}
By the hypotheses of Theorem D, $\rho(\wh{f})$ contains points $(0,a)$ and $(0,b)$, with $a < 0 < b$. If $\max \txt{pr}_1(\rho(\wh{f})) >0$, then by Theorem \ref{teo.frint} there is a periodic point $p\in\T^2$ such that pr$_1(\rho(p,\wh{f})) >0$, and then Theorem D holds. Therefore, we may assume there exist extremal points of $\rho(\wh{f})$ of the form $(0,e_1)$, $(0,e_2)$, with $e_1 < 0 < e_2$, and by Proposition \ref{prop.mz} there are ergodic measures $\mu_1$, $\mu_2$ for $f$, with $\rho(\mu_1, \wh{f})=(0,e_1)$ and $\rho(\mu_2, \wh{f}) = (0,e_2)$. 

Consider the displacement function $\varphi:\T^2\ra\R^2$, $\varphi(x)= \wh{f}(\wh{x}) - \wh{x}$, where $\wh{x}\in \pi^{-1}(x)$ and $\pi:\R^2\ra\T^2$ denotes the canonical projection. Let $\varphi_1 = \txt{pr}_1 \circ\varphi$. As $\int \varphi_1 d\mu_2 = 0$, by Atkinson's Lemma \ref{atk} there exists a $\mu_2$-total measure set $X_1\subset\T^2$ such that, for any $x\in X_1$, there is a sequence of integers $n_i$ such that 
\begin{equation}    \label{eq.65a}
 f^{n_i}(x) \ra x \ \ \ \ \ \txt{and} \ \ \ \ \ \sum_{j=0}^{n_i-1} \varphi_1(f^{n_i}(x)) \ra 0 \ \ \ \ \ \txt{as} \ i\ra\infty.
\end{equation}
By Birkhoff's Theorem there exists a $\mu_2$-total measure set $X_2\subset\T^2$ such that, for any $x\in X_2$, 
\begin{equation}    \label{eq.66a}
\frac{1}{n_i}\sum_{j=0}^{n_i-1} \varphi(f^{n_i}(x))  = \frac{ \wh{f}^{n_i}(\wh{x}) - \wh{x} }{ n_i} \longrightarrow  \int\varphi d\mu_2 = (0,e_2) \ \ \ \ \txt{as} \ i\ra\infty,
\end{equation}
where $\wh{x}\in\pi^{-1}(x)$.

Let $x\in X_1\cap X_2$ and $\wh{x}\in\pi^{-1}(x)$. Consider a flow box $B$ for the foliation $\cl{F}$ containing $x$, and let $\wh{B}$ be a lift of $B$ containing $\wh{x}$. By (\ref{eq.65a}) and (\ref{eq.66a}) there exist integers $i,m\in\N$ such that $\wh{f}^{n_i}(\wh{x}) \in T_2^m(B)$. As the isotopy $\wh{f}_t$ is transverse to $\wh{\cl{F}}$, there exists an arc $\wh{\gamma}$ going from $\wh{x}$ to $\wh{f}^{n_i}(\wh{x})$, positively transverse to $\wh{\cl{F}}$ (and also homotopic wfe Rel$(\txt{sing}(\wh{\cl{F}}))$ to the isotopy arc $(\wh{f}_t)_{t\in[0,n_i]}$). 

Let $\gamma=\pi(\wh{\gamma})\subset\T^2$. The point $f^{n_i}(x)\in B$ might be joined to a point $p\in \gamma\cap B$ by an arc $\gamma'\subset B$ positively transverse to $\cl{F}$ (see Fig. \ref{fig.box}). If $\gamma_1$ is the subarc of $\gamma$ from $p$ to $f^{n_i}(x)$, then the loop $c_1= \gamma_1\cdot \gamma'$ is vertical, positively transverse to $\cl{F}$, disjoint from $\txt{sing}(\cl{F})$ and oriented upwards.

\begin{figure}[h]        
\begin{center} 
\includegraphics{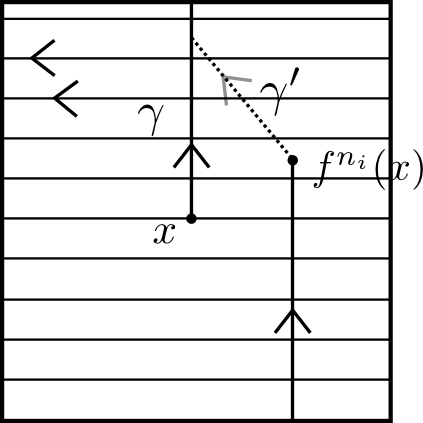}
\caption{}
\label{fig.box}
\end{center}  
\end{figure}

In a symmetric way, using the ergodic measure $\mu_1$, one constructs a vertical loop $c_2\subset\T^2$ which is oriented downwards, positively transverse to $\cl{F}$ and disjoint from $\txt{sing}(\cl{F})$.
\end{prueba}

\begin{lema}   \label{claim.td}
Let $c\subset\T^2\minus\txt{sing}(\cl{F})$ be any essential vertical loop (which exists by Lemma \ref{lema.fijos1}). Then, if 
$$\sup \{ \,  \txt{pr}_1 (\wh{f}^n(x) - x) \ : \  x\in\R^2, n > 0 \, \} > \txt{diam}_1(c) + 2,$$
the set $\txt{sing}(\cl{F})$ is inessential.
\end{lema}
\begin{prueba}
Let $U$ be the connected component of $\T^2\setminus\txt{sing}(\cl{F})$ that contains $c$. Consider a connected component $\wt{U}\subset \R^2$ of the preimage of $U$ by the projection map $\pi:\R^2\ra \T^2$. The set $\txt{sing}(\cl{F})$ is inessential if $U$ is fully essential, and to show that $U$ is fully essential, it suffices to see that $\wt{U}$ intersects $T_1(\wt{U})$.

Consider a lift $\tl{c}$ of $c$ contained in $\tl{U}$. Suppose there is $x\in \R^2$ and $n_0>0$ such that pr$_1(\wh{f}^{n_0}(x)-x) > \txt{diam}_1(c) + 2$. Consider an integer translate of $x$ contained in $(T_1^{-1}(\tl{c}) , \tl{c})$, and denote it also by $x$. We then have that $\wh{f}^{n_0}(x)\in R( T_1(\tl{c}))$. As $T_1(\tl{c}))\subset T_1(\wt{U})$, we then have that the isotopy path $\{\wh{f}_t(x)\}_{t\in[0,n_0]}$ is contained in $\wt{U} \cap T_1(\wt{U})$, and then $\wt{U}\cap T_1(\wt{U})\neq\empt$, as desired.  
\end{prueba}

The following proposition is a consequence of an improved version of a theorem of Moore (\cite{daver}, Theorems 13.4 and
25.1). For a proof, see \cite{st} (Proposition 1.6).

\begin{proposicion}     \label{prop.ines2}
Let $K\subset\T^2$ be a compact inessential filled set such that $f(K)=K$. Then there is a continuous surjection $h:\T^2\ra\T^2$ and a homeomoprhism $f':\T^2\ra\T^2$ such that 
\begin{itemize}
\item $h$ is homotopic to the identity,
\item $hf = f' h$,
\item $K'=h(K)$ is totally disconnected, and
\item $h|_{\T^2\minus K}: \T^2\minus K \ra \T^2\minus K'$ is a homeomorphism. 
\end{itemize} 
\end{proposicion}

By this proposition we have that if $K= \txt{Fill}(\txt{sing}(\cl{F}))$ is inessential, the homotopy $f'_t := h\circ f_t$ on $\T^2$ from $\txt{Id}|_{\T^2}$ to $f':= f'_1$ and the foliation $\cl{F}' = h(\cl{F})$ of $\T^2$ with singularities are such that $(f'_t)_t$ is transverse to $\cl{F}'$, and sing$(\cl{F}')$ is totally disconnected. As the map $h$ is homotopic to the identity, $\| \wh{h} - \txt{Id}\| <\infty$ for all lift $\wh{h}:\R^2\to\R^2$. By this and as $hf=f'h$ it holds the following.

\begin{proposicion}   \label{prop.ines}
 Let $\wh{f}'$ be the canonical lift of $f'$ with respect to the isotopy $(f'_t)_t$. Then:
\begin{enumerate}
\item $\rho(x,\wh{f})=v$ if and only if $\rho(h(x), \wh{f}')=v$,
\item Let $x\in\R^2$. There is $N>0$ such that $\txt{pr}_1( \wh{f}^{nN}(x) - x) > n$ for all $n>0$ if and only if there is $N'>0$ such that $\txt{pr}_1( \wh{f}'^{nN'}(x) - x) > n$ for all $n>0$. 
\end{enumerate}
\end{proposicion}

\subsection{The circloid $\cl{C}\subset\R\times\T^1$.}       \label{sec.bh}

Consider the lift $\wt{\cl{F}}$ of $\cl{F}$ to $\R\times\T^1$. We will show that there exists an essential circloid $\cl{C}\subset\R\times\T^1$ which is the union of leaves and singularities of $\wt{\cl{F}}$. We begin by showing that the leaves of $\wt{\cl{F}}$ are uniformly bounded horizontally.

\begin{claim}     \label{limit.horiz}
There is $m>0$ such that for every leaf $l\in\wh{\cl{F}}$, we have $\txt{diam}(\txt{pr}_1(l)) < m$.
\end{claim}
\begin{prueba}
By Lemma \ref{lema.fijos1}, there exist two closed vertcal loops $c_1$,$c_2\subset \T^2$, such that they are positively transverse to $\cl{F}$, $c_1$ is oriented upwards, and $c_2$ is oriented downwards. Consider lifts $\wh{c}_1,\wh{c}_2\subset\R^2$ of $c_1$ and $c_2$, respectively, such that $\wh{c}_1 \subset L(\wh{c}_2)$, and such that the open strip $(\wh{c}_1,\wh{c}_2)\subset\R^2$ contains a fundamental domain $D$ of the projection $\R^2\ra\T^2$. Finally, consider lifts $\wh{c}_0$ and $\wh{c}_3$ of $c_2$ and $c_1$, respectively, such that $\wh{c}_0 \subset L(\wh{c}_1)$ and $\wh{c}_2 \subset L(\wh{c}_3)$.

We will show that every leaf $l\in\wh{\cl{F}}$ that intersects the fundamental domain $D$ is contained in the open strip $(\wh{c}_0,\wh{c}_3)$. Setting $m=\txt{diam}_1((\wh{c}_0,\wh{c}_3))$, this will prove the claim. Suppose then that a leaf $l\in\wh{\cl{F}}$ intersects $D$ in a point $x$. Parametrize $l:\R\to\R^2$ according to its orientation, and let $l_1= l|_{[0,\infty)}$ and $l_2= l|_{(-\infty,0)}$. As $c_1$ and $c_2$ are positively transverse to $\cl{F}$, the leaves $\wh{c}_1$ and $\wh{c}_2$ are positively transverse to $\wh{\cl{F}}$, and also $\hat{c}_1$ is oriented upwards and $\hat{c}_2$ is oriented downwards. Therefore, the curve $l_1\subset l\in\cl{F}$ is contained in $(\wh{c}_1,\wh{c}_2)$. Analogously, as $\hat{c}_0$ is oriented downwards and $\hat{c}_3$ is oriented upwards, the curve $l_2$ is contained in $(\wh{c}_0,\wh{c}_3)$, and therefore $l=l_1\cup l_2$ is contained in $(\wh{c}_0,\wh{c}_3)$. As we mentioned, this proves the claim.   
\end{prueba}

We now state a result from \cite{j1} concerning the construction of circloids. We say that a set $U\subset\R\times\T^1$ is an \textit{upper generating set} if $U$ is bounded to the left and $\ol{U}$ is essential. Its \textit{associated lower component} $\cl{L}(U)$ is the connected component of $\R\times\T^1\minus\ol{U}$ that is unbounded to the left. Similarily, we call a set $L\subset\R\times\T^1$ a \textit{lower generating set} if $L$ is bounded to the right and $\ol{L}$ is essential, and its \textit{associated upper component} $\cl{U}(L)$ is the connected component of $\R\times\T^1\minus \cl{U}$ that is unbounded to the right. Observe that if $U$ is upper (lower) generating, then $\cl{L}(U)$ is lower (upper) generating, and then the expressions $\cl{U}\cl{L}(U)$, $\cl{L}\cl{U}(L)$, $\cl{LUL}(U)$, etc. make sense. 

\begin{lema}[Lemma 3.2 from \cite{j1}]    \label{lema.circloids}
Suppose that $U$ is an upper generating set. Then $\cl{C}^-(U):= \R\times\T^1\minus (\cl{UL}(U)\cup\cl{LUL}(U))$ is a circloid. Similarly, if $L$ is a lower generating set, then $\cl{C}^+:= \R\times\T^1\minus(\cl{LU}(L)\cup \cl{ULU}(L))$ is a circloid. 
\end{lema}

In the following claim we use this lemma to find a circloid $\cl{C}\subset\R\times\T^1$ formed by leaves and singularities of $\wt{\cl{F}}$.

\begin{claim}
There exists an essential circloid $\cl{C} \subset \R\times\T^1$ which is a union of leaves and singularities of $\wt{\cl{F}}$.
\end{claim}
\begin{prueba}
Let $C\subset\R\times\T^1$ be a vertical straight circle. Consider its saturation by the foliation $\wt{\cl{F}}$, i.e., the set $S(C)\subset\R\times\T^1$ which is the union of all the leaves and singularities of $\wt{\cl{F}}$ that intersect $C$. The set $\ol{S(C)}$ is then compact, because by Claim \ref{limit.horiz}, $S(C)$ is bounded. It is also essential, as it contains $C$, and therefore it is an upper generating set. The foliation $\wt{\cl{F}}$ may be embedded in a flow $\phi_t$ \cite{whi33,whi41}, and the set $S(C)$ is by definition totally invariant by $\phi_t$. Then, the sets $\cl{L}(S(C))$, $\cl{UL}(S(C))$, etc. and their complements are also totally invariant by $\phi_t$. 

By Lemma \ref{lema.circloids}, the set $\cl{C}:= \cl{C}^-(S(C))$ is a circloid, and by the above it is totally invariant by $\phi_t$. In other words, $\cl{C}$ is a union of leaves and singularities of $\wt{\cl{F}}$, as desired. 
\end{prueba}

We say that a leaf $\gamma \subset\R\times \T^1$ is \textit{bounded} if any lift of $\gamma$ to $\R^2$ is bounded.

\begin{claim} \label{claim.bounded}
If $\cl{C}$ contains a singularity of $\wt{\cl{F}}$, then the leaves of $\wt{\cl{F}}$ contained in $\cl{C}$ are bounded. 
\end{claim}
\begin{prueba}
We know by Claim \ref{limit.horiz} that such leaves are bounded horizontally. By contradiction, suppose on the contrary that there exists a leaf $\gamma$ of $\wt{\cl{F}}$ contained in $\cl{C}$ which is unbounded vertically. It is therefore easy to see that the omega-limit $\omega(\gamma)\subset\R\times\T^1$ is essential. As $\cl{C}$ is compact, $\omega(\gamma)$ is contained in $\cl{C}$. As $\cl{C}$ contains a singularity of $\wt{\cl{F}}$, $\omega(\gamma)\neq\gamma$ and $\omega(\gamma)$ is a proper essential subset of $\cl{C}$, which contradicts the fact that $\cl{C}$ is a circloid. This proves the claim. 
\end{prueba}

For further reference, we restate Claim \ref{claim.bounded} explicitly in terms of $\wh{\cl{F}}$.

\begin{claim}   \label{coro.bounded}
If $\cl{C}$ contains a singularity of $\wt{\cl{F}}$, the leaves of $\wh{\cl{F}}$ contained in $\pi^{-1}(\cl{C})\subset\R^2$ are bounded (where $\pi:\R^2\ra\R\times\T$ denotes the canonical projection).
\end{claim}

\subsection{Definition of the bound $M=M(f)$ from Theorem D.}     \label{sec.M}

Recall that by Lemma \ref{claim.td}, if $c\subset\T^2\minus\txt{sing}(\cl{F})$ is any essential vertical loop, and if
$$\sup \{ \txt{pr}_1 (\wh{f}^n(x) - x) \ : \  x\in\R^2, n\geq 0 \} > \txt{diam}_1(c) + 2,$$
we have that $\txt{sing}(\cl{F})$ is inessential. By a \textit{vertical} Brouwer line we mean a line $\ell\subset\R^2$ such that its (canonical) projection to $\R\times\T^1$ is an essential simple closed curve.

The definition of $M$ is done separately in two cases. \\
\\
\textit{Case 1: there exists a vertical Brouwer line $\ell$ for $\wh{f}$.} In this case define 
$$ M = \txt{diam}_1 (\ell) + 3.$$
What will be used of this definition is that it guarantees that, if there are $z$ and $n$ such that pr$_1(\wh{f}^n(z)-z)>M$, then $\txt{sing}(\cl{F})$ is inessential, and also 
\begin{equation}   \label{eq.111}
\wh{f}^n(\ell)\cap R(T_1^3(\ell))\neq\empt.
\end{equation}

\noindent \textit{Case 2: there is no vertical Brouwer line for $\wh{f}$.} In this case we may assume that if $(0,a)$,$(0,b)\in\rho(\wh{f})$ are as in Theorem D, then for any rational $r\in(a,b)$, the vector $(0,r)\in\rho(\wh{f})$ is realized by a periodic orbit. This is because of the following result, which is proved in \cite{dav} although it was already known in the folklore and is essentially contained in \cite{lc2}.

\begin{proposicion}    
If there is $r\in(a,b)$ such that $(0,r)\in\rho(\wh{f})$ is not realized by a periodic orbit, then there is a vertical Brouwer line for $\wh{f}$.
\end{proposicion}

As we are in the case that there is no vertical Brouwer line for $\wh{f}$, we may then assume there is $p\in\R^2$ be such that $\pi(p)\in\T^2$ is periodic for $f$ and such that 
\begin{equation}   \label{eq.a.3} 
\rho(\pi(p),\wh{f})=(0,r), \ \ r<0.
\end{equation}
Now, fix a connected component $\cl{C}_0\subset\R^2$ of $\pi^{-1}(\cl{C})$. Consider an integer translate of $p$, denoted also by $p$, such that there is a straight vertical line $l_1$ oriented upwards with $\cl{C}_0\subset L(l_1)$ and $\{\wh{f}_t(p) \, : \, t\in\R\}\subset R(l_1)$.

Let $n_1>0$ be such that there is a straight vertical line $l_2$ oriented upwards with 
\begin{equation}   \label{eq.ccc}
\{\wh{f}_t(p) \, : \, t\in\R\} \subset L(l_2) \ \ \ \txt{and} \ \ \ R(l_2) \supset T_1^{n_1}(\cl{C}_0).
\end{equation}

For $i\in\Z$, denote 
\begin{equation}   \label{eq.a.2}
\cl{C}_i= T_1^{n_1 i}(\cl{C}_0),
\end{equation}
and note that
\begin{equation}       \label{eq.4a}
\{\wh{f}_t(p) \, : \, t\in\R \} \subset (l_1,l_2) \subset R(\cl{C}_0)\cap L(\cl{C}_1).
\end{equation}
Fix $c\subset\T^2\minus\txt{sing}(\cl{F})$ an essential vertical loop (which exists by Lemma \ref{lema.fijos1}), and define
$$ M = \txt{diam}_1(\cl{C}) + \txt{diam}_1(c) + 2n_1.$$
What will be used of this definition is that it guarantees that if there are $z$ and $n$ such that pr$_1(\wh{f}^n(z)-z)>M$, then $\txt{sing}(\cl{F})$ is inessential and
\begin{equation}   \label{eq.a.1} 
\wh{f}^n(\cl{C}_i)\cap \cl{C}_{i+2} \neq\empt
\end{equation}
for any $i$.

\subsection{Collapsing the connected components of $\txt{sing}(\wh{\cl{F}})$.}      \label{sec.col}

The main part of the proof of Theorem D is contained in Sections \ref{sec.ss} and \ref{sec.cs}. To prove the theorem, in those sections it will be assumed that there is $x\in\R^2$ and $n>0$ such that pr$_1(\wh{f}^n(x)-x) > M$. In such case, by the construction of $M$ in last section the set sing$(\cl{F})$ is inessential. Therefore, by Proposition \ref{prop.ines2} and by the remark following it, we have that there is a surjection $h:\R^2\to\R^2$ such that the isotopy $f'_t=h\circ f_t$ is transverse to the foliation $\cl{F}' = h(\cl{F})$, and sing$(\cl{F}')$ is totally disconnected. 

Consider the canonical lift $\wh{f}'$ of $f'=f'_1$ with respect to the isotopy $(f'_t)_t$. Note that if $l\subset\T^2$ is an essential closed curve such that a lift $\ell$ of $l$ to $\R^2$ is a vertical Brouwer curve for $\wh{f}$, then $l\cap \txt{Fix}(f)=\empt$ (unless $l\subset\txt{Fix}(f)$, case in which $\wh{f}(\ell)$ is a translate of $\ell$ and Theorem D easily follows). Let $\wh{h}:\R^2\to\R^2$ be a lift of $h$. As $h|_{\T^2\minus\txt{Fix}(f)}$ is a homeomorphism, a curve $\ell$ is a vertical Brouwer curve for $\wh{f}$ iff $\wh{h}(\ell)$ is a vertical Brouwer curve for $\wh{f}'$. We conclude that if there are $x,n$ such that pr$_1(\wh{f}^n (x)-x) >M$, the properties of $M$ that will be used, namely (\ref{eq.111}) and (\ref{eq.a.1}), are still valid for $\wh{f}'$:
\begin{itemize}
\item if there is a vertical Brouwer curve $\ell$ for $\wh{f}$, then $\wh{f}'^n(\wh{h}(\ell))\cap T_1^3(\wh{h}(\ell))\neq\empt$,
\item if there is no vertical Brouwer curve for $\wh{f}$, then $\wh{f}'^n(\wh{h}(\cl{C}_i))\cap \wh{h}(\cl{C}_{i+2}) \neq\empt$.
\end{itemize}

By this and by Proposition \ref{prop.ines}, from now on we may make the following assumption.

\begin{asuncion}    \label{asunc.td}
The set $\txt{sing}(\cl{F})$ is totally disconnected. 
\end{asuncion}

\subsection{Case that $\cl{C}$ does not contain singularities.}    \label{sec.ss}


Observe that in this case, the lifts $\cl{C}_i$ from (\ref{eq.a.2}) are Brouwer curves for $\wh{f}$. Through this section, we denote
$$\ell = \cl{C}_0.$$

Consider the constant $M>0$ constructed in Section \ref{sec.M}. To prove Theorem D in this case, through all of this section we assume that there is $x\in\R^2$ and $N_1>0$ such that pr$_1(\wh{f}^{N_1}(x)-x) > M$ and we will show that there is $z\in\R^2$ and $N>0$ such that pr$_1(\wh{f}^{nN}(z)-z) > n$ for all $n\in\N$. 

We begin by giving an idea of the proof.

\subsubsection{Idea of the proof.}    \label{sec.ideass}

Recall that, by hypothesis, $\rho(\wh{f})$ contains two points $(0,a)$ and $(0,b)$ such that $a<0<b$. We will work under the following assumption, which will be justified in Section \ref{sec.per}. 

\begin{asuncion}     \label{asuncion.per}
Every rational point contained in the segment $\{0\}\times (a,b)$ is realized by a periodic orbit. 
\end{asuncion}

Under this assumption, we will construct a set $\cl{V}\subset R(\ell)$ such that it separates $R(\ell)$ and an iterate of $\cl{V}$ is contained in $R(T_1(\ell))$. Denote by $\cl{V}_{i,j}$ its integer translates $\cl{V}_{i,j}= T_1^i T_2^j(\cl{V})$, $i,j\in\Z$. We will define the property of a curve $\gamma$ having \textit{good intersection} with one of the sets $\cl{V}_{i,j}$. By now, we can think of that property as meaning that the arc $\gamma$ `traverses' $\cl{V}_{i,j}$ (cf. Fig. \ref{fig.V1}). 

We will prove that there is $N_2>0$ such that, if a curve $\gamma$ has good intersection with $\cl{V}_{i,j}$ for some $i,j$, then $\wh{f}^{N_2}(\gamma)$ has good intersection with $\cl{V}_{i+1,j'}$, for some $j'\in\Z$ (cf. Fig. \ref{fig.V2}). We will verify also that $\ell$ has good intersection with $\cl{V}$. An easy induction then will give us that, for any $n>0$ there is $j_n\in\Z$ such that $\wh{f}^{n N_2}(\ell)$ has good intersection with $\cl{V}_{n,j_n}$. In particular, 
$$ \wh{f}^{n N_2}(\ell) \cap R(T_1^n(\ell))  \supset \wh{f}^{n(N_2)}(\ell) \cap V_{n,j_n}   \neq\empt \ \ \ \ \ \forall\, n\in\N.$$
This will easily imply the desired conclusion, namely, there are $z\in\R^2, N>0$ such that pr$_1(\wh{f}^{nN}(z)-z) >n$ for all $n>0$.

\subsubsection{Assumption \ref{asuncion.per}.}    \label{sec.per}

In this section we show how Assumption \ref{asuncion.per} is justified by the results of \cite{dav}. We will see that if there is a rational point in the segment $\{0\} \times [a,b]\subset\rho(\wh{f})$ which is not realized by a periodic orbit, then the conclusions of Theorem D hold.

Suppose there is $p/q\in\Q\cap [a,b]$, with $p$ and $q$ coprime integers, and such that $(0,p/q)\in \rho(\wh{f})$ is not realized by a periodic orbit. We now show that we may assume that $(0,p/q)=(0,0)$. Let $\wh{g} = T_2^{-p}\wh{f}^q$, and note that $(0,0)\in \rho(\wh{g})$ is not realized by a periodic orbit of $f$ (see Section \ref{sec.rotset} for the basic properties of the rotation set). Note also the following:
\begin{itemize}
\item if there are $x\in\R^2$, $n\in\N$ s.t. pr$_1(\wh{g}^n(x)-x)> M$, then pr$_1(\wh{f}^{n'}(x)-x)> M$ with $n'=nq$;
\item if there are $x\in\R^2$, $N\in\N$ s.t. pr$_1(\wh{g}^{nN}(x)-x)> n$ for all $n\in\N$, then pr$_1(\wh{f}^{nN'}(x)-x)> n$ for all $n\in\N$, with $N'=qN$.
\end{itemize}
Therefore, without loss of generality we may assume that the rational point in $\rho(\wh{f})$ which is not realized by a periodic orbit is $(0,0)$.

In Proposition 5.2 of \cite{dav} it is shown\footnote{That proposition is stated under the hypothesis that $\rho(\wh{f})$ is a segment. However, that is used only guarantee the existence of $x,y\in\T^2$ with $\rho(x,\wh{f})=(0,a)$, $\rho(y,\wh{f})=(0,b)$ and $a<0<b$, which is within our hypotheses.} that the fact that $\rho(\wh{f})$ contains $(0,a),(0,b)$, $a< 0 < b$, implies that there exists a compact vertical leave $l \subset\T^2$ of the foliation $\cl{F}$ (this was already known in the folklore).

For $i\in\N$, fix a lift $\hat{l}_0\subset\R^2$ of the curve $l$ and let $\hat{l}_1=T_1(\hat{l}_0)$. Note that $\hat{l}_0$ and $\hat{l}_1$ are Brouwer curves for $\wh{f}$.

In \cite{dav} it is shown that, if there is $n>0$ such that $\wh{f}^n(\hat{l}_0)\cap R(\hat{l}_{1}) \neq\empt$, then $\max(\txt{pr}_1(\rho(\wh{f}))) >0$ (cf. Main Lemma 6.14 and Claim 6.15 in that article, where it is proved a stronger result). In this case, Frank's Theorem \ref{teo.frint} implies that there is a periodic point for $f$ which rotates rightwards, that is, there are $z\in\R^2$ and $N>0$ such that pr$_1(\wh{f}^{nN}(z)-z) >n$ for all $n\in\N$. 

On the other hand, by definition of $M=M(f)>0$ we have that if there are $x\in\R^2$ and $n>0$ such that pr$_1(\wh{f}^n(x)-x) >M$, then $\wh{f}^n(\cl{C}_i)\cap T_1(\cl{C}_{i})\neq\empt$ for all $i$ (cf. Section \ref{sec.M}), which in particular implies that $\wh{f}^n(\hat{l}_0)\cap\hat{l}_{1}\neq\empt$ (recall that $\cl{C}$ is also a compact vertical leaf of $\cl{F}$). 

We conclude that if there are $x\in\R^2$, $n>0$ such that pr$_1(\wh{f}^n(x)-x) >M$, then $\wh{f}^n(\hat{l}_0)\cap R(\hat{l}_1)\neq\empt$ and then there are $z\in\R^2$ and $N>0$ such that pr$_1(\wh{f}^{nN}(z)-z) >n$ for all $n\in\N$, in the case that $(0,0)\in\rho(\wh{f})$ is not realized by a periodic orbit. That is, Theorem D holds in the case that $\cl{C}$ does not contain singularities and $(0,0)\in\rho(\wh{f})$ is not realized by a periodic orbit, as we wanted.

\subsubsection{Statement of the main lemma and proof of Theorem D when $\cl{C}$ contains no singularities.}

The following main and technical lemma gives us the precise properties of the set $\cl{V}$ mentioned in $\S$\ref{sec.ideass}. 

Recall that $\ell=\cl{C}_0$ is a Brouwer curve for $\wh{f}$ and that $L_{\infty}^i$, $R_{\infty}^i$ denote the stable and unstable sets, respectively, of the maximal invariant set in each strip $(\ell_i,\ell_{i+1}) = (T_1^i(\ell), T_1^{i+1}(\ell))$ (cf. Section \ref{sec.conjest}).

\begin{figure}[h!]
\begin{center}    
\includegraphics{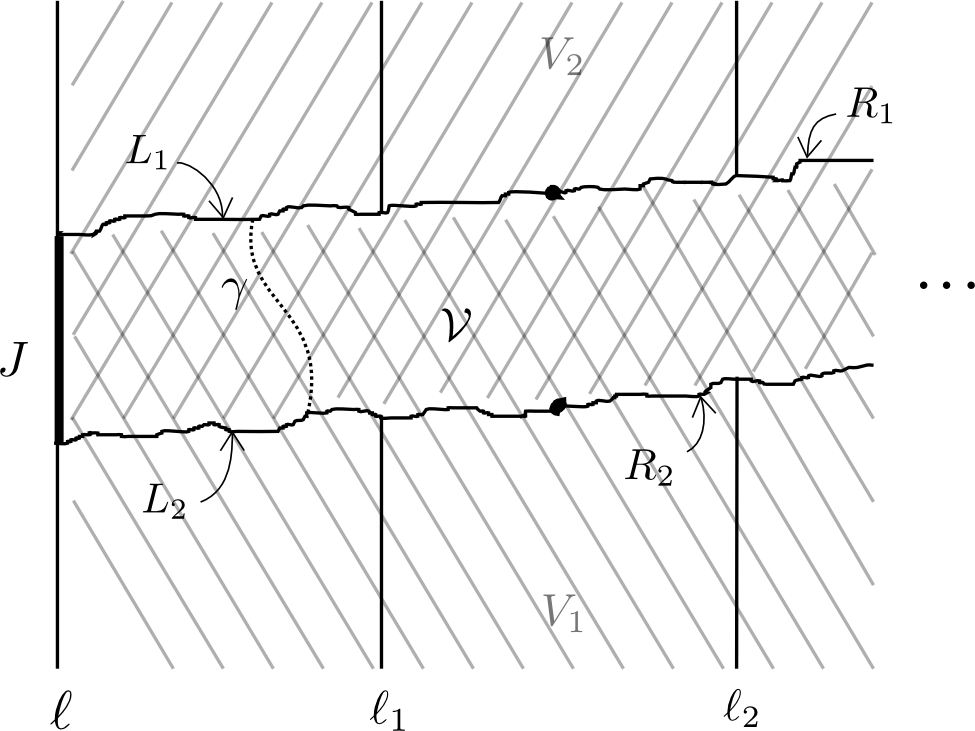}
\end{center}
\caption{The sets $L_i$, $R_i$, $\mathcal{V}$, $J$ from Lemma \ref{conjV1}, and an illustration of an arc $\gamma$ which has good intersection with $\mathcal{V}$.}
\label{fig.V1}
\end{figure}

\begin{lema}      \label{conjV1}
The following hold (cf. Fig. \ref{fig.V1}):
\begin{enumerate}
\item There exist sets $L_1, L_2, R_1,R_2 \subset R(\ell)$ such that:
		\begin{enumerate}
		\item for $i=1,2$, $L_i\subset L_{\infty}^1 \minus (L_{\infty}^0 \cup R^1_{\infty})$, and $\wh{f}^{n}(L_i)\subset (\ell_1,\ell_2)$ for all $n>0$ sufficiently large, 			\label{it.V4}
		\item $R_i\subset R_{\infty}^1$, and then $\wh{f}^{-n}(R_i)\subset R(\ell_1)$ for all $n\geq 0$,               \label{it.V5}
		\item $\rho(\ol{L}_i,\wh{f})=\{r_i\}$ for some $r_1\neq r_2$, and in particular $\ol{L}_1\cap\ol{L}_2=\empt$,                     \label{it.V8}
		\item the sets $L_i$ and $R_i$ are connected, the $\ol{L}_i$ are compact, the $R_i$ are unbounded to the right, and the sets $F_i=L_i\cup R_i$ separate $R(\ell)$,    \label{it.V2} 
		\end{enumerate}
\item There is a set $\cl{V}\subset R(\ell)$ such that:
		\begin{enumerate}
		\item $\cl{V}$ is of the form $\cl{V} = ( V_1 \cap V_2 ) \cup R_1$, for some open sets $V_1,V_2\subset R(\ell)$ such that $L_1\subset V_2$, $L_2\subset V_1$, and $		 \pr V_i \subset 	 \ell\cup F_i$,    \label{it.V10}
		\item $\pr \cl{V} \subset F_1\cup F_2 \cup \ell$, $\empt\neq V_1\cup V_2\subset \txt{int}(\cl{V})$, and $\cl{V}$ separates $R(\ell)$,    \label{it.V1}
		\item there is $n_1>0$ such that $\wh{f}^{n_1}(\ol{\cl{V}})\subset R(\ell_1)$,                 \label{it.V7}
    \item $\ol{\cl{V}}\cap [\ell,\ell_1]$ is compact,                      \label{it.V9}
    \item the set $J:=\pr\cl{V}\cap\ell$ is a non-degenerate arc such that $J \cap L_{\infty}^0 = \empt$, $J(0)\in\ol{L}_2$, and $J(1)\in\ol{L}_1$.          \label{it.V6}
		\end{enumerate}
\end{enumerate} 
\end{lema}

The proof of this lemma will be given in Section \ref{sec.conjv1}.

\begin{notacion}
For $i\in\N_0$ and $j\in\Z$ we denote $\ell_i=T_1^i(\ell)$, and $\cl{V}_{i,j}=T_1^i T_2^j (\cl{V})$. Similarily, for $k=1,2$, we set $L_k^{i,j}=T_1^i T_2^j (L_k)$, $R_k^{i,j}= T_1^i T_2^j(R_k)$ and $F_k^{i,j}=T_1^i T_2^j(F_k)$.
\end{notacion}

The main induction step in this section is given by Proposition \ref{N2curvalibre} below. In order to state it we need the following definition.

\begin{definicion}       \label{def.goodint}
Let $i\in\N_0$, $j\in\Z$. We say that a curve $\gamma$ has \textit{good intersection} with $\cl{V}_{i,j}$ if the following hold (see Fig. \ref{fig.V1}):
\begin{itemize}
\item $\gamma\cap R_{\infty}^{i+1}=\empt$,
\item one endpoint of $\gamma$ lies in $\ol{L}_1^{i,j}$ and the other in $\ol{L}_2^{i,j}$,
\item $\txt{int}(\gamma) \subset \txt{int}(\cl{V}_{i,j}) \cup J$, where $J\subset\ell$ is the arc from Lemma \ref{conjV1}-\ref{it.V6}.
\end{itemize}
\end{definicion}

\begin{figure}[h!]
\begin{center}    
\includegraphics{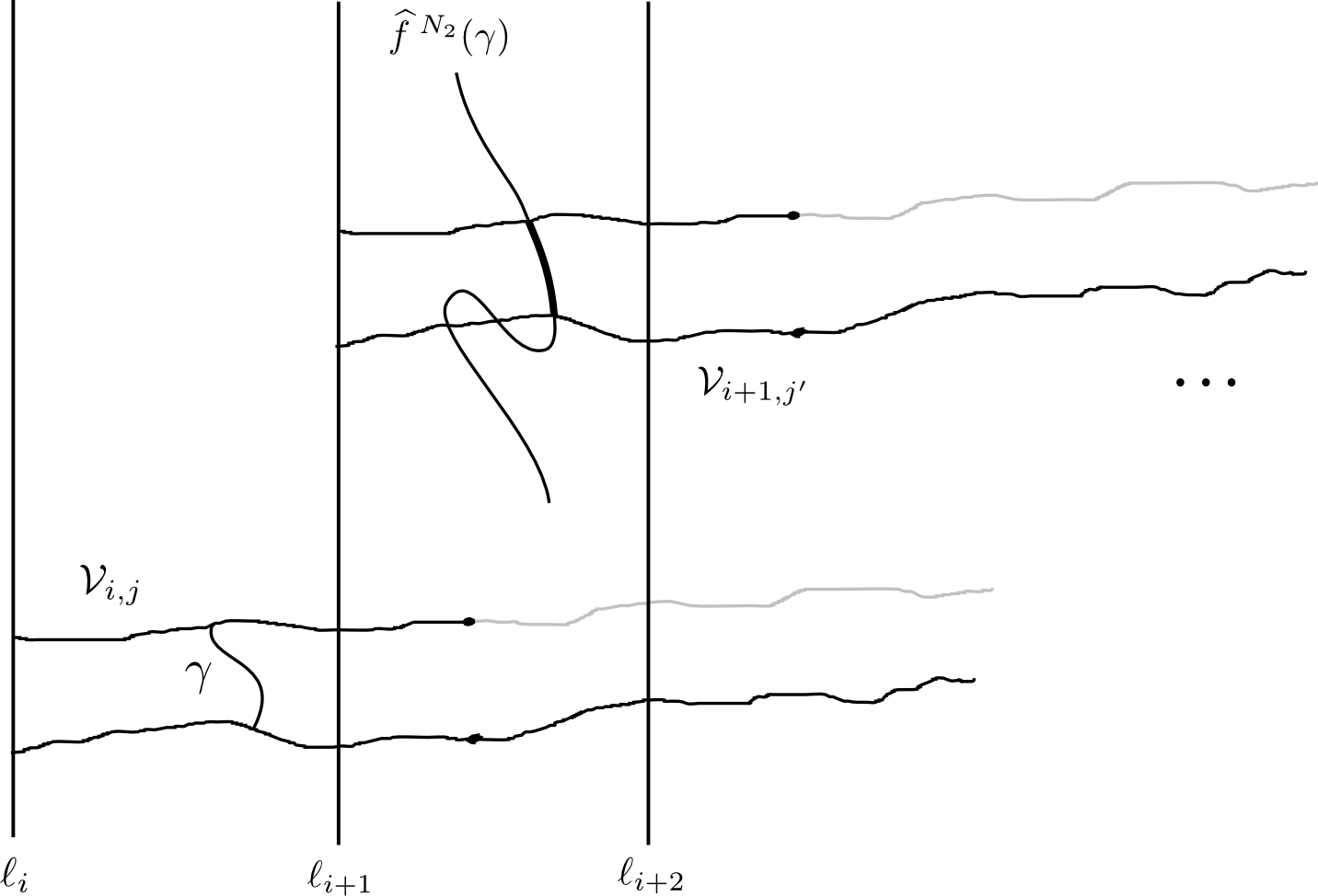}
\end{center}
\caption{Illustration of Proposition \ref{N2curvalibre}.}
\label{fig.V2}
\end{figure}

Note that, in particular, if an arc $\gamma$ has good intersection with $\cl{V}_{i,j}$ for some $i,j$, then $\gamma\subset \ol{R}(\ell_i)$.

\begin{proposicion}    \label{N2curvalibre}
There is $N_2>0$ such that:
\begin{enumerate}
\item if an arc $\gamma$ has good intersection with $\cl{V}_{i,j}$ for some $i,j\in\Z$, then $\wh{f}^{N_2}(\gamma)$ contains an arc which has good intersection with $\cl{V}_{i+1,j'}$, for some $j'\in\Z$ (see Fig. \ref{fig.V2}),    \label{it.qq}
\item the arc $J\subset\ell$ from Lemma \ref{conjV1}-\ref{it.V6} is such that $\wh{f}^{N_2}(J)$ has good intersection with $\cl{V}_{1,j}$, for some $j\in\Z$.             \label{paso1curvalibre}
\end{enumerate}
\end{proposicion}

Before proving this proposition, we use it to prove Theorem D.

\begin{prueba}[Proof of Theorem D when $\cl{C}$ contains no singularities.]
Let $n_0>0$ be such that $\min \{\txt{pr}_1(\ell_{n_0}))\} > \max \{ \txt{pr}_1(\ell) \} + 1$. By proposition \ref{N2curvalibre} we easily get by induction that $\wh{f}^{n_0N_2}(J)$ contains an arc $\beta_{1}$ which has good intersection with $\cl{V}_{n_0,j_1}$, for some $j_1\in\Z$. In particular, $\beta_{1}\subset \ol{R}(\ell_{n_0})$.

Let 
$$\ol{\beta}_{1}=\wh{f}^{-n_0 N_2}(\beta_{1})\subset J.$$ 
Analogously, by induction using item 1 of Proposition \ref{N2curvalibre} it easily follows that there is a sequence of arcs $\ol{\beta}_{n}\subset \ol{\beta}_{1}$ such that $\ol{\beta}_k\subset\ol{\beta}_l$ if $k>l$, and such that for all $n \geq 1$, the arc $\wh{f}^{n \cdot n_0N_2}(\ol{\beta}_{n})$ has good intersection with $\cl{V}_{n\cdot n_0,j_n}$ for some $j_n\in\Z$, and in particular, $\wh{f}^{n\cdot n_0 N_2}(\ol{\beta}_n)\subset \ol{R}(\ell_{n\cdot n_0})$.  

If $z\in\cap_{n\geq 0} \ol{\beta}_{n}$, then for every $n\geq 1$, $\wh{f}^{n \cdot n_0 N_2}(z)\in\ol{R}(\ell_{n \cdot n_0})$, and in particular, by the choice of $n_0$, 
$$\txt{pr}_1(\wh{f}^{n \cdot n_0N_2}(z)-z) > n \ \ \ \ \txt{ for all } n > 0.$$
Setting $N= n_0 N_2$, this finishes the proof.
\end{prueba}

We now prove Proposition \ref{N2curvalibre}. 

\begin{prueba}[Proof of Proposition \ref{N2curvalibre}.]
Without loss of generality, assume that $i=j=0$. Let $K= \txt{diam}_2(\cl{V} \cap (\ell,\ell_1))+1,$ which by item \ref{it.V9} from Lemma \ref{conjV1} is finite, and note that, by the periodicity of $\wh{f}$, 
$$K= \txt{diam}_2(\cl{V}_{i,j} \cap (\ell_i,\ell_{i+1}))+1 \ \ \ \forall\, i,j.$$

Consider the sets $L_1,L_2$ from Lemma \ref{conjV1} and let $n_1$ be as in item \ref{it.V7} of such lemma. By item \ref{it.V4} we have that there is $n_2>n_1$ such that
$$\wh{f}^n( L_1\cup L_2) \subset (\ell_1,\ell_2) \ \ \ \forall n\geq n_2.$$
By item \ref{it.V8} we have that $\rho(\ol{L}_i,\wh{f})=r_i$ for some $r_1\neq r_2$, and therefore there is $N_2>n_2$ such that
$$ d(\txt{pr}_2(\wh{f}^{N_2}(\ol{L}_1)), \txt{pr}_2(\wh{f}^{N_2}(\ol{L}_2)) ) >K.$$
As $\gamma$ has good intersection with $\cl{V}$, we have that an endpoint of $\wh{f}^{N_2}\gamma$, say $\wh{f}^{N_2}\gamma(0)$, belongs to $\wh{f}^{N_2}(L_1)$ and $\wh{f}^{N_2}\gamma(1)\in\wh{f}^{N_2}(L_2)$, and by definition of $K$ we have then that, for some $j_1\in\Z$, the set pr$_2(\cl{V}_{1,j_1})$ is contained in the interval between pr$_2(\wh{f}^{N_2}\gamma(0))$ and pr$_2(\wh{f}^{N_2}\gamma(1))$. As the sets $R_i$ are contained in $R_{\infty}^1\subset R(\ell_2)$ (item \ref{it.V5}) and as the sets $F_i^{1,j_1} = L_i^{1,j_1}\cup R_i^{1,j_1}$ separate $R(\ell_1)$ (by item \ref{it.V1} and by the periodicity of $\wh{f})$, we obtain that, for $i=1,2$, the points $\wh{f}^{N_2}\gamma(0)$ and $\wh{f}^{N_2}\gamma(1)$ belong to different connected components of $R(\ell_1)\minus F_i^{1,j_1}$.

Also, as $\gamma$ has good intersection with $\cl{V}$, we have $\wh{f}^{N_2} \gamma \subset \wh{f}^{N_2}(\ol{\cl{V}})$, and as $N_2\geq n_2\geq n_1$, by item \ref{it.V7} we obtain 
\begin{equation}    \label{eq.60a}
\wh{f}^{N_2} \gamma \subset   R(\ell_1). 
\end{equation}

Also by the fact that $\gamma$ has good intersection with $\cl{V}$, $\gamma\cap R_{\infty}^1=\empt$, and then $\wh{f}^{N_2}\gamma\cap R_{\infty}^1=\empt$ by the invariance of $R_{\infty}^1$. As $R_{\infty}^2\subset R_{\infty}^1$,
\begin{equation}   \label{eq.aa8}
 \wh{f}^{N_2}\gamma \cap R_{\infty}^2 = \empt.
\end{equation}
By this, by (\ref{eq.60a}), and as the points $\wh{f}^{N_2}\gamma(0)$, $\wh{f}^{N_2}\gamma(1)$ belong to different connected components of $R(\ell_1)\minus F_i^{1,j_1}$, $i=1,2$, we conclude that the arc $\wh{f}^{N_2}\gamma$ must intersect both sets $L_1^{1,j_1}\minus R_{\infty}^2$ and $L_2^{1,j_1}\minus R_{\infty}^2$.

Consider a subarc $\bar{\gamma} \subset \wh{f}^{N_2}\gamma$ with one endpoint in $L_1^{1,j_1}\minus R_{\infty}^2$, the other endpoint in $L_2^{1,j_1}\minus R_{\infty}^2$, and minimal with respect to this property. We claim that $\txt{int}(\ol{\gamma}) \subset \txt{int}(\cl{V}_{1,j_1})$. To see this, note that by item \ref{it.V10}, the set $\cl{V}_{1,j_1}$ is of the form $\cl{V}_{1,j_1}=(V_1\cap V_2)\cup R_1^{1,j_1}$ for some open sets $V_1,V_2$ such that $L_1^{1,j_1}\subset V_2\cup R_2$, $L_2^{1,j_1}\subset V_1\cup R_1$ and $\pr V_i \subset \ell_1\cup F_i^{1,j_1}$. By the minimality of $\bar{\gamma}$ and by (\ref{eq.60a}), (\ref{eq.aa8}) we have
$$ \txt{int}(\bar{\gamma}) \cap \pr V_i \subset \txt{int}(\bar{\gamma}) \cap L_i^{1,j_1} = \empt\ \ \ \ \txt{for $i=1,2$}. $$
Also, by definition $\ol{\gamma}(0)\in L_1^{1,j_1} \subset V_2$, $\ol{\gamma}(1)\in L_2^{1,j_1} \subset V_1$ and therefore 
$$\txt{int}(\ol{\gamma}) \subset V_1\cap V_2 \subset \txt{int}(\cl{V}_{1,j_1}),$$
as claimed. By (\ref{eq.a88}) we actually have int$(\ol{\gamma})\subset\txt{int}(\cl{V}_{1,j_1})\minus R^2_{\infty}$, and then $\ol{\gamma}$ has good intersection with $\cl{V}_{1,j_1}$ which proves item 1 of the proposition. 

We now prove item 2. As the arc $J$ is such that $J\subset\ell$, $J(0)\in\ol{L}_2$ and $J(1)\in\ol{L}_1$ then $J$ has good intersection with $\cl{V}$. By the previous item, $\wh{f}^{N_2}J$ contains an arc which has good intersection with $\cl{V}_{1,j_1}$, for some $j_1\in\Z$.
\end{prueba}

\subsubsection{Proof of Lemma \ref{conjV1}; construction of the sets  $L_i$, $R_i$ and $\cl{V}$.}    \label{sec.conjv1}

\paragraph{Some preliminary results.}

\noindent In this section we use Handel's theorem \ref{teo.handel} to construct the sets $L_1, L_2$. Recall that by hypothesis on Theorem D, $\rho(\wh{f})$ contains two vectors $(0,a)$, $(0,b)$, with $a<0<b$.

\begin{proposicion}    \label{propcurvalibre1}
For all $r\in [a,b]$ there exists a closed $\wh{f}$-invariant set $K_r\subset (\ell, \ell_1)$ such that $\rho(K_r,\wh{f})= \{(0,r)\}$.
\end{proposicion}

Handel's theorem \ref{teo.handel} is a theorem for annulus maps. Lemma \ref{lemacurvalibre3} below will give us that we may assume we are dealing with an annulus map, to the one we will be able to apply Handel's theorem.

Let $\tl{f}:\R\times\T^1\ra\R\times\T^1$ be the lift of $f$ to $\R\times\T^1$ such that $\pi\circ\wh{f} = \tl{f}\circ\pi$, where $\pi:\R^2\ra\R\times \T^1$ denotes the canonical projection. Let $\ol{\ell}=\pi(\ell)$ and $\ol{\ell}_1=\pi(\ell_1)$. Let $m>0$ be such that $\txt{pr}_1( (\wh{f}^{-1}(\ell), \wh{f}(\ell_1)) ) \subset [-m,m]\subset\R$, and denote by $A_m$ the compact annulus given by $A_m = [-m,m] \times \T^1\subset \R\times \T^1$.

\begin{lema}   \label{lemacurvalibre3}
There exists a homeomorphism homotopic to the identity $g:A_m\ra A_m$ such that:
\begin{enumerate}
\item $g|_{(\ol{\ell},\ol{\ell}_1)} = \tl{f}|_{(\ol{\ell},\ol{\ell}_1)}$.  \label{itemcl3}
\item For all $x$ in $R(\ol{\ell}_1)$, the $\omega$-limit $\omega(x,g)$ is contained in $\{m\}\times \T^1$, and for all $x$ in $L(\ol{\ell})$, the $\alpha$-limit $\alpha(x,g)$ is contained in $\{-m\}\times \T^1$. \label{itemcl1}
\item If $\wh{g}:[-m,m]\times\R \ra [-m,m]\times \R$ is the lift of $g$ such that $\wh{g}|_{(\ell,\ell_1)} = \wh{f}|_{(\ell,\ell_1)}$, then $\rho(\{-m\}\times\T^1,\wh{g})=\{(0,s_1)\}$ and $\rho(\{m\} \times \T^1,\wh{g}) = \{(0,s_2)\}$, for some $s_1,s_2\in\Z$ with $s_1 < \min \txt{pr}_2(\rho(\wh{f})) \leq a < b \leq  \max \txt{pr}_2(\rho(\wh{f})) < s_2$.               \label{itemcl2}
\end{enumerate}
\end{lema}

\begin{prueba}
For $0<r_1<r_2$, denote the complex annulus $\A_{r_1,r_2}=\{ re^{i\theta} \in\C  \, : \, r_1 \leq r \leq r_2  \}$. Consider the homeomorphism $h:A_m\to \A_{e^{-m},e^m}$ given by $h(x,y)=e^{-x}e^{iy}$, which sends horizontal segments $\{y=y_0\}$ in $A_m$ to radial segments $\{\theta= e^{y_0} \}$ in $\A_{e^{-m},e^m}$, and sends the border components $\{-m\}\times\T^1$ and $\{m\}\times\T^1$ to $\{r=e^m\}$ and $\{r=e^{-m}\}$, respectively. Up to a change of coordinates (isotopic to the identity), we may assume that the circle $\tl{f}(\ol{\ell}_1)$ is a straight vertical circle, which corresponds by $h$ to a circle $c_{s}=\{r=s\}\subset \A_{e^{-m},e^m}$. 

For $r>0$, denote $\D_r=\{ z \in \C \, : \, |z| < r\}$. Let $F=h\tl{f} h^{-1}: \A_{e^{-m},e^m}\ra \A_{e^{-m},e^m}$ and note that, as $\tl{f}^2(\ol{\ell}_1)\subset R(\tl{f}(\ol{\ell}_1))$, it holds $F(c_{s})=F(\pr\D_{s})\subset \D_{s}$. By a Theorem of Moser (\cite{moser}), if $t\in(r_1,s)$ there exists a homeomorphism $\phi: \A_{e^{-m},s} \to \A_{e^{-m},s}$ such that $\phi|_{\pr \A_{e^{-m},s}} = \txt{Id}$, and such that $\phi(F(c_{s})) = \pr{\D}_{t}$. Composing with a twist map with support on an annulus $U\subset \A_{e^{-m},s}$ which contains $\pr \D_{t}$ we may also suppose that $\phi F|_{c_s}$ fixes angles, that is, we may suppose that $\phi$ is such that for any $z\in c_{s}$, arg$(\phi  F(z))= \txt{arg}(z)$, where arg$(z)$ denotes the argument of $z$. Let $G_1= \phi F$ and note that, as $\phi$ has support on $\A_{e^{-m},s}$ and as $\tl{f}(\ol{\ell}_1)=h^{-1}(c_s)$,
\begin{equation}  \label{eq.aa9}
 G_1|_{h(L(\ol{\ell}_1))} = F|_{h(L(\ol{\ell}_1))},
\end{equation}
and then the map $g_1:=h^{-1} G_1 h:A_m\to A_m$ is such that 
$$g_1|_{L(\ol{\ell}_1)}\equiv \tl{f}|_{L(\ol{\ell}_1)}.$$
Note that, as $G_1|_{c_s}$ fixes angles, $g_1|_{\tl{f}(\ol{\ell}_1)}$ fixes the second coordinate. By an adequate choice of the mentioned twist map, we may assume that, not only $g_1|_{\tl{f}(\ol{\ell}_1)}$ fixes the second coordinate, but also, if $\wh{g}_1:[-m,m]\times\R \to [-m,m]\times\R$ is the lift of $g_1$ such that $\wh{g}_1|_{L(\ell_1)} \equiv \wh{f}|_{L(\ell_1)}$, we have 
\begin{equation}    \label{eq.z1z}
\txt{pr}_2(\wh{g}_1(\wh{z})) = \txt{pr}_2(\wh{z}) + s_2
\end{equation}
for all $\wh{z}\in \wh{f}(\ell_1)$ and where $s_2\in\Z$, $s_2 >  \max \txt{pr}_2(\rho(\wh{f}))\geq b >0$. 

Consider the restriction $G_1|_{c_{r_3}}: c_{r_3}\ra \A_{e^{-m},e^m}$, and radially extend it to an injective map $G_2: \A_{e^{-m},s}  \ra \A_{e^{-m},s}$, so that for every $re^{i\theta}\in \A_{e^{-m},s}$, $G_2(re^{i\theta})=\tl{r}e^{i \theta}$, with $\tl{r}\leq r$, and with equality iff $r=e^{-m}$. Then, extend $G_2$ continuously to $\A_{e^{-m},e^m}$ as $G_2(x)= G_1(x)$ for $x\in \A_{s,e^m}$.

Observe that as $G_2|_{\A_{s,e^m}} \equiv G_1$, by (\ref{eq.aa9}) we have $G_2|_{h(L(\ol{\ell}_1))} \equiv F|_{h(L(\ol{\ell}_1))}$, and then the map $g_2:=h^{-1}G_2 h:A_m\to A_m$ is such that 
\begin{equation}   \label{eq.z2z}
g_2|_{L(\ol{\ell}_1)} \equiv \tl{f}|_{L(\ol{\ell}_1)}.
\end{equation}

Let $\wh{g}_2:[-m,m]\times\R \to [-m,m]\times\R$ be the lift of $g_2$ such that $\wh{g}_2|_{L(\ol{\ell}_1)}\equiv \wh{f}|_{L(\ol{\ell}_1)}$. As $g_2|_{R(\ol{\ell}_1)} \equiv g_1|_{R(\ol{\ell}_1)}$ (by definition) and by (\ref{eq.z1z}) we have that pr$_2(\wh{g}_2(x)) = \txt{pr}_2(x)$ for $x\in R(\ell_1)$. Also, as $G_2(\{r=e^{-m}\})= \{r=e^{-m}\}$,
\begin{equation}   \label{eq.z4z}
\wh{g}_2(x)=x+ (0,s_2) \ \ \ \txt{for $x\in\{m\}\times\R$}.
\end{equation}

Being the map $G_2$ a radial contraction on $\A_{e^{-m},s}$, we have that for any $z\in \A_{e^{-m},s}$, $\omega(z,G_2)\subset \{r=e^{-m}\}$. This implies 
\begin{equation}   \label{eq.z3z}
\omega(x,g_2)\subset \{m\}\times\T^1 \ \ \ \txt{for any $x\in \ol{R}(\ol{\ell}_1)$}.
\end{equation}

In a symmetric way, we may modify $g_2$ to obtain a homeomorphism $g:A_m\ra A_m$ such that:
\begin{itemize}
\item $g|_{R(\ol{\ell})} = g_2|_{R(\ol{\ell})}$,
\item $\alpha(x,g)\subset \{-M\}\times \T^1$ for any $x\in \ol{L}(\ol{\ell})$,
\end{itemize}
and such that, if $\wh{g}:\R^2\ra\R^2$ is the lift of $g$ such that $\wh{g}|_{(\ell,\ell_1)} \equiv \wh{g}_2|_{(\ell,\ell_1)} \equiv \wh{f}|_{(\ell,\ell_1)}$, then 
\begin{itemize}
\item $\wh{g}(x)=x+ (0,s_1)$,
\end{itemize}
for all $x\in\{-m\}\times\R$ and where $s_1\in\Z$, $s_1 < \min \txt{pr}_2(\rho(\wh{f})) \leq a < 0$. By this, and by eqs. (\ref{eq.z2z}), (\ref{eq.z4z}), (\ref{eq.z3z}) we have that $g$ satisfies items 1 to 3 of the lemma.
\end{prueba}

We are now ready to prove Proposition \ref{propcurvalibre1}, using Handel's theorem \ref{teo.handel}.

\begin{prueba}[Proof of Proposition \ref{propcurvalibre1}]
Let $g:A_m\ra A_m$ be the annulus homeomorphism given by Lemma \ref{lemacurvalibre3}, and let $\wh{g}: [-m,m] \times \R \ra  [-m,m] \times \R$ be the lift of $g$ such that $\wh{g}|_{(\ell,\ell_1)} \equiv \wh{f}|_{(\ell,\ell_1)}$. Consider the sets $\rho(\wh{g})$ and $\rho_{point}(\wh{g})$ as subsets of $\{0\}\times\R\subset\R^2$. Let $E_g$ be the maximal invariant set of $(\ol{\ell},\ol{\ell}_1)$ for $g$, and let $B=\{-m,m\}\times\T^1$. As the curves $\ol{\ell}$ and $\ol{\ell}_1$ are free for $g$, items (\ref{itemcl1}) and (\ref{itemcl2}) of Lemma \ref{lemacurvalibre3} imply that
\begin{equation}  \label{eqcl7}
\rho_{point}(\wh{g})=\rho_{point}(E_g\cup B, \wh{g}) = \rho_{point}( E_g,\wh{g}) \cup\{(0,s_1),(0,s_2)\}.
\end{equation}
As $\wh{g}|_{(\ell,\ell_1)} = \wh{f}|_{(\ell,\ell_1)}$, 
\begin{equation}  \label{eqcl8}
\rho_{point}(E_g, \wh{g})= \rho_{point}(E_{\tl{f}},\wh{f}),
\end{equation}
where $E_{\tl{f}}$ is the maximal invariant set of $(\ol{\ell},\ol{\ell}_1)$ for $\tl{f}$. Observe that, by Assumption \ref{asuncion.per}, 
$$\rho_{point}(E_{\tl{f}},\wh{f}) \supset (\{0\}\times [a,b]) \cap\Q^2,$$ 
and then equations (\ref{eqcl7}) and (\ref{eqcl8}) imply 
\begin{equation}
\rho_{point}(\wh{g}) = \rho_{point}(E_{\tl{f}},\wh{f}) \cup \{(0,s_1),(0,s_2)\}  \supset  (\{0\}\times [a,b]) \cap\Q^2.
\end{equation}

By Handel's theorem \ref{teo.handel}, $\rho_{point}(\wh{g})$ is closed, and therefore $\rho_{point}(\wh{g}) \supset  \{0\}\times [a,b]$. Also by that theorem, for all $r\in [a,b]$ there exists a compact set $K_r\subset A_m$ that is $g$-invariant, and such that $\rho(z,\wh{g})=r$ for all $z\in K_r$. 

By item 2 from Lemma \ref{lemacurvalibre3}, the forward $g$-iterates of $\ell_1$ accumulate on $\{m\}\times\T^1$, the backward $g$-iterates of $\ell$ accumulate on $\{-m\}\times\T^1$, and by item 3 of the same lemma, $\rho(x,\wh{g})=\{(0,s_1)\}$ for any $x\in\{-m\}\times\T^1$ and $\rho(x,\wh{g})=\{(0,s_2)\}$ for all $x\in\{m\}\times\T^1$, where $s_1 < a < b < s_2$. The fact that $K_r$ is $g$-invariant and $r\in[a,b]$ implies then that $K_r\subset (\ol{\ell},\ol{\ell}_1)$. As $g|_{(\ol{\ell},\ol{\ell}_1)} = \wt{f}|_{(\ol{\ell},\ol{\ell}_1)}$, we have that $K_r$ is $\tl{f}$-invariant and $\rho(K_r,\wh{f})= \rho(K_r,\wh{g}) = \{(0,r)\}$. This proves the proposition.
\end{prueba}

We now use Proposition \ref{propcurvalibre1} to construct the sets $L_1,L_2$. Such construction will be carried out basically with the use of Corollary \ref{corolariocl1} below.

\begin{lema}   \label{lemacurvalibre1}
For all $p\in [a,b]$ there exists $z\in L_{\infty}^1\cap R_{\infty}^1$ such that, if $C$ is the connected component of $L_{\infty}^1 \cap R(\ell)$ that contains $z$, then:
\begin{enumerate}
\item $\ol{C}$ is compact.
\item $\rho(\ol{C},\wh{f})=\{(0,p)\}$.
\item $\ol{C}\cap \ell\neq\empt$.
\item $\ol{C} \cap L_{\infty}^0 = \empt$, and therefore there is $n>0$ such that $\wh{f}^{n}(C)\subset R(\ell_1)$. 
\end{enumerate}
\end{lema}

\begin{corolario}     \label{corolariocl1}
There exists an arc $I\subset\ell$, and two points $x,y\in L_{\infty}^1\cap R_{\infty}^1$ such that, if $L_x$ and $L_y$ are the connected components of $L_{\infty}^1 \cap R(\ell)$ that contain $x$ and $y$, respectively, then: 
\begin{enumerate}
\item $\ol{L}_x, \ol{L}_y$ are compact.
\item $\txt{int}(I) \cap L_{\infty}^0 = \empt$, and the endpoints of $I$ lie in $L_{\infty}^0$.   
\item $\ol{L}_i\cap \txt{int}(I) \neq\empt$, for $i=x,y$.
\item $\ol{L}_i\cap L_{\infty}^0=\empt$ and there is $n_1>0$ such that $\wh{f}^{n_1}(\ol{L}_1\cup \ol{L}_2)\subset R(\ell_1)$.
\item $\rho(\ol{L}_i, \wh{f})=\{(0,r_i)\}$, for $i=x,y$ and some $r_x\neq r_y$. 
\end{enumerate}
\end{corolario}

\begin{prueba}
By Lemma \ref{lemacurvalibre1}, for any $p\in[a,b]$ there is $z_p\in L_{\infty}^1 \cap R_{\infty}^1$ such that, if $C_{z_p}$ denotes the connected component of $L_{\infty}^1\cap R(\ell)$ that contains $z_p$, then 
\begin{equation}   \label{eqcl2}
\rho(\ol{C}_{z_p},\wh{f})=\{(0,p)\},
\end{equation} 
\begin{equation}  \label{eqcl1}
\ol{C}_{z_p}\cap L_{\infty}^0=\empt,
\end{equation}
\begin{equation}   \label{eq.aaa}
\ol{C}_{z_p}\cap\ell\neq\empt
\end{equation} 
and there is $n_p>0$ such that 
\begin{equation}    \label{eq.aa1}
\wh{f}^{n_p}(\ol{C}_{z_p}) \subset R(\ell_1).
\end{equation}

Let $U=\ell\minus L_{\infty}^0$. As $L_{\infty}^0$ is closed, $U$ is a countable union of open arcs. By (\ref{eq.aaa}) and as the set $[a,b]$ is uncountable, there must be points $r,s\in[a,b]$, $r\neq s$, points $z_r,z_s\in L_{\infty}^1\cap R_{\infty}^1$ and a connected component $U_0$ of $U$ such that
\begin{equation}  \label{eqcl3}
\ol{C}_{z_r}\cap U_0\neq\empt \ \ \txt{ and } \ \ \ol{C}_{z_s}\cap U_0\neq\empt,
\end{equation}
and, letting $N_1=\max\{n_r,n_s\}$, by (\ref{eq.aa1}) we have
\begin{equation}  \label{eqcl4}
\wh{f}^{N_1}(\ol{C}_{z_r}\cup \ol{C}_{z_s})\subset R(\ell_1).
\end{equation}

Let $I=U_0$. By definition, we have that $I$ satisfies item 2 of the corollary. Define $x=z_r$ and $y=z_s$, so $L_x=C_{z_r}$ and $L_y=C_{z_s}$. By item 1 of Lemma \ref{lemacurvalibre1}, the sets $\ol{L}_x$ and $\ol{L}_y$ are compact, and item 1 of the corollary holds. By eqs. (\ref{eqcl1}) and (\ref{eqcl4}), $x$ and $y$ satisfy item $(4)$ of the corollary. By eqs. (\ref{eqcl3}) and (\ref{eqcl2}), items $3$ and $5$, respectively, hold. 
\end{prueba}

To prove Lemma \ref{lemacurvalibre1} we will need the following.

\begin{lema}     \label{lemacurvalibre2}
There does not exist a compact connected set $K$ such that:
\begin{enumerate}
\item $K \cap L_{\infty}^0\neq\empt$, and $K \cap R_{\infty}^2\neq\empt$,
\item $\rho(K, \wh{f})$ consists of a point. 
\end{enumerate}
\end{lema}

\begin{figure}[h] 
\begin{center} 
\includegraphics{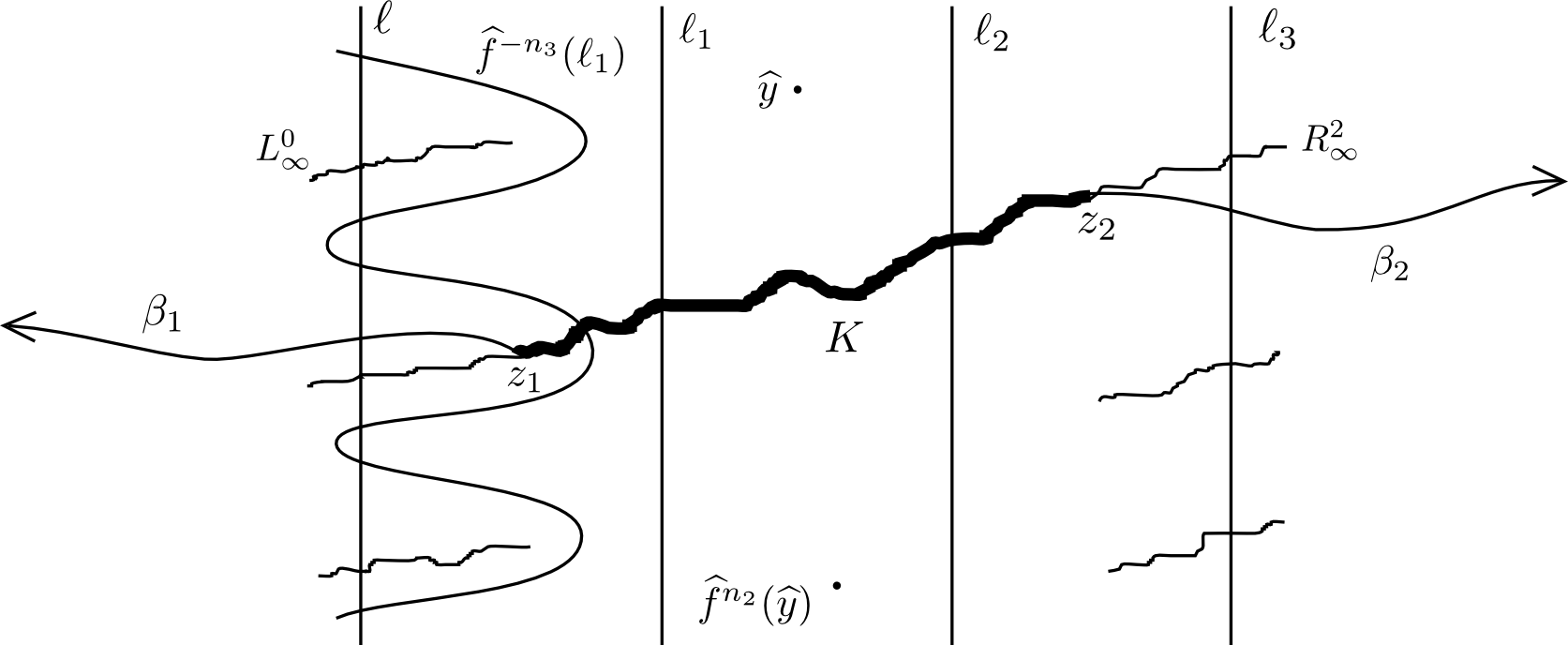}
\caption{Proof of Lemma \ref{lemacurvalibre2}.}
\label{fig.lemaK}  
\end{center}  
\end{figure}

\begin{prueba}
We proceed by contradiction. Suppose that $K$ is a compact connected set such that $K \cap L_{\infty}^0\neq\empt$, $K \cap R_{\infty}^2\neq\empt$, and $\rho(K, \wh{f}) = \{(0,p_0)\}$, for some $p_0 \neq 0$ (see Fig. \ref{fig.lemaK}). We treat the case $p_0\geq 0$, the case $p_0<0$ being similar. 

By Assumption \ref{asuncion.per} there is a periodic point $y\in\T^2$ such that $\rho(y,\wh{f})=(0,c)$, with $c<0$. We choose a lift $\wh{y}\in\R^2$ of $y$ such that $\wh{y}\in(\ell_1,\ell_2)$, and such that $\wh{y}$ is above $K$. As $K$ is compact, as $\rho(K,\wh{f})=\{(0,p_0)\}$ and as $\rho(y,\wh{f})=(0,c)$, $c<0$, there is $n_1>0$ such that $\wh{f}^{n}(\wh{y})$ is below $\wh{f}^n(K)$ for all $n\geq n_1$. Observe that the whole orbit of $\wh{y}$ must be contained in $(\ell_1,\ell_2)$, as the curves $\ell_1$ and $\ell_2$ are Brouwer lines for $\wh{f}$, $y$ is periodic, and $\rho(y,\wh{f})$ is a vertical vector. 

Let $z_1\in K\cap L_{\infty}^0$. Then by definition of $L_{\infty}^0$, $z\in L(\wh{f}^{-n_1}(\ell_1))$. Let $\beta_1:(-\infty,0]\ra\R^2$ be a proper immersion such that:
\begin{itemize}
\item $\beta_1\subset L(\wh{f}^{-n_1}(\ell_1))$,
\item $\beta_1(0)=z_1$, and
\item $-\infty < \inf\txt{pr}_2(\beta_1) < \sup\txt{pr}_2(\beta_1) < \infty$.
\end{itemize}
Now, let $z_2\in K\cap R_{\infty}^2$, and let $\beta_2:[0,\infty)\ra\R^2$ be a proper immersion such that:
\begin{itemize}
\item $\beta_2(0)=z_2$, 
\item $\beta_2\subset R(\ell_2)$, and
\item $-\infty < \inf\txt{pr}_2(\beta_2) < \sup\txt{pr}_2(\beta_2) < \infty$.
\end{itemize}
The set $\beta_1\cup K\cup \beta_2$ is therefore bounded vertically, unbounded horizontally, and separates $\R^2$. The complement of $\beta_1\cup K \cup \beta_2$ has exactly one connected component $V$ unbounded from above and bounded from below. By definition, the point $\wh{y}\in(\ell_1,\ell_2)$ is above $K$, and as $(\beta_1\cup\beta_2) \cap (\ell_1,\ell_2)=\empt$, we have that $y\in V$. 

By the periodicity of $\wh{f}$, $\wh{f}^{n_1}(V)$ is also unbounded from above and bounded from below. By construction of $\beta_1$ and $\beta_2$ we have that $\wh{f}^{n_1}(\beta_1)\subset L(\ell_1)$ and $\wh{f}^{n_1}(\beta_2)\subset R(\ell_2)$, that is
$$ \wh{f}^{n_3}(\beta_1\cup \beta_2) \cap (\ell_1,\ell_2) = \empt.$$
As $\wh{f}^{n_1}(y)\in (\ell_1,\ell_2)$ is below $\wh{f}^{n_1}(K)$ and as $\wh{f}^{n_1}(\beta_1\cup\beta_2)\cap(\ell_1,\ell_2)=\empt$, we have that $\wh{f}^{n_1}(y)$ belongs to a connected component of $\R^2\setminus \wh{f}^{n_1}(\beta_1\cup K\cup \beta_2)$ which is unbounded from below. This contradicts the fact that $\wh{f}^{n_1}(y)$ belongs to the set $\wh{f}^{n_1}(V)$, which is bounded from below, and this contradiction finishes the proof of the lemma.
\end{prueba}

We end up this section with the proof of Lemma \ref{lemacurvalibre1}. 

\begin{figure}[h] 
\begin{center} 
\includegraphics{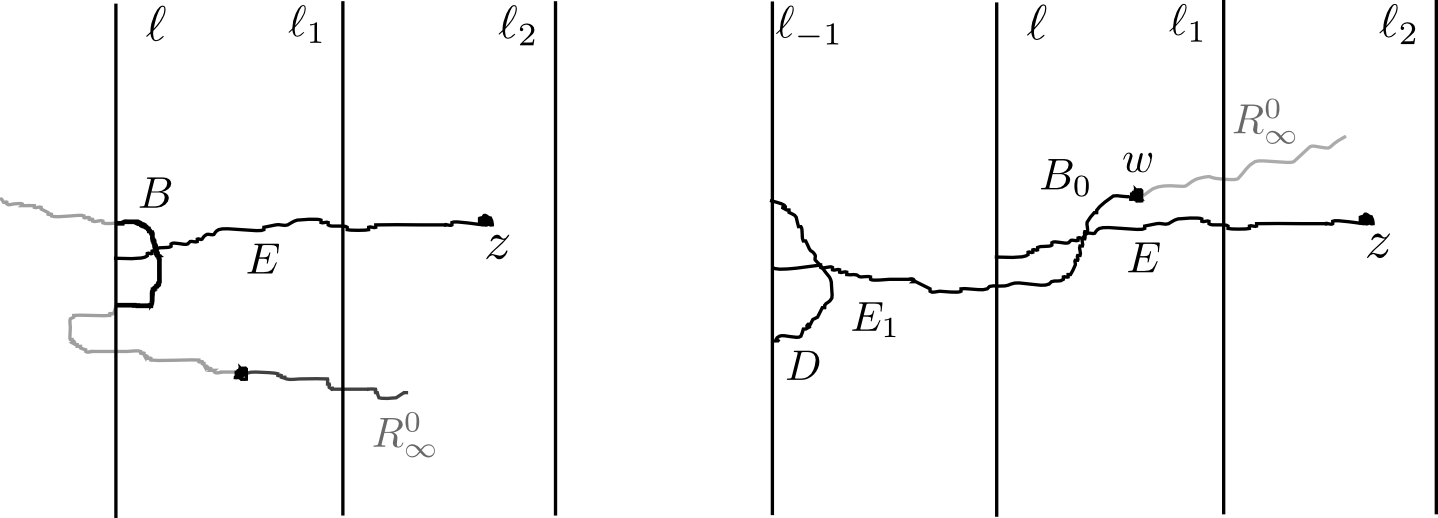}
\caption{Proof of Lemma \ref{lemacurvalibre1}. Left: Case 2.1. Right: Case 2.2-(ii).}
\label{fig.lemaL}  
\end{center}  
\end{figure}

\begin{prueba}[Proof of Lemma \ref{lemacurvalibre1}]
Let $p\in[a,b]$. By Proposition \ref{propcurvalibre1} and by the periodicity of $\wh{f}$, there exists a closed $\wh{f}$-invariant set $K\subset(\ell_1,\ell_2)$ such that $\rho(K,\wh{f})=\{(0,p)\}$. Observe that as $K\subset(\ell_1,\ell_2)$ is $\wh{f}$-invariant, $K\subset L_{\infty}^1\cap R_{\infty}^1$.

Let $z\in K$, and let $E$ be the connected component of $L_{\infty}^1\cap R(\ell)$ that contains $z$. By Lemma \ref{lema.Lcomp} we have that 
\begin{equation}   \label{eq.aa2}
\txt{$\ol{E}$ is compact, $\rho(\ol{E},\wh{f})=\{(0,p)\}$, and $\ol{E}\cap \ell\neq\empt$}.
\end{equation}

The sought set $C$ will be constructed from $E$. The construction is divided in the following cases. \\
\\
\textit{Case 1:} $\ol{E}\cap L_{\infty}^0=\empt$.

In this case, $\wh{f}^n(\ol{E})\subset R(\ell_1)$ for some $n>0$, by the definition of $L_{\infty}^0$ and by the compacity of $\ol{E}$. Letting $x=z$, by (\ref{eq.aa2}) we have that $x$ satisfies items $(1)$ to $(4)$ of the lemma, with $C=E$.\\
\\
\textit{Case 2:} $\ol{E}\cap L_{\infty}^0\neq\empt$. 

We divide this case in two subcases.\\
\\
\textit{Case 2.1: For every connected component $B$ of $L_{\infty}^0\cap \ol{R}(\ell)$ that intersects $\ol{E}$, it holds $B\cap R_{\infty}^0 =\empt$ (cr. Fig. \ref{fig.lemaL}).}

In this case, for any such component $B$ there is $m_B>0$ such that  $\wh{f}^{-m_B}(B)\subset L(\ell)$. Then, by the compacity of $L_{\infty}^0\cap\ol{E}$, there is $m>0$ such that 
\begin{equation}  \nonumber
\wh{f}^{-m}(\ol{E})\cap (L_{\infty}^0\cap \ol{R}(\ell)) = \empt
\end{equation}
and therefore
\begin{equation} \label{eqcl5}
\ol{\wh{f}^{-m}(E)\cap R(\ell)}\cap L_{\infty}^0 = \empt.
\end{equation}

Let $E_0$ be the connected component of $\wh{f}^{-m}(E)\cap R(\ell)$ that contains $f^{-m}(z)$. Then $\ol{E}_0$ is compact because $\ol{E}$ is, and by (\ref{eqcl5}) we have
$$\ol{E}_0\cap L_{\infty}^0=\empt,$$
so there is $n>0$ such that $f^n(\ol{E}_0)\subset R(\ell_1)$. By Lemma \ref{lema.Lcomp} $\ol{E}_0\cap\ell\neq\empt$ and $\rho(\ol{E}_0,\wh{f})=\{(0,p)\}$. Therefore, letting $x=\wh{f}^{-m}(z)$, we have that $x$ satisfies the conclusions of the lemma, with $C=E_0$.\\
\\
\textit{Case 2.2: There is a connected component $B_0$ of $L_{\infty}^0 \cap \ol{R}(\ell)$ that intersects $\ol{E}$, and such that $B_0\cap R_{\infty}^0\neq \empt$ (see Fig. \ref{fig.lemaL}).}

Let $E_1$ be the connected component of $L_{\infty}^0\cap R(\ell_{-1})$ that contains $B_0$. We have three possibilities for $E_1$: 
\begin{enumerate}[(i)]
\item $\ol{E}_1\cap L_{\infty}^{-1}=\empt$.
\item $\ol{E}_1\cap L_{\infty}^{-1} \neq\empt$.
\end{enumerate}
To deal with case (i), let $w\in B_0\cap R_{\infty}^0\subset L_{\infty}^0\cap R_{\infty}^0$. Proceeding as in Case 1 above, and using the periodicity of $\wh{f}$ we have that the point $x=T_1(w)\in L_{\infty}^1\cap R_{\infty}^1$ satisfies the conclusions of the lemma, with $C= T_1(B_0)$. 

We now deal with case (ii). Let $D$ be a connected component of $L_{\infty}^{-1}\cap \ol{R}(\ell_{-1})$ that intersects $\ol{E}_1$. Recall that through all of Section ref{sec.ss} we are assuming that there are $z\in\R^2,n\in\N$ such that pr$_1(\wh{f}^n(z)-z)>M$, and by the construction of $M$, $\wh{f}^n(\ell_{-1})\cap R(\ell_2)\neq\empt$ (cf. equation (\ref{eq.111})). Then, by Lemma \ref{lema.Lcomp} and as $\ol{E}_1\cup D\cup E \subset L_{\infty}^1$, $\rho(\ol{E}_1\cup D\cup E,\wh{f})$ consists of a point. As $\rho(\ol{E},\wh{f})=\{(0,p)\}$, we then have
\begin{equation}   \label{eqcl6}
\rho(\ol{E}_1\cup D\cup \ol{E},\wh{f})=\{(0,p)\}. 
\end{equation}

The set $\ol{E}_1\cup D \cup \ol{E}$ is compact, intersects $L_{\infty}^{-1}$ (because $D$ does), and intersects also $R_{\infty}^1$ (as $z\in E$). This, together with (\ref{eqcl6}), contradicts Lemma \ref{lemacurvalibre2}. Therefore, case (ii) cannot occur, and this finishes the proof of the lemma. 
\end{prueba}

\paragraph{Definition of the sets $L_i$, $R_i$ and $\cl{V}$.} \mbox{} \\

We are now ready to define the sets $L_i$, $R_i$, $\cl{V}$ and verify that they satisfy Lemma \ref{conjV1}.

\noindent \textbf{The sets $L_i$ and $R_i$.} By Corollary \ref{corolariocl1}, there exists an arc $I\subset\ell$, and two points $x,y\in L_{\infty}^1\cap R_{\infty}^1$ such that the connected components $L_x$ and $L_y$ of $L_{\infty}^1 \cap R(\ell)$ that contain $x$ and $y$, respectively, are such that: 
\begin{itemize}
\item $\ol{L}_x$ and $\ol{L}_y$ are compact,
\item $\ol{L}_i\cap I \neq\empt$, for $i=x,y$,
\item $\ol{L}_i\cap L_{\infty}^0=\empt$, for $i=x,y$,
\item $\rho(\ol{L}_i, \wh{f})=\{(0,r_i)\}$, for $i=x,y$ and some reals $r_x\neq r_y$,
\item $\txt{int}(I) \cap L_{\infty}^0 = \empt$, and the endpoints of $I$ lie in $L_{\infty}^0$.   
\end{itemize}

Define $L_1$ as a connected component of $L_x\minus R_{\infty}^1$ such that $\ol{L}_1\cap I \neq\empt$, and similarily, define $L_2$ as a connected component of $L_y\minus R_{\infty}^1$ such that $\ol{L}_2\cap I\neq \empt$. By the above properties, we have $L_i\subset L_{\infty}^1\minus L_{\infty}^0$ and $\rho(\ol{L}_i,\wh{f})=r_i$, $r_1\neq r_2$, and then the sets $L_i$ satisfy items \ref{it.V4} and \ref{it.V8} from Lemma \ref{conjV1}.

Define $R_1$ and $R_2$ to be connected components of $R_{\infty}^1$ that intersect $\ol{L}_1$ and $\ol{L}_2$, respectively.\footnote{Note that, in principle, the set $R_{\infty}^1$ might be connected, and in such case $R_1=R_2$.} Then item \ref{it.V5} from Lemma \ref{conjV1} holds. 

As $\ol{L}_i\cap I\neq\empt$ and the $R_i$ are unbounded to the right (cf. Lemma \ref{discos}), the sets $F_i=L_i\cup R_i$ separate $R(\ell)$. By definition the $L_i$ and $R_i$ are connected, and by the above mentioned properties the $\ol{L}_i$ are compact. Thus, item \ref{it.V2} from Lemma \ref{conjV1} holds. \\

\noindent \textbf{The set $\cl{V}$.} Item \ref{it.V10} from Lemma \ref{conjV1} states that $\cl{V}$ is of the form $\cl{V}=(V_1\cap V_2) \cup R_1$, for some open sets $V_1,V_2\subset R(\ell)$. We construct now such sets.

\begin{claim}  \label{claimcl2}
Let $i\in\{1,2\}$ and let $p,q\in\ell\minus \ol{F}_i$, $x\in\ell\cap \ol{F}_i$ be points such that $p<x<q$, where $<$ denotes the order of $\ell$ induced by its upwards orientation. Then, $p$ and $q$ belong to different connected components of $\ol{R}(\ell) \minus \ol{F}_i$.
\end{claim}
\begin{prueba}
This follows easily from the fact that $F_i$ is connected and separates $R(\ell)$. 
\end{prueba}

A consequence of Claim \ref{claimcl2} is that, as the sets $F_i$ separate $R(\ell)$ and as $\ol{L}_1\cap \ol{L}_2=\empt$, either $\min \ol{L}_1\cap \ell > \max \ol{L}_2\cap\ell$ or $\min \ol{L}_2\cap \ell > \max \ol{L}_1\cap\ell$. Without loss of generality, suppose it holds the former;
\begin{equation}   \label{eq.aa6}
\min \ol{L}_1\cap \ell > \max \ol{L}_2\cap\ell,
\end{equation} 
that is, in some sense $L_1$ is `above' $L_2$.

As $\ol{L}_1\cap \ell$ is compact (because $\ol{L}_1$ is), and as $R_1 \subset R_{\infty}^1 \subset R(\ell_1)$, we have that $\ol{F}_1\cap\ell$ is compact, and therefore there is only one connected component $V_1$ of $R(\ell)\minus F_1$ whose closure contains a subcurve of $\ell$ that is unbounded from below. Analogously, there is only one connected component $V_2$ of $R(\ell)\minus F_2$ whose closure contains a subcurve of $\ell$ that is unbounded from above (see Fig. \ref{fig.V1}). 

Note that, by definition, $\pr V_i\subset F_i\cup\ell$, for $i=1,2$. Therefore, to show that the set
$$\cl{V} := ( V_1\cap V_2)\cup R_1$$
satisfies item \ref{it.V10} from Lemma \ref{conjV1}, it suffices to show the following:
\begin{claim}  \label{claimcl1}
$L_1\subset V_2$, and $L_2\subset V_1$.
\end{claim}
\begin{prueba}[Proof of Claim \ref{claimcl1}.]
We will prove that $L_1 \subset V_2$, the proof that $L_2\subset V_1$ being symmetric. First observe that by (\ref{eq.aa6}) and by definition of $V_2$, it follows that $L_1\cap V_2\neq\empt$. As $L_1$ is by definition connected and disjoint from $R_2$, it holds 
\begin{equation}   \nonumber
L_1 \subset V_2.
\end{equation}
\end{prueba}

Note that by (\ref{eq.aa6}) and by definition of $V_1,V_2$, it follows that $\empt\neq V_1\cap V_2 \subset \txt{int}(\cl{V})$, $\ol{\cl{V}}\cap \ell\neq\empt$ and as $\cl{V}$ is unbounded to the right, $\cl{V}$ separates $R(\ell)$. Also by definition, $\pr \cl{V} \subset F_1\cup F_2 \cup \ell$, and therefore item \ref{it.V1} from Lemma \ref{conjV1} holds. 

As the sets $F_i$ separate $R(\ell)$ and as $\ol{F}_i\cap [\ell,\ell_1] = \ol{L}_i\cap [\ell,\ell_1]$ is compact, we have that $\ol{V}_1\cap [\ell,\ell_1]$ is bounded from above and $\ol{V}_2\cap[\ell,\ell_1]$ is bounded from below. Therefore $\ol{\cl{V}}\cap[\ell,\ell_1] \subset (\ol{V}_1\cap\ol{V}_2) \cap [\ell,\ell_1]$ is compact, which proves item \ref{it.V9}.

Observe that, also by (\ref{eq.aa6}) and by definition of $\cl{V}$, the arc $I$ mentioned above and given by Corollary \ref{corolariocl1} contains a non-degenerate arc $J$ such that $J = \pr\cl{V}\cap\ell$, and such that $J(0)\in \ol{L}_2$, $J(1)\in \ol{L}_1$ and $J\cap L_{\infty}^0=\empt$. Thus, item \ref{it.V6} holds. 

By last, we prove item \ref{it.V7}, namely, there is $n_1>0$ such that $\wh{f}^{n_1}(\cl{V})\subset R(\ell_1)$. To this end, it suffices to show that $\ol{\cl{V}}\cap L_{\infty}^0=\empt$. Suppose this is not the case, and let $x\in \ol{\cl{V}}\cap L_{\infty}^0$. Consider the connected component $C$ of $L_{\infty}^0$ that contains $x$. As $\pr\cl{V}\subset J\cup F_1\cup F_2$, and as the sets $J, F_1,F_2$ do not intersect $L_{\infty}^0$, we then have that $C$ is contained in int$(\cl{V})\subset R(\ell)$. This contradicts the fact that $C$ is unbounded to the left (Lemma \ref{lema.Lcomp}). We therefore must have $\ol{\cl{V}}\cap L_{\infty}^0=\empt$, and item \ref{it.V7} holds.

\subsection{Case that $\cl{C}$ contains singularities}       \label{sec.cs}

As we did in the case that $\cl{C}$ does not contain singularities, we will prove that if $M=M(f)>0$ is the constant constructed in Section \ref{sec.M} and if there are $x\in\R^2$ and $N_1>0$ such that pr$_1(\wh{f}^{N_1}(x)-x)>M$, then there are $z\in\R^2$ and $N>0$ such that pr$_1(\wh{f}^{nN}(z)-z) > n$ for all $n\in\N$.

We begin by giving an idea of the proof.

\subsubsection{Idea of the proof}     \label{sec.ideasing}

As we mentioned, we are assuming that there is $x\in\R^2$ and $N_1>0$ such that pr$_1(\wh{f}^{N_1}(x)-x) > M$. Recall from (\ref{eq.a.1}) that this implies $\wh{f}^{N_1}(\cl{C}_{i})\cap\cl{C}_{i+2}\neq\empt$ for any $i$, and in particular there are leaves $\gamma\subset\cl{C}_{-2}$ and $\beta\subset\cl{C}_0$ such that $\wh{f}^{N_1}\gamma\cap\beta\neq\empt$. Recall also from (\ref{eq.a.3}) that there is $p\in\R^2$ such that $\pi(p)$ is periodic for $f$ and $\rho(\pi(p),\wh{f})=(0,c)$, for some $c<0$. We will first choose an integer translate of $p$, denoted also $p$, such that $p\in R(\cl{C}_{-2})\cap L(\cl{C}_0)$ and $p$ is above $\wh{f}^{N_1}(\gamma)$. 

We will show that, as the isotopy $(f_t)_t$ is positively transverse to the foliation $\wh{\cl{F}}$, all the iterates $\wh{f}^{N_1+n}(\gamma)$, $n>0$, also intersect $\beta$, which might be interpreted as $\wh{f}^{N_1}(\gamma)$ being `anchored' to $\beta$. 

\begin{figure}[h]        
\begin{center} 
\includegraphics{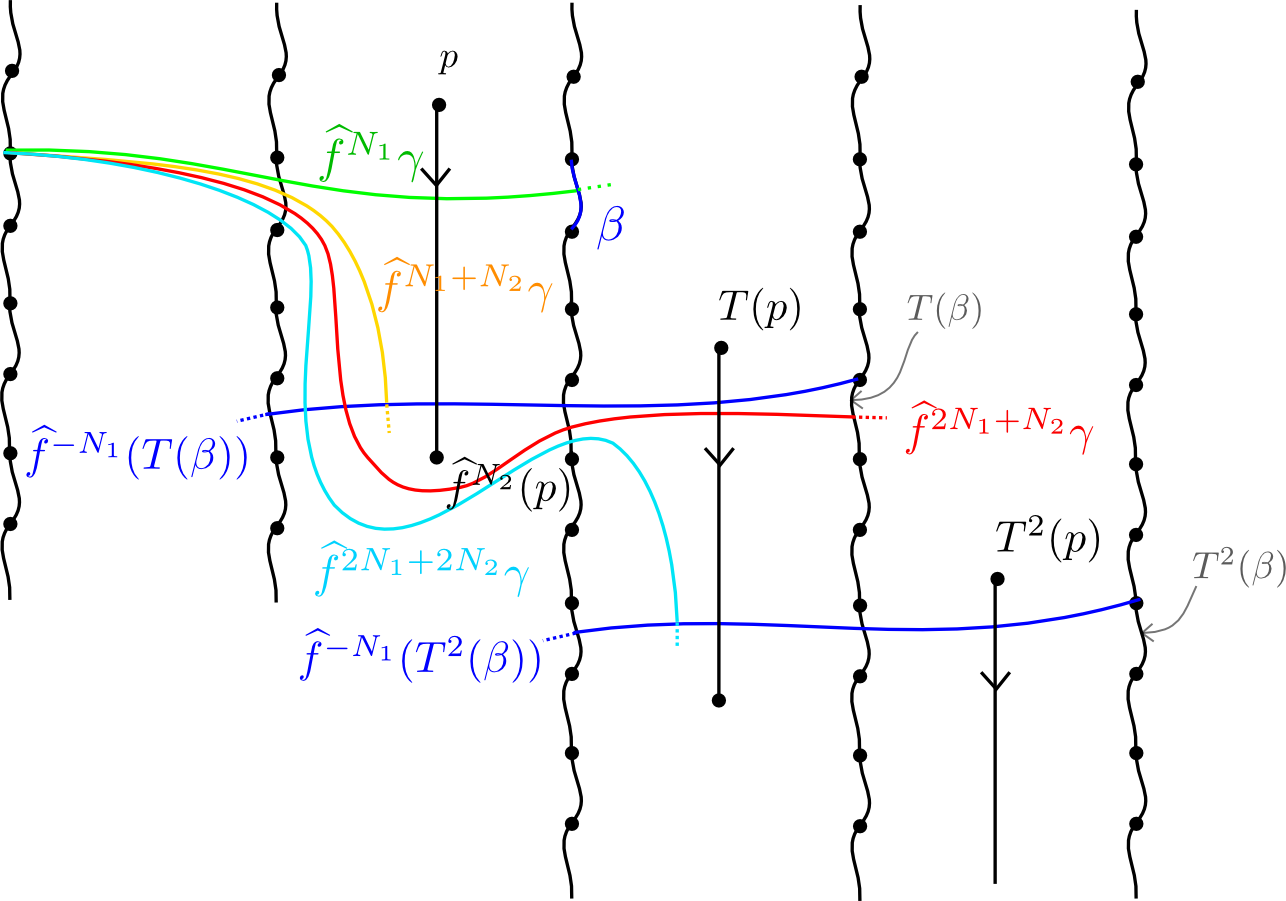}
\caption{}
\label{fig.idea}
\end{center}  
\end{figure}

We will choose an integer translation $T:\R^2\to\R^2$ such that $T(\beta)\subset \cl{C}_1$ and $T(\beta)$ is below $\beta$. Observe that by the periodicity of $\wh{f}$, $\wh{f}^{-N_1}(T(\beta))\cap \cl{C}_{-1}\neq\empt$. We will then choose $N_2>0$ such that $\wh{f}^{N_2}(p)$ is below $\wh{f}^{-N_1}(T(\beta))$ (see Fig. \ref{fig.idea}). As $\wh{f}^{N_1}\gamma$ is `anchored' to $\beta$, the iterates of $p$ by $\wh{f}$ will `push' the curve $\wh{f}^{N_1}\gamma$ downwards, and we will show that
$$\wh{f}^{N_1+N_2}\gamma \cap \wh{f}^{-N_1}(T(\beta))\neq\empt,$$
which in turn implies 
$$\wh{f}^{2N_1+N_2}\gamma \cap \cl{C}_1 \supset \wh{f}^{2N_1+N_2}\gamma \cap T(\beta)\neq\empt.$$

We want to continue this process inductively, in order to prove the following: 
\begin{lema}   \label{lema.induc}
$\wh{f}^{(i+1)N_1 + iN_2}\gamma \cap \cl{C}_i \neq\empt \ \ \ \forall \, i\in\N.$
\end{lema}

This lemma clearly implies that $\max(\txt{pr}_1(\rho(\wh{f}))) > 0$. This in turn implies that int$(\rho(\wh{f}))\neq\empt$, and then by Theorem \ref{teo.frint} we have that there is a periodic point $q$ for $f$ such that pr$_1(\rho(q,\wh{f})) >0$. This will yield Theorem D. 

Let us explain the main ideas in Lemma \ref{lema.induc}. As we saw, we have $\wh{f}^{2N_1+N_2}\gamma \cap T(\beta)\neq\empt$. The first main step in the proof will be to show that the iterate $\wh{f}^{2N_1+N_2}\gamma$ is \textit{well positioned} with respect to the point $T(p)$ (cf. Definition \ref{def.bp}), which very roughly speaking means that $\wh{f}^{2N_1+N_2}\gamma$ is in some sense `below' $T(p)$. Also, the fact that $(f_t)_t$ is positively transverse to $\wh{\cl{F}}$ will imply that $\wh{f}^{2N_1+N_2+n}\gamma\cap T(\beta)\neq\empt$ for all $n\geq 0$, and then $\wh{f}^{2N_1+N_2}\gamma$ is `anchored' to $T(\beta)$. Therefore, the forward iterates of $T(p)$ by $\wh{f}$ will `push' the arc $\wh{f}^{2N_1+N_2}\gamma$ downwards, and we will show that 
$$\wh{f}^{2N_1+2N_2}\gamma\cap\wh{f}^{-N_1}(T^2(\beta))\neq\empt.$$
The second main step in the proof will be to show that this intersection is a \textit{good intersection}, meaning that $\wh{f}^{2N_1+2N_2}\gamma$ intersects $\wh{f}^{-N_1}(T^2(\beta))$ in a way such that $\wh{f}^{3N_1+2N_2}\gamma$ is well positioned with respect to $T^2(p)$. 

In this way, by induction we will obtain Lemma \ref{lema.induc}, which as we mentioned, implies Theorem D.

The precise definitions of being well positioned and having good intersection are given in Section \ref{sec.wpgi}. To introduce such concepts, some previous definitions are given in \ref{sec.plg}. The induction steps for the proof of Lemma \ref{lema.induc} are stated in Section \ref{sec.indstep}, and in \ref{sec.lemaind} it is shown how the lemma follows from them. Finally, the proof of the induction steps is carried out in sections \ref{sec.bp-ib} and \ref{sec.ib-bp}.

\subsubsection{The points $p_i$, the arcs $\lambda$, $\wt{\gamma}$, the leaves $\beta_i$ and the integer $N_2$.}        \label{sec.plg}

In this section we give some definitions that will be used through the rest of $\S$\ref{sec.cs}. 

\begin{figure}[h]        
\begin{center} 
\includegraphics{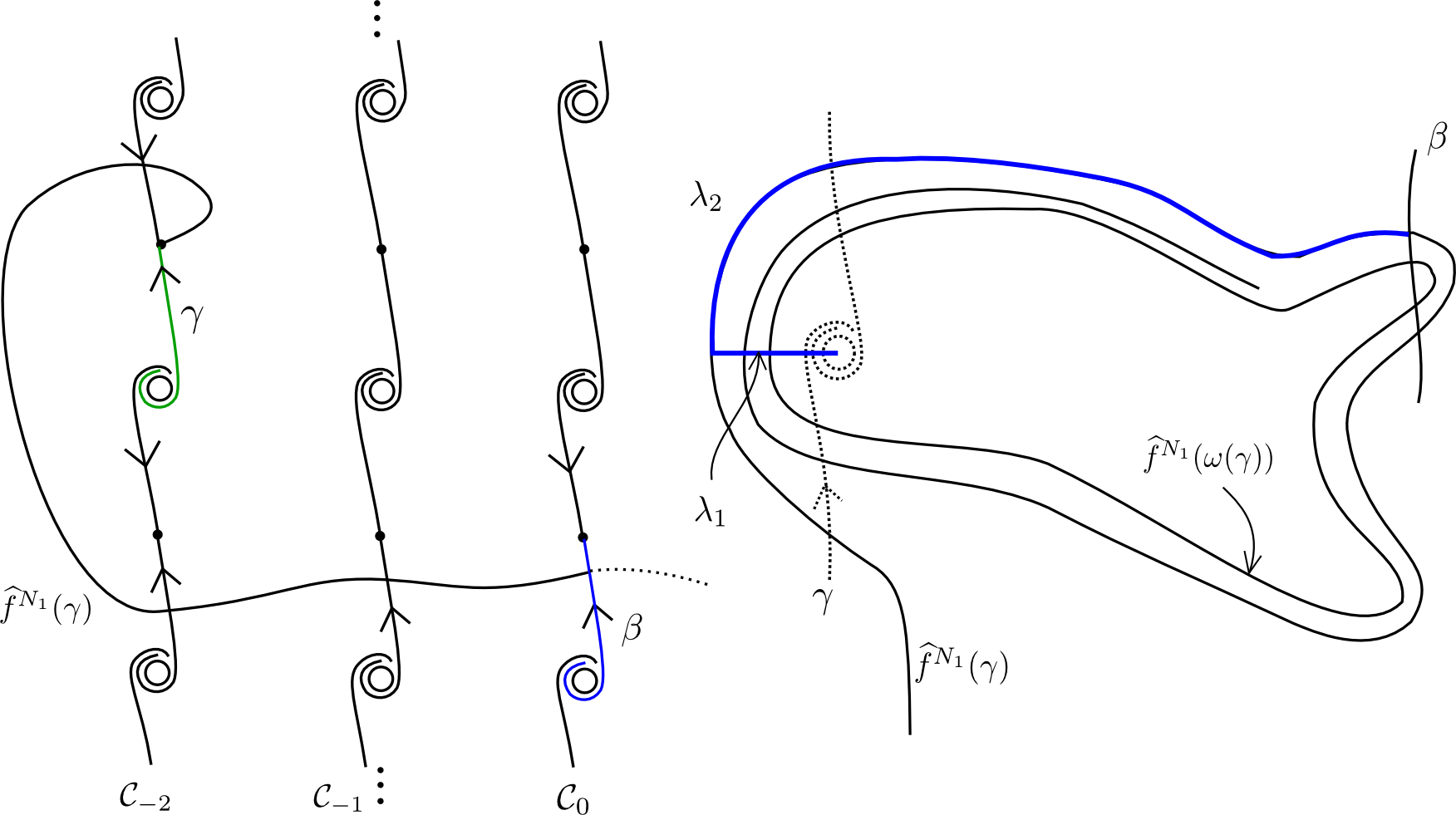}
\caption{Definition of $\beta$. Left: case that $\omega(\gamma)$ consists of a singularity. Right: case that neither $\omega(\gamma)$ nor $\alpha(\gamma)$ consist of a singularity.}
\label{fig.def-beta}
\end{center}  
\end{figure}

\noindent \textbf{The arc $\lambda$ and the leaf $\beta$.} Recall that $\wh{f}^{N_1}(\cl{C}_{i})\cap \cl{C}_{i+2}\neq\empt$ for all $i$ (see (\ref{eq.a.1})), and in particular $\wh{f}^{N_1}(\cl{C}_{-2})\cap \cl{C}_0\neq\empt$. Let $\gamma$ be a leaf of $\wh{\cl{F}}$ contained in $\cl{C}_{-2}$ such that 
$$\wh{f}^{N_1}(\gamma)\cap \cl{C}_0 \neq\empt.$$
We will now chose a leaf $\beta$ of $\wh{\cl{F}}$ contained in $\cl{C}_0$ such that $\wh{f}^{N_1}(\gamma)\cap\beta\neq\empt$, and an arc $\lambda$ going from sing$(\gamma)$ to $\wh{f}^{N_1}(\gamma)\cap\beta$. Such choice will be done separately for two possibilities on $\gamma$:
\vspace{2mm}\\
\textit{Case 1: Either $\alpha(\gamma)$ or $\omega(\gamma)$ consist of a singularity.} \\
If $\omega(\gamma)$ consists of a singularity (see Fig. \ref{fig.def-beta}), define $\beta$ as the leaf that contains the last point of intersection of $\wh{f}^{N_1}\gamma$ with $\cl{C}_0$, that is, $\beta$ is the leaf of $\wh{\cl{F}}$ that contains $\gamma(t_*)$, where
$$t_* = \max \{ t\in(0,1) \, : \, \wh{f}^{N_1}\gamma(t)\in\cl{C}_0 \}.$$ 
Define then the arc
$$\lambda = (\wh{f}^{N_1}\gamma|_{[t_*,1]})^{-1}.$$ 
Analogously, if $\omega(\gamma)$ does not consist of a singularity but $\alpha(\gamma)$ does, define $\beta$ as the first point of intersection of $\wh{f}^{N_1}(\gamma)$ with $\cl{C}_0$; that is, let $\beta$ be the leaf of $\wh{\cl{F}}$ that contains $\gamma(t_*)$, where
$$t_* = \min \{ t\in(0,1) \, : \, \wh{f}^{N_1}\gamma(t)\in\cl{C}_0 \},$$
and define
$$ \lambda = \wh{f}^{N_1}\gamma|_{[0,t_*]}.$$
Note that, in both cases, $\lambda(0)\in\txt{sing}(\gamma)$, $\lambda|_{[0,1)}\subset L(\cl{C}_0)$ and $\lambda(1)\in\beta\subset\cl{C}_0$. \\
\\
\textit{Case 2: Neither $\alpha(\gamma)$ or $\omega(\gamma)$ consist of a singularity.}\\
Note that by Corollary \ref{coro.bounded} we have that $\alpha(\gamma)$ and $\omega(\gamma)$ are compact, and then by Theorem \ref{poinc.ben} we have in this case that $\alpha(\gamma)$ and $\omega(\gamma)$ are either loops or generalized cycles of connections. 

Thus, there may not exist neither a first or a last point of intersection of $\wh{f}^{N_1}(\gamma)$ with $\cl{C}_0$ (see Fig. \ref{fig.def-beta}). We will define the auxiliary arc $\lambda$, and then we will define $\beta$ as the leaf of $\wh{\cl{F}}$ containing $\lambda(1)$. 

Let $s$ be any singularity contained in Fill$(\omega(\gamma))$, such that there is no other singularity $s'\in\txt{Fill}(\omega(\gamma))$ with $\txt{pr}_2(s')=\txt{pr}_2(s)$ and $\txt{pr}_1(s')<\txt{pr}_1(s)$. Let $\lambda_1$ be the straight horizontal arc going leftwards from $s$ to a point of $\wh{f}^{N_1}(\gamma)$ sufficiently close to $\wh{f}^{N_1}(\omega(\gamma))$ so that $\wh{f}^{-N_1}(\lambda_1)\cap\cl{C}_0=\empt$. Let $\lambda_2$ be a subarc of $\wh{f}^{N_1}\gamma$ going from $\lambda_1(1)$ to a point $z\in\wh{f}^{N_1}\gamma\cap\cl{C}_0$ such that $\lambda_2|_{[0,1)}\subset L(\cl{C}_0)$.

Let then $\lambda=\lambda_1\cdot\lambda_2$, and define $\beta$ as the leaf of $\wh{\cl{F}}$ that contains $\lambda(1)$.\\ 

\noindent \textbf{The points $p_i$ and the leaves $\beta_i$.} Recall from (\ref{eq.a.3}) and (\ref{eq.ccc}) that there is $p\in\R^2$ such that $\pi(p)\in\T^2$ is periodic and rotates downwards, $l_1,l_2$ are straight vertical lines oriented upwards such that $\cl{C}_0\subset L(l_1)$ and $R(l_1) \subset \{\wh{f}_t(p)\, : \, t\in\R\} \subset L(l_2)$, and $n_1\in\N$ is such that $\cl{C}_1= T_1^{n_1}(\cl{C}_0)\subset R(l_2)$. Fix an integer translate $p_0$ of $p$ of the form $p_0=T_2^i T_1^{-n_1} \in R(\cl{C}_{-1})\cap L(\cl{C}_0)$, with $i$ such that:
\begin{itemize}
\item $p_0$ is above $\lambda$, and 
\item $\{\wh{f}_t(p_0) \, : \, t\in (-\infty,0] \}$ is above $\{ \wh{f}_t(\beta) \, : \, t\in[-N_1,0] \}$.
\end{itemize}
Note that the second condition on $p_0$ may indeed be satisfised, because $\pi(p_0)\in\T^2$ is periodic and rotates downwards.

Define then
$$\ell_1 = T_1^{-n_1} l_1, \ \ \ \ \ \ell_2= T_1^{-n_1} l_2,$$
and note that 
\begin{equation}   \label{eq.78}
\{\wh{f}_t(p_0) \, : \, t\in\R \} \subset (\ell_1,\ell_2)\subset R(\cl{C}_{-1})\cap L(\cl{C}_0).
\end{equation}

Section \ref{sec.delta} will be devoted to the proof of the following.

\begin{figure}[h]        
\begin{center} 
\includegraphics{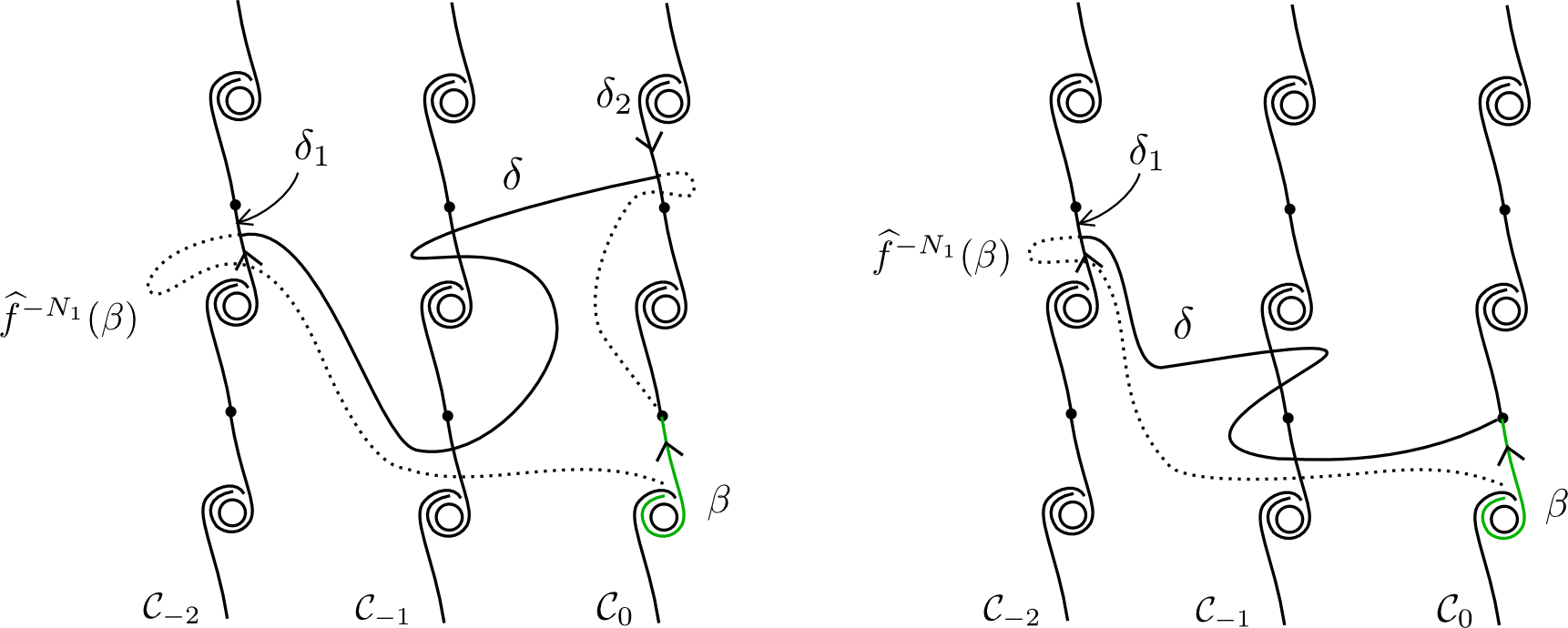}
\caption{The arc $\delta\subset\wh{f}^{-N_1}(\beta)$. Left: case that $\delta(1)$ belongs to a leaf $\delta_2\subset\cl{C}_0$. Right: case that $\delta(1)$ is a singularity.}
\label{fig.delta}
\end{center}  
\end{figure}

\begin{lema}         \label{lema.delta1}
There exists an arc $\delta$ such that (see Fig. \ref{fig.delta}):
\begin{enumerate}
\item $\txt{int}(\delta) \subset \wh{f}^{-N_1}\beta$,
\item $\delta(0)\in\cl{C}_{-2}$, $\delta(1)\in\cl{C}_{0}$, $\txt{int}({\delta})\subset R(\cl{C}_{-2})\cap L(\cl{C}_{0})$,
\item $\delta$ leaves a leaf $\delta_1\subset \cl{C}_{-2}$ on $t=0$ by the right (cf. Definition \ref{def.btl}), 
\item one of the following holds:
	\begin{enumerate}
	\item $\delta(1) \in \txt{sing}(\wh{\cl{F}})\cap\cl{C}_0$. 
	\item $\delta$ arrives in a leaf $\delta_2\subset \cl{C}_{0}$ on $t=1$ by the right.
	\end{enumerate}  
\end{enumerate}
\end{lema}


\begin{notacion}     \label{notacion1}
For not having to continuously refer to $\delta$ separately for cases 4(a) and 4(b) from last lemma, we define $\wt{\delta}_2=\empt$ in case item 4(a) holds, and $\wt{\delta}_2=\delta_2$ in case item 4(b) holds. 
\end{notacion}

Define $n_2>0$ such that 
\begin{equation}  \label{eq.1a}
\wh{f}^{-N_1}\lambda \cup \lambda \cup \{ \wh{f}_t(\beta) \, : \, t\in[-N_1,0] \}  \ \ \txt{is above} \ \  T_1^{n_1}T_2^{-n_2} ( \{ \wh{f}_t(\beta) \, : \, t\in[-N_1,0] \} \cup \delta_1\cup \wt{\delta}_2 ).
\end{equation} 

\begin{figure}[h]        
\begin{center} 
\includegraphics{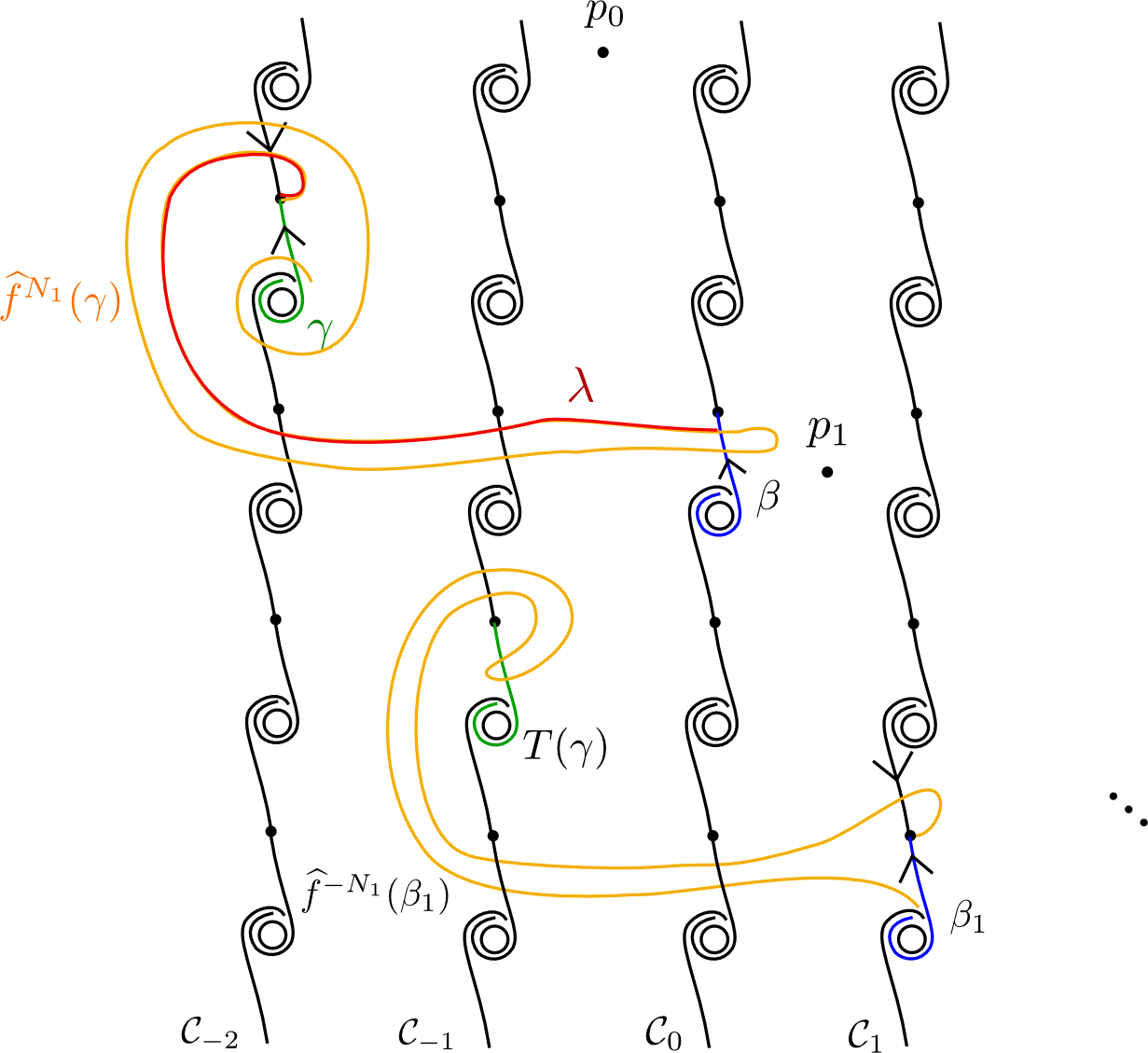}
\caption{$\wh{f}^{-N_1}\beta_1$ is below $\lambda$, and $p_0$ is above $\lambda$.}
\label{fig.T}
\end{center}  
\end{figure}

Denote 
\begin{equation}      \label{eq.20a}
T=T_1^{n_1}T_2^{-n_2},
\end{equation}
and for $i\geq 0$ let 
\begin{equation}   \label{def.pi}
\beta_i= T^i(\beta) \ \ \ \txt{and} \ \ \ p_i=T^i(p_0),
\end{equation}
(see Fig. \ref{fig.T}). Note that by (\ref{eq.78}), the paths $\{\wh{f}_t(p_i)\}$ are disjoint.\\

\noindent \textbf{The integer $N_2$.} The definition of $N_2$ is given by the following.

\begin{claim}    \label{claim.N2}
There exists $N_2>0$ such that:
\begin{enumerate}
\item $\wh{f}^{N_2}(p_0)$ is below $\wh{f}^{-N_1}(\beta_1)$, and 
\item $\wh{f}^{-k(N_1+N_2)}(p_k)$ is above $p_0$ for all $k\in\N$.
\end{enumerate}
\end{claim}
\begin{prueba}
Denote by $\pi:\R^2\ra\T^2$ the canonical projection. As $\pi(p_0)\in\T^2$ is periodic and rotates downwards, and as $\wh{f}^{-N_1}(\beta_1)$ is bounded (because by $\beta$ is bounded by Corollary \ref{coro.bounded}) the existence of $N_2^1>0$ satisfying Item 1 is trivial. Now we see that there is $N_2^2>0$ satisfying Item 2. Taking $N_2=\max\{N_1^1,N_2^2\}$, this will prove the claim. 

Also by the fact that $\pi(p_0)\in\T^2$ is periodic and rotates downwards, there is $N>0$ such that $\wh{f}^{-n}(p_{1})$ is above $p_0$, for all $n\geq N$. By induction this gives us that $\wh{f}^{-kn}(p_k)$ is above $p_0$ for all $n\geq N$ and $k\in\N$. Choosing $N_2^2=N$, we get in particular that $\wh{f}^{-k(N_1+N_2^2)}(p_k)$ is above $p_0$ for all $k\in\N$. That is, $N_2^2$ satisfies Item 2, as desired.  
\end{prueba}

\begin{remark}
Taking $N_2$ as a multiple of the period of $\pi(p_0)$ for $f$, we may assume that pr$_1(\wh{f}^{N_2}(p_i)) = \txt{pr}_1(p_i)$, for all $i$. 
\end{remark}

\noindent \textbf{The arc $\wt{\gamma}$.} By last, we define $\tl{\gamma}$ as 
$$\wt{\gamma} = \wh{f}^{-N_1}\lambda.$$
Observe that, as $\wh{f}^{N_1}(\wt{\gamma})=\lambda$, by the definition of $p_0$,
$$\wh{f}^{N_1}\tl{\gamma}  \ \ \   \txt{is below $p_0$,}$$
and by (\ref{eq.1a})
\begin{equation}      \label{eq.21a}
\wh{f}^{N_1}\tl{\gamma} \ \ \  \txt{is above $\wh{f}^{-N_1}(T(\beta))$}.
\end{equation}
Note also that for all $n$, $\wh{f}^{n}\tl{\gamma}(0)=\wh{f}^{n-N_1}\lambda(0)=\lambda(0) \in \txt{sing}(\wh{\cl{F}}) \cap \cl{C}_{-2}$, and $\wh{f}^{N_1}\tl{\gamma}(1)=\lambda(1)\in\beta\subset\cl{C}_0$.

\subsubsection{`Well positioned' and `good intersection'.}     \label{sec.wpgi}

Recall that for a bounded leaf $\ell\in\wh{\cl{F}}$, $\txt{sing}(\ell)$ denotes the set of singularities of $\wh{\cl{F}}$ contained in $\txt{Fill}(\alpha(\ell)) \cup \txt{Fill}(\omega(\ell))$. Also, recall that by Corollary \ref{coro.bounded}, every leaf $\ell$ of $\wh{\cl{F}}$ contained in $\cup_i \cl{C}_i$ is bounded, and therefore for such $\ell$, sing$(\ell)\neq\empt$.

Let $\delta$ be the arc given by Lemma \ref{lema.delta1}. For each $n\in \N_0$ we now define 
$$F_n=\left\{ \wh{f}^{-k(N_1+N_2)}(p_{n+k}) \, : \,  k\in\N_0 \right\}  \ \bigcup \   \cup_{k\in\N_0} \, \txt{sing}(\beta_{n+k})  \    \bigcup \  \cup_{k\in\N_0} \, T^k \left( \txt{sing}(\delta_1)\cup \txt{sing}(\wt{\delta}_2) \right),$$
and proceed to the definitions of `well positioned' and `good intersection'. Consider the arc $\wt{\gamma}$ defined in $\S$\ref{sec.plg}.

\begin{figure}[h]        
\begin{center} 
\includegraphics{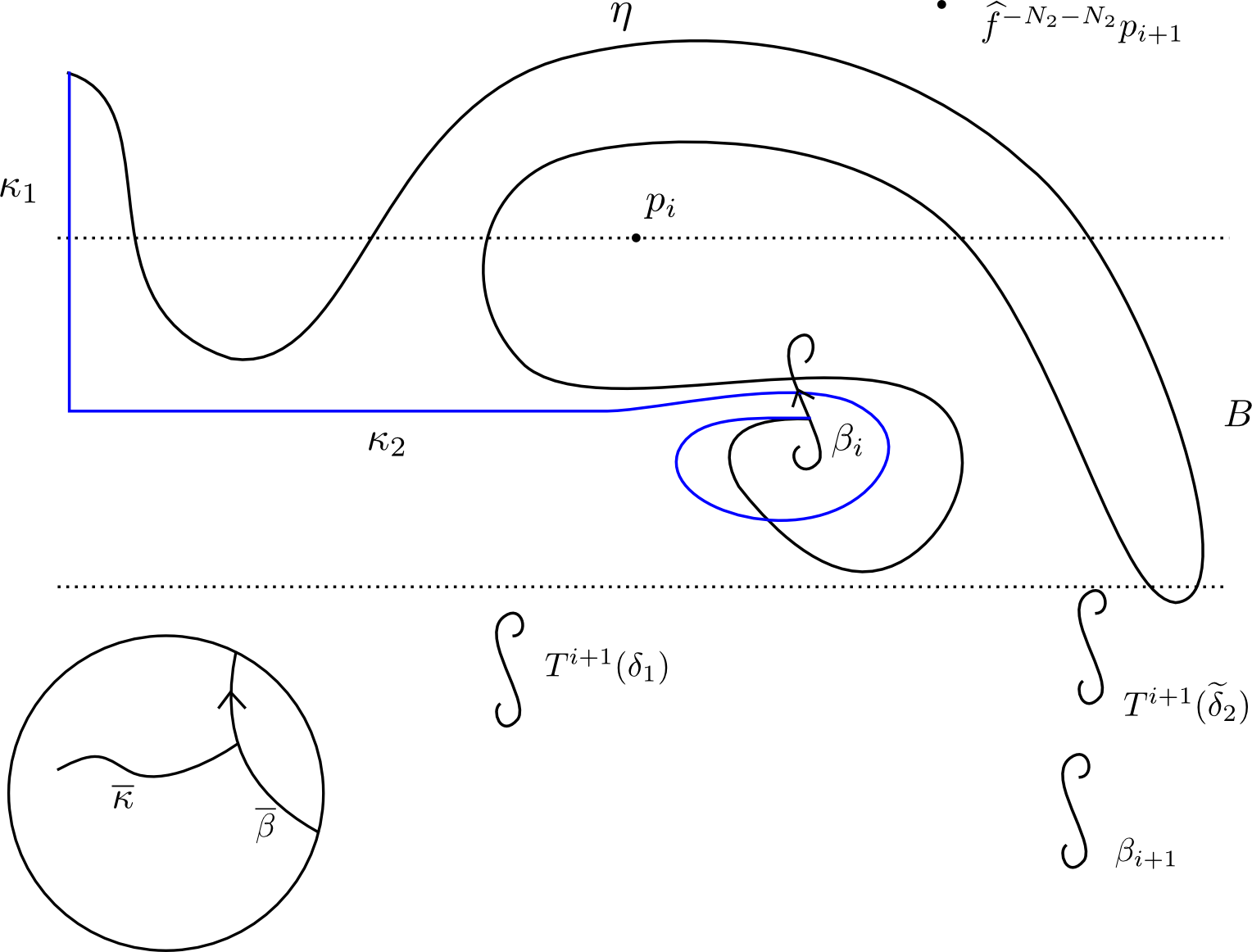}
\caption{Illustration of an arc which is well positioned with respect to $p_i$. Left-bottom corner: the lifts $\ol{\beta},\ol{\kappa}$ to the universal cover of $\R^2\minus (\txt{sing}(\wh{\cl{F}})\cap F_i)$.}
\label{fig.bp}
\end{center}  
\end{figure}

\begin{definicion}         \label{def.bp}
Let $i\in\N_0$. We say that an iterate $\wh{f}^n(\tl{\gamma})$ of $\tl{\gamma}$ is \textit{well positioned with respect to} $p_i$ if $\wh{f}^n(\tl{\gamma})$ contains a subarc $\eta$ such that (see Fig. \ref{fig.bp}):
\begin{enumerate}
\item $\eta(0)=\wh{f}^n(\tl{\gamma}(0))=\tl{\gamma}(0)\in\txt{sing}(\wh{\cl{F}}) \cap \cl{C}_{-2}$,
\item $\eta(1)\in \beta_{i}$, and
\item $\eta$ is homotopic wfe Rel$(F_i)$ to an arc $\kappa$ of the form $\kappa=\kappa_1\cdot\kappa_2$, where:
\begin{enumerate}
\item $\kappa_1$ is a vertical arc from $\eta(0)$ to a point in the strip 
$$B_i= \{ x\in\R^2 \, : \,  x \txt{ is below $p_i$ and above $\beta_{i+1} \cup T^{i+1}(\delta_1\cup \wt{\delta}_2)$} \},$$ and
\item $\kappa_2$ is an arc contained in $B\minus F_i$ from $\kappa_1(1)$ to $\eta(1)\in\beta_i$ such that, if $\ol{\kappa}_2$ is a lift of $\kappa_2$ to the universal cover of $\R^2\minus (\txt{sing}(\wh{\cl{F}})\cap F_i)$, and if $\ol{\beta}$ is a lift of $\beta_i$ containing $\ol{\kappa}_2(1)$, then $\ol{\kappa}_2(0) \in L(\ol{\beta})$.
\end{enumerate}
\end{enumerate}
\end{definicion}

As we mentioned in $\S$\ref{sec.ideasing}, last definition gives us that $\wh{f}^n\wt{\gamma}$ contains an arc $\eta$ which is in some sense ``below'' $p_i$. Also, $\eta$ intersects $\beta_i$ ``by the left''.

Consider the arc $\delta$ from Lemma \ref{lema.delta1}.

\begin{figure}[h]        
\begin{center} 
\includegraphics{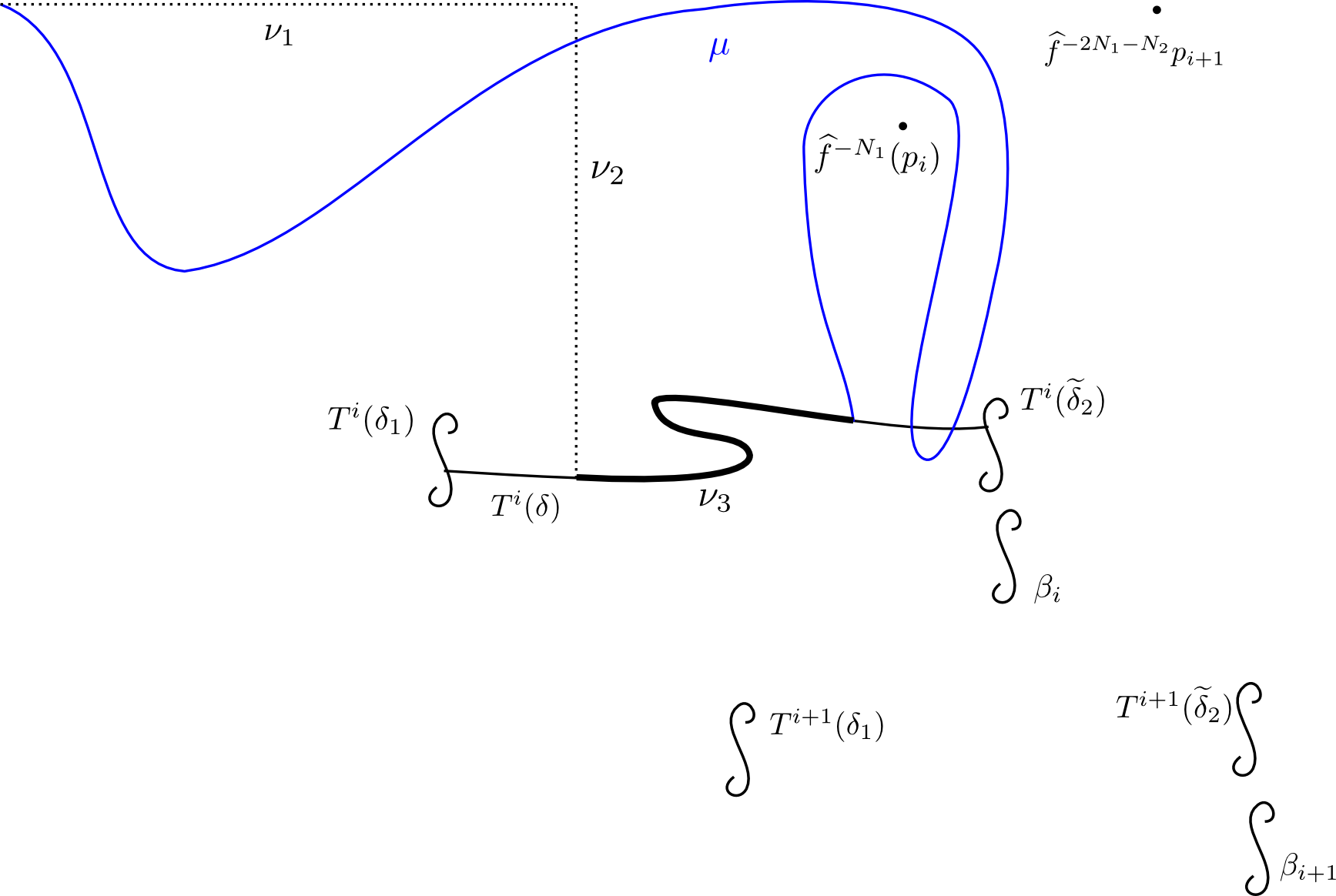}
\caption{Illustration of an arc which has good intersection with $\wh{f}^{-N_1}\beta_i$.}
\label{fig.ib}
\end{center}  
\end{figure}

\begin{definicion}         \label{def.ib}
Let $i\in\N$. We say that an iterate $\wh{f}^n(\tl{\gamma})$ of $\tl{\gamma}$ has \textit{good intersection with} $\wh{f}^{-N_1}(\beta_i)$ if $\wh{f}^n(\tl{\gamma})$ contains a subarc $\mu$ such that (see Fig. \ref{fig.ib}):
\begin{enumerate}
\item $\mu(0)=\wh{f}^n(\tl{\gamma}(0))=\tl{\gamma}(0) \in\txt{sing}(\wh{\cl{F}})\cap \cl{C}_{-2}$,
\item $\mu(1)\in T^i(\delta) \subset \wh{f}^{-N_1}(\beta_{i})$, and
\item $\mu$ is homotopic wfe Rel$(\wh{f}^{-N_1}(F_i))$ to an arc $\nu$ of the form $\nu=\nu_1\cdot\nu_2\cdot\nu_3$, where $\nu_1$ is a horizontal arc starting in $\mu(0)$ and ending in the `strip' $R(\cl{C}_{i-2}) \cap L(\cl{C}_{i-1})$, $\nu_2$ is a vertical arc from $\nu_1(1)$ downwards to a point of $T^i(\delta)$, and $\nu_3$ is the arc contained in $T^i(\delta)$ going from $\nu_2(1)$ to $\mu(1)$.
\end{enumerate}
\end{definicion}

\subsubsection{Induction steps.}     \label{sec.indstep}

The following two lemmas are the main induction steps in the proof of Theorem D, for the case that $\cl{C}$ has singularities. 

\begin{lema}   \label{bp-ib}
Let $i\in\N_0$. If an iterate $\wh{f}^n\tl{\gamma}$ of $\tl{\gamma}$ is well positioned with respect to $p_i$, then $\wh{f}^{n+N_2}\tl{\gamma}$ has good intersection with $\wh{f}^{-N_1}\beta_{i+1}$.
\end{lema}

\begin{lema}    \label{ib-bp}
Let $i\in\N$. If an iterate $\wh{f}^n\tl{\gamma}$ of $\tl{\gamma}$ has good intersection with $\wh{f}^{-N_1}\beta_{i}$, then $\wh{f}^{n+N_1}\tl{\gamma}$ is well positioned with respect to $p_i$. 
\end{lema}

\subsubsection{Lemmas \ref{bp-ib} and \ref{ib-bp} imply Theorem D.}    \label{sec.lemaind}

As we mentioned in $\S$\ref{sec.ideasing}, to prove Theorem D it suffices to prove Lemma \ref{lema.induc}. To prove such lemma from the induction lemmas \ref{bp-ib} and \ref{ib-bp}, it suffices to show that $\wh{f}^{N_1}(\wt{\gamma})$ is well positioned with respect to $p_0$.

Consider the arc $\lambda$ defined in $\S$\ref{sec.plg}, and recall that by definition, $\lambda(0)\in \txt{sing}(\gamma)$, $\lambda|_{[0,1)}\subset L(\cl{C}_0)$, and $\lambda(1)\in\beta_0\subset\cl{C}_0$. By definition $\wt{\gamma}=\wh{f}^{-N_1}\lambda$, and then $\wh{f}^{N_1}(\wt{\gamma})=\lambda$. Defining $\eta=\wh{f}^{N_1}(\wt{\gamma})=\lambda$ we have that $\eta$ satisfies items 1 and 2 from the definition of being well positioned with respect to $p_0$. 

We now see that $\eta$ satisfies item 3. By definition, $p_0$ is above $\lambda=\eta$. By item 2 from Claim \ref{claim.N2} we have that $\wh{f}^{-i(N_1+N_2)}p_i$ is above $p_0$ for all $i>0$, and then $\eta$ is below $\wh{f}^{-i(N_1+N_2)}p_i$ for all $i$. Also, from (\ref{eq.1a}) we have that $\lambda=\eta$ is above the sets sing$(\beta_1)$, sing$(T(\delta_1))$, sing$(T(\wt{\delta}_2))$, and then $\eta$ is above sing$(\beta_i)$, sing$(T^i(\delta_1))$, sing$(T^i(\wt{\delta}_2))$ for all $i\geq 1$. This easily implies that $\eta$ is homotopic wfe Rel($F_0$) to an arc $\kappa=\kappa_1\cdot\kappa_2$, with $\kappa_1$ straight vertical and $\kappa_2$ contained in the strip $B_0$ from Definition \ref{def.bp}. We have then item 3(a) from such definition. 

To verify item 3(b), consider a lift $\ol{\kappa}_2$ of $\kappa_2$ to the universal cover of $\R^2\minus (\txt{sing}(\wh{\cl{F}})\cap F_0)$, and let $\ol{\beta}$ be a lift of $\beta$ containing $\ol{\kappa}_2(1)$. We want to see that $\ol{\kappa}_2(0)\in L(\ol{\beta})$. Let $\ol{\kappa}_1$ be the lift of $\kappa_1$ such that $\ol{\kappa}_1(1)=\ol{\kappa}_2(0)$ and note that as $\kappa(0)\in\txt{sing}(\wh{\cl{F}})\minus F_0$, then $\ol{\kappa}_1(0)$ is a singularity of the lift of the foliation $\wh{\cl{F}}\minus (\txt{sing}(\wh{\cl{F}})\cap F_0)$. Suppose by contradiction that $\ol{\kappa}_2(0)\in R(\ol{\beta})$. As $\kappa_1\cap\beta=\empt$ we then have $\ol{\kappa}_1(0)=\ol{\kappa}_1\cdot\ol{\kappa}_2(0)\in R(\ol{\beta})$. Let $(\ol{f}_t)$ be the canonical lift of the isotopy $(\wh{f}_t)$. By Lemma \ref{lemabrouwer1} we have that $\ol{f}^{-N_1}(\ol{\kappa}_1\cdot\ol{\kappa}_2)\cap\ol{\beta}\neq\empt$, and then $\wh{f}^{-N_1}(\kappa_1\cdot\kappa_2)\cap\beta\neq\empt$. As $\eta$ is homotopic wfe Rel($F_0$) to $\kappa_1\cdot\kappa_2$, we get that there is a lift $\ol{\eta}$ of $\eta$ such that $\ol{\eta}(0) = \ol{\kappa}_1\cdot\ol{\kappa}_2(0)$ and $\ol{\eta}(1)=\ol{\kappa}_1\cdot \ol{\kappa}_2(1)$, and then $\wh{f}^{-N_1}\ol{\eta}\cap\ol{\beta} \neq\empt$. Then, 
$$\wh{f}^{-N_1}\lambda\cap\beta=\wh{f}^{-N_1}\eta\cap\beta\neq\empt,$$
which contradicts the definition of $\lambda$. Therefore we must have $\ol{\kappa}_1\cdot\ol{\kappa}_2(0)\in L(\ol{\beta})$, and then $\ol{\kappa}(0)\in L(\ol{\beta})$ because $\kappa_1\cap\beta=\empt$, as desired.

\subsubsection{Proof of lemma \ref{bp-ib}}     \label{sec.bp-ib}

To simplify the notation, we will prove Lemma \ref{bp-ib} for the case that $i=0$. Namely, we will prove the following.

\begin{figure}[h]        
\begin{center} 
\includegraphics{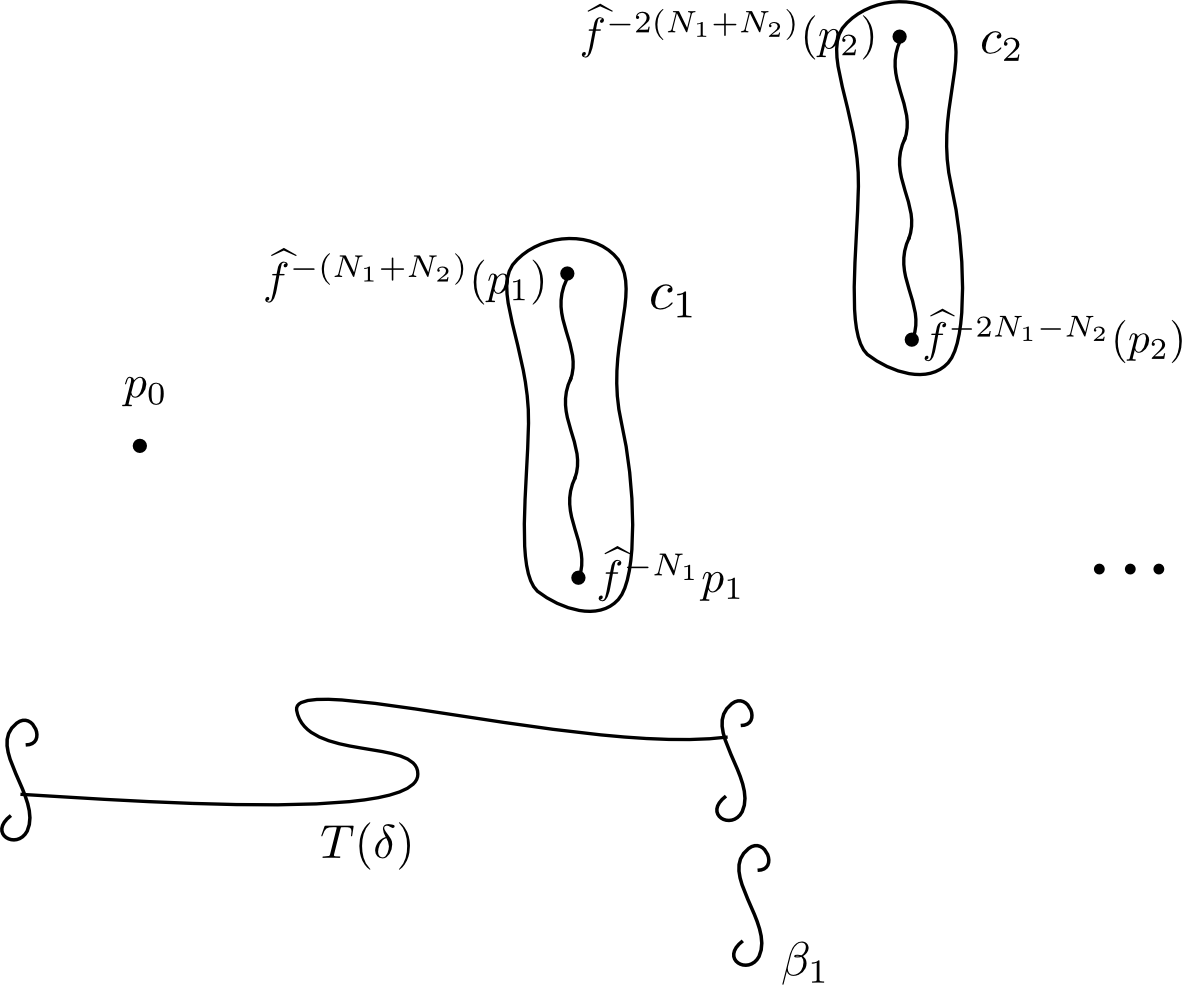}
\caption{The curves $c_i$.}
\label{fig.ci}
\end{center}  
\end{figure}

\begin{lema}     \label{bp-ib1}
If an iterate $\wh{f}^n\tl{\gamma}$ of $\tl{\gamma}$ is well positioned with respect to $p_0$, then $\wh{f}^{n+N_2}\tl{\gamma}$ has good intersection with $\wh{f}^{-N_1}\beta_{1}$.
\end{lema}
In the proof of this lemma, the fact that $i=0$ will have nothing special, and therefore the proof of Lemma \ref{bp-ib} will follow by an identical argument. The rest of this section is devoted to the proof of Lemma \ref{bp-ib1}.\\
\\
\textbf{The loops $c_i$, $i\geq 1$.} We now fix a family of loops $c_i$ that will be used in the following lemma. Recall that the points $p_i$ were defined in (\ref{def.pi}) in a way that the paths $\{\wh{f}_t(p_i) \, : \, t\in\R \}$ are pairwise disjoint. Therefore we can choose simple closed curves $c_i$, for $i\in\N$, such that (see Fig. \ref{fig.ci}):
\begin{itemize}
\item $\{ \wh{f}_t( \wh{f}^{-i(N_1+N_2)}(p_{i})) \, : \, t\in[0,N_2] \} \subset \txt{int}(c_{i})$, where int$(c_i)$ denotes the bounded connected component of $\R^2\minus c_i$,
\item $\txt{Fill}(c_i) \cap \wh{f}_t(F_0) = \{\wh{f}_t(\wh{f}^{-i(N_1+N_2)}(p_i)) \}$ for all $t\in [0,N_2]$,
\item $\txt{Fill}(c_i) \subset R(\cl{C}_{i-1})\cap L(\cl{C}_{i})$.
\end{itemize}
Also, note that by definition of $p_0$, $\{\wh{f}_t(p_0)\, : \, t\in(-\infty,0]\}$ is above $\{ \wh{f}_t \beta \, : \, t\in [-N_1,0]\}$. By periodicity and as $T(\delta)\subset \wh{f}_{-N_1}\beta_1$, we may choose $c_1$ such that
\begin{itemize}
\item $\txt{Fill}(c_1)$ is above $T(\delta)$. 
\end{itemize}

For the next lemma, we also define $E= (F_0 \minus \{p_0\} ) \cup\{\tl{\gamma}(0) \}$, and for $t\in\R$ set
$$E^t = \wh{f}_t(E).$$

\begin{lema}      \label{lema.fam-Pit}
There exist a continuous family $\{\Pi_t\}_{t\in[0,N_2]}$ of covering maps $\Pi_t:\D \ra \R^2\minus E^t$, and an isotopy $( \ol{f}_t )_{t\in[0,N_2]}$ of $\D$, starting in $\ol{f}_0=\txt{Id}$, such that:
\begin{enumerate}

\item For any $0\leq t \leq N_2$, the diagram             \label{ita2}
\begin{diagram}    \label{diag.1}
\D              &  \rTo^{\ol{f}_t}                    & \D  \\
\dTo^{\Pi_0}    &                                     &  \dTo_{\Pi_t}  \\
\R^2\minus E    & \rTo_{\wh{f}_t|_{\R^2\minus E}}     & \R^2\minus E^t
\end{diagram}
commutes.

\item For all $t\in[0,N_2]$ we have $\Pi_0|_{\D\minus C} = \Pi_t|_{\D\minus C}$, where $C= \Pi_0^{-1} \left( \cup_{k\in\N} \, \txt{int}_E(c_k) \right)$ and $\txt{int}_E(c_k)$ denotes the bounded connected component of $(\R^2\minus E)\minus c_k$.       \label{ita1}

\item Let $*\in \{ \delta_1, \wt{\delta}_2,\beta \}$. Then, if $\ol{*}$ is a lift of $*$ to $\D$ by $\Pi_0$, $\ol{*}$ is a Brouwer curve for $\ol{f}_{N_2}$.     \label{ita3}

\end{enumerate}  
\end{lema}
\begin{prueba}
Let $P=\{ \wh{f}^{-k(N_1+N_2)}(p_{k}) \, : \, k\in\N\}$. Define
$$E_1= E \minus P = E\cap \txt{sing}(\wh{\cl{F}}),$$
and let $\pi_1: \D \ra \R^2 \minus E_1$ be a covering map. As $E_1 \subset \txt{sing}(\wh{\cl{F}})$, we may consider the canonical lift $(\wt{f}_t)$ by $\pi_1$ of the isotopy $\wh{f}_t|_{\R^2\minus E_1}$. Also, the Brouwer foliation $\wh{\cl{F}}|_{\R^2\minus E_1}$ (with singularities) lifts by $\pi_1$ to a foliation of $\D$ (also with singularities) which is transverse to the isotopy $\wt{f}_t$. 
For every $t$, the following diagram commutes:
\begin{diagram}    
\D              &  \rTo^{\wt{f}_t}                    & \D  \\
\dTo^{\pi_1}    &                                     &  \dTo_{\pi_1}  \\
\R^2\minus E_1  & \rTo_{\wh{f_t}|_{\R^2\minus E_1}}  & \R^2\minus E_1
\end{diagram}
For $t\in[0,N_2]$, let $P^t = \wh{f}_t(P)$ and define
$$D^t= \D\minus \pi_1^{-1}(P^t).$$
As $E=E_1\cup P$, considering the restrictions $\wh{f}_t|_{\R^2\minus E}$, $\pi_1|_{D^0}$ and $\tl{f}_t|_{D^0}$, we obtain the following commutative diagram
\begin{diagram}    
D^0                   &  \rTo^{\wt{f}_t|_{D^0 }}              & D^t \\
\dTo^{\pi_1|_{D^0} }         &                                       &  \dTo_{\pi_1|_{D^t }}  \\
\R^2\minus E          & \rTo_{\wh{f}_t|_{\R^2\minus E}}      & \R^2\minus E^t
\end{diagram}
for $t\in[0,N_2]$.

Using standard techniques from plane topology, one may construct an isotopy 
$(J_t)_{t\in[0,N_2]}$ on $\R^2$ starting in the identity and such that 
\begin{itemize}
\item $J_t(x)=x$  for all $t\in[0,N_2]$ and all $x\in \R^2\minus \cup_{k\in\N} \, \txt{int}(c_{k})$, and
\item $J_t(p_{k})=\wh{f}_t(p_{k})$ for all $t\in[0,N_2]$ and all $k\in\N$.
\end{itemize}

Consider the restriction $( J_t|_{\R^2\minus E_1} )_{t\in[0,N_2]}$, which is an isotopy on $\R^2\minus E_1$, and consider the canonical lift $( \wt{J}_t )_{t\in[0,N_2]}$ by $\pi_1$ of $J_t|_{\R^2\minus E_1}$. Defining
\begin{equation}     \label{eq.bp-ib1}
F_t=\wt{J}^{-1}_t \circ \wt{f}_t
\end{equation}
we obtain an isotopy $(F_t )_{t\in[0,N_2]}$ on $\D$. By the properties of $J_t$, $F_t(D^0)=D^0$ for all $t$, and then we have an isotopy $(F_t|_{D^0})_{t\in[0,N_2]}$ on $D^0$.

Consider a covering map $\pi_2:\D \ra D^0$, and the canonical lift $(\ol{f}_t)$ by $\pi_2$ of the isotopy $(F_t|_{D^0})$, so
\begin{equation}    \label{eq.bp-ib2}
\pi_2 \circ \ol{f}_t = F_t|_{D^0} \circ \pi_2.
\end{equation}
For $t\in[0,N_2]$, define the covering map $\pi_{2}^t :\D \ra D^t$ by
$$\pi_{2}^t = \wt{J}_t|_{D^0} \circ \pi_2.$$ 
As $J_0= \txt{Id}$ we have that $\pi_{2}^0=\pi_2$. The maps $\pi_{2}^t$ vary continuously, and we have 
\begin{align*}
 \pi_{2}^t \circ \ol{f}_t 
        &=  \wt{J}_t|_{D^0} \circ \pi_2 \circ \ol{f}_t = \wt{J}_t|_{D^0} \circ F_t|_{D^0}
              \circ \pi_2    \ \ \   \ \ \ \ \ \txt{by (\ref{eq.bp-ib2}}) \\
        &= \wt{J}_t|_{D^0} \circ \wt{J}_t^{-1}|_{D^t} \circ \wt{f}_t|_{D^0} \circ\pi_2  
				 \ \ \   \ \ \ \ \ \ \ \ \ \ \ \txt{by (\ref{eq.bp-ib1}})    \\
		    &= \txt{Id}_{D^t} \circ \wt{f}_t|_{D^0} \circ \pi_2 = \wt{f}_t|_{D^0} \circ \pi_2.   
\end{align*}

Thus, for every $t\in[0,N_2]$, the following diagram commutes
\begin{diagram}    
\D                              &  \rTo^{\ol{f}_t }                                & \D   \\
\dTo^{\pi_2}                    &                                                  & \dTo_{\pi_{2}^t}  \\
D^0     												& \rTo_{\tl{f}_t|_{D^0}}   												 & D^t
\end{diagram}

If we define $\Pi_t:\D\ra\R^2\minus E^t$ for $t\in[0,N_2]$ by 
$$\Pi_t= \pi_1|_{D^t} \circ \pi_{2}^t,$$
we have that the maps $\Pi_t$ vary continuously, and $\Pi_0=\pi_1|_{D^0}\circ \pi_2^0 = \pi_1|_{D^0} \circ\pi_2$. Therefore combining the last two commutative diagrams we obtain that 
$$\Pi_t \circ \ol{f}_t = \wh{f}_t|_{\R^2\minus E} \circ \Pi_0 \ \ \ \ \ \ \txt{for $t\in[0,N_2]$},$$
which proves Item 1 of the lemma\footnote{So far the intermediate covering $\pi_1$ might seem unnecesary. However, it will be used in the proof of item 3.}. 

To prove Item 2, first note that, as 
$$J_t(x) = x \ \ \ \ \ \forall \, t\in[0,N_2] \ \ \forall \, x\in \R^2\minus \cup_{k\in\N} \, \txt{int}(c_{k}),$$
we have that 
$$\wt{J}_t(x)=x \ \ \ \ \ \forall \, t\in[0,N_2] \ \ \forall \, x\in\D\minus C_1,$$ 
where $C_1= (\pi_1)^{-1} ( \cup_{k\in\N}\, \txt{int}(c_{k}) ) $. By this, if 
$$C= \pi_2^{-1}(C_1\cap D^0) = \Pi_0^{-1}( \cup_{k\in\N}\, \txt{int}_{E}(c_{k}) ),$$ 
we have that
\begin{equation}          \label{eq.2a}
\pi_{2}^t(x) = \wt{J}_t|_{D^0}\circ\pi_2 (x) = \pi_2(x) \ \ \ \ \ \forall \, t\in[0,N_2] \ \ \forall \, x \in \D\minus C.
\end{equation}
Observe also that $\D\minus C_1 \subset D^t$ for all $t\in[0,N_2]$, and in particular 
\begin{equation}            \label{eq.3a}
\pi_1|_{D^t} \circ \pi_2(x) = \pi_1|_{D^0} \circ \pi_2(x) \ \ \ \ \ \forall \, t\in[0,N_2] \ \ \forall \, x\in\D\minus C.
\end{equation}
Thus, for all $t\in[0,N_2]$ and all $x\in \D\minus C$ we have
\begin{align*}
\Pi_t(x) = \pi_1|_{D^t}\circ \pi_{2}^t(x) &= \pi_1|_{D^t} \circ \pi_2(x)  \ \ \ \ \ \ \ \ \txt{by } (\ref{eq.2a})   \\
                &= \pi_1|_{D^0} \circ \pi_2(x) \ \ \ \ \ \ \ \ \txt{by } (\ref{eq.3a})      \\
								&= \Pi_0(x). 
\end{align*}
This proves Item (2).

To prove Item \ref{ita3}, we will now see that, if $\ol{\beta}$ is a lift by $\Pi_0$ of $\beta$, then $\ol{\beta}$ is a Brouwer line for $\ol{f}_{N_2}$. The same argument will hold for $\delta_1$ and $\wt{\delta}_2$. 

As we saw above, the isotopy $(\wh{f}_t|_{\R^2\minus E_1})$ lifts by $\pi_1$ to its canonical lift $(\wt{f}_t)$ on $\D$, and the Brouwer foliation $\wh{\cl{F}}$ lifts to a foliation of $\D$ with singularities which is transverse to $\wt{f}_t$. In particular, if $\beta'$ is a lift by $\pi_1$ of $\beta$, then $\beta'$ is a Brouwer line for $\wt{f}_{N_2}$. Observe also that, as $\beta$ is disjoint from $P$, then $\pi_1^{-1}(\beta)=(\pi_1|_{D^0})^{-1}(\beta)$. 

Now, as $\beta'$ is a Brouwer line for $\wt{f}_{N_2}$, we have that if $\beta''$ is a lift of $\beta'$ by $\pi_2^0$, then $\ol{f}_{N_2}(\beta'')\cap\beta'' = \empt$. 

Thus, any lift $\beta''$ of $\beta$ by $\Pi_0=\pi_1|_{D^0}\circ\pi_2^0$ is disjoint from its image by $\ol{f}_{N_2}$, and to prove Item \ref{ita3} it suffices to show that $\ol{f}_{N_2}(\beta'')$ cannot be contained in $L(\beta'')$. However, this is a simple consequence of the fact that the maps $\pi_1|{D^0}$ and $\pi_2^t$ preserve orientation, and that $\beta'$ is a Brouwer line for $\wt{f}_{N_2}$. This proves item \ref{ita3}, and finishes the proof of the lemma.
\end{prueba}
\vspace{2mm}
\textbf{The arc $S$ and the loop $c_0$.} Recall our notation from Lemma \ref{lema.fam-Pit} for $E= (F_0 \minus \{p_0\})\cup \{ \wt{\gamma}(0) \}$. By construction, the path $( \wh{f}_t(p_0) )_{t \in [0, N_2]}$ is contained in the straight vertical strip $(\ell_1,\ell_2)\subset R(\cl{C}_{-1}) \cap L(\cl{C}_{0})$ (see (\ref{eq.78})). As $E$ is disjoint from $R(\cl{C}_{-1}) \cap L(\cl{C}_{0})$, we have that
\begin{equation}        \nonumber 
( \wh{f}_t(p_0) )_{t\in [0,N_2]} \ \ \txt{is homotopic wfe Rel$(E)$ to a straight vertical arc $S$}. 
\end{equation}

Let $c_0$ be a closed loop oriented counterclock-wise containing $\wt{\gamma}(0)$ in its interior, and such that $c_0\subset L(\cl{C}_{-1})$. Then, by definition of the sets $F_i$, Fill$(c_0)\cap F_i = \empt$ for all $i\geq 0$. 
\vspace{2mm} \\
\textbf{Lifting to $\D$ by $\Pi_0$}.
First recall that as we are assuming that $\wh{f}^n(\wt{\gamma})$ is well positioned with respect to $p_0$, there exists a subarc $\eta\subset\wh{f}^n(\wt{\gamma})$ which is homotopic wfe to an arc $\kappa=\kappa_1\cdot\kappa_2$ satisfying the properties from definition \ref{def.bp}. 

\begin{figure}[h]        
\begin{center} 
\includegraphics{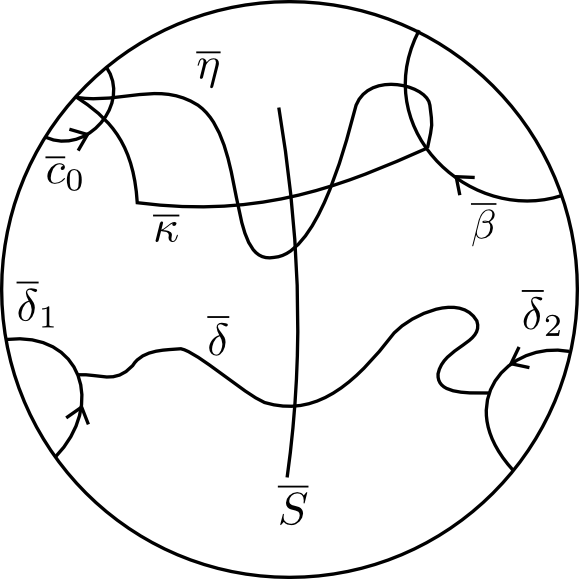}
\caption{Lifting by $\Pi_0$.}
\label{fig.levant}
\end{center}  
\end{figure}

Fix a lift $\ol{\kappa}$ of $\kappa|_{(0,1]}$ to $\D$ by $\Pi_0: \D\to \R^2\minus E$. Let $\ol{\beta}$ be a lift of $\beta$ to $\D$ by $\Pi_0$ that contains the point $\ol{\kappa}(1)$ (see Fig. \ref{fig.levant}). Let $q_3 =\min\{t\in(0,1] \, : \, \kappa(t) \in c_0\}$, and define $\ol{c}_0$ as a lift of $c_0$ such that $\ol{c}_0$ contains $\ol{\kappa}(q_3)$. Note that, as $c_0$ is oriented counterclock-wise, 
\begin{equation}        \label{eq.5a}
\ol{\kappa}(t)\in L(\ol{c}_0) \ \ \ \ \txt{ for} \ t\in(0,q_3). 
\end{equation}

Define then $\ol{\eta}$ as the lift of $\eta|_{(0,1]}$ by $\Pi_0$ such that $\ol{\eta}(1)=\ol{\kappa}(1)$, and note that, as $\eta$ and $\kappa$ are homotopic wfe Rel$(F_0)$, and by (\ref{eq.5a}),
\begin{equation}     \label{eq.8a}
\ol{\eta}(t)\in L(\ol{c}_0), \ \ \ \ \txt{ for $t$ small enough.}
\end{equation}
Observe that as $\wh{f}^{N_2}p_0$ is below $T(\delta)$ (see Claim \ref{claim.N2}) and by definition of $\kappa$ and $S$, the arc $S$ intersects $\kappa$ and $T(\delta)$.

Let 
\begin{equation}     \label{eq.10a}
q_4 = \min \{ t\in (0,1] \, : \kappa(t)\in S\},
\end{equation}
and let $q_5$ be such that $S(q_5) = \kappa(q_4).$ Let $\ol{S}$ be the lift of $S$ by $\Pi_0$ such that $\ol{S}(q_5)\in \ol{\kappa}(q_4).$ Note that, as $S$ is homotopic wfe Rel($E$) to $(\wh{f}_t(p_0))_{t\in[0,N_2]}$, we have that 
$$\ol{S}(1)=\ol{f}_{N_2}(\ol{S}(0)).$$

Let $q_6,t_*\in[0,1]$ be such that $S(q_6)=T\delta(t_*)$, and let $\ol{\delta}$ be the lift of $T(\delta)$ by $\Pi_0$ such that $\ol{S}(q_6) = \ol{\delta}(t_*)$. Also, let $\ol{\delta}_1$ be a lift of $\delta_1$ by $\Pi_0$ such that 
\begin{equation}   \label{eq.d111}
\ol{\delta}(0)\in \ol{\delta}_1. 
\end{equation}
In case $\wt{\delta}_2\neq\empt$, let $\ol{\delta}_2$ be a lift of $\wt{\delta}_2$ by $\Pi_0$ such that 
\begin{equation}   \label{eq.c1}
\ol{\delta}(1)\in\ol{\delta}_2,
\end{equation}
and in case $\wt{\delta}_2=\empt$ define $\ol{\delta}_2=\empt$.

Finally, let $\ol{f}_t:\D\ra\D$ and $\Pi_t:\D\ra \R^2\minus E^t$ be as in Lemma \ref{lema.fam-Pit}.\\
\\
\textbf{Main claims}. The following two claims will finish the proof of Lemma \ref{bp-ib1}.

\begin{claim}    \label{cl.bp-ib1}
If $\ol{f}_{N_2}(\ol{\eta}) \cap \ol{\delta} \neq \empt$, then $\wh{f}^{N_2}(\eta)$ has good intersection with $\wh{f}^{-N_1}(\beta_{1})$.
\end{claim}

\begin{claim}      \label{cl.bp-ib2}
$\ol{f}_{N_2}(\ol{\eta}) \cap \ol{\delta}\neq\empt$.
\end{claim}

To prove these claims, we begin with the following lemma.

\begin{lema}       \label{lema.fam-Pit2}
The following hold:
\begin{enumerate}
\item For all $t\in [0,N_2]$ we have $\Pi_t|_A = \Pi_0|_A$, where $A = \ol{L}(\ol{c}_0) \cup \ol{\kappa} \cup \ol{S} \cup \ol{\delta} \cup \ol{\delta}_1\cup\ol{\delta}_2$.   \label{itb1}
\item Let $\ast \in \{ \delta,\delta_1,\wt{\delta}_2 \}$. If $\ol{f}_{N_2}(\ol{\eta}) \cap \ol{*}\neq\empt$ then $\wh{f}^{N_2}(\eta)\cap \ast\neq\empt$.   \label{itb3}
\item $\forall\, n\in\Z \ \ \exists\, a>0$ such that $(\ol{f}_1)^n\ol{\eta}(t) \in L(\ol{c}_0)$, for $0< t <a$.     \label{itb2}
\item Each of the arcs $\ol{S}$, $\ol{\delta}$, $\ol{\delta}_1$ and $\ol{\delta}_2$ is contained in $R(\ol{c}_0)$.    \label{itb4}
\end{enumerate}
\end{lema}
\begin{prueba}
To prove item \ref{itb1} just note that $\Pi_0(A) \cap \txt{int}(c_k) =\empt$ for all $k>0$, by definition of the curves $c_k$, and therefore, by Item \ref{ita1} from Lemma \ref{lema.fam-Pit} we have that $\Pi_t|_A = \Pi_0|_A$, for all $t\in [0,N_2]$.

By this, item 1 of Lemma \ref{lema.fam-Pit} implies that
\begin{equation}   \label{eq.cc1}
\Pi_t \circ \ol{f}_t|_A = \wh{f}_t \circ \Pi_0|_A \ \ \ \forall\, t\in[0,N_2]
\end{equation}
which implies item 2. Also by Lemma \ref{lema.fam-Pit},
\begin{equation}   \label{eq.cc2}
\Pi_0 \circ \ol{f}_{-t} = (\wh{f}_{t}|_{\R^2\minus E})^{-1} \circ \Pi_t \ \ \ \forall\, t\in[0,N_2].
\end{equation}
Item 3 is an easy consequence of (\ref{eq.cc1}), (\ref{eq.cc2}) and the fact that for all $t\in\R$ the arc $\wh{f}_t\eta$ is issued from $\eta(0)=\wt{\gamma}(0)\in\txt{int}(c_0)$, as $c_0$ is oriented counterclock-wise.

To prove Item \ref{itb4}, note that by definition of $c_0$, the set Fill$(c_0)$ does not intersect any of the sets $T(\delta)$, $\delta_1$, $\wt{\delta}_2$ and $S$, and then the lifts $\ol{S}$, $\ol{\delta}$, $\ol{\delta}_1$ and $\ol{\delta}_2$ are contained in $R(\ol{c}_0)$. 
\end{prueba}

\begin{figure}[h]        
\begin{center} 
\includegraphics{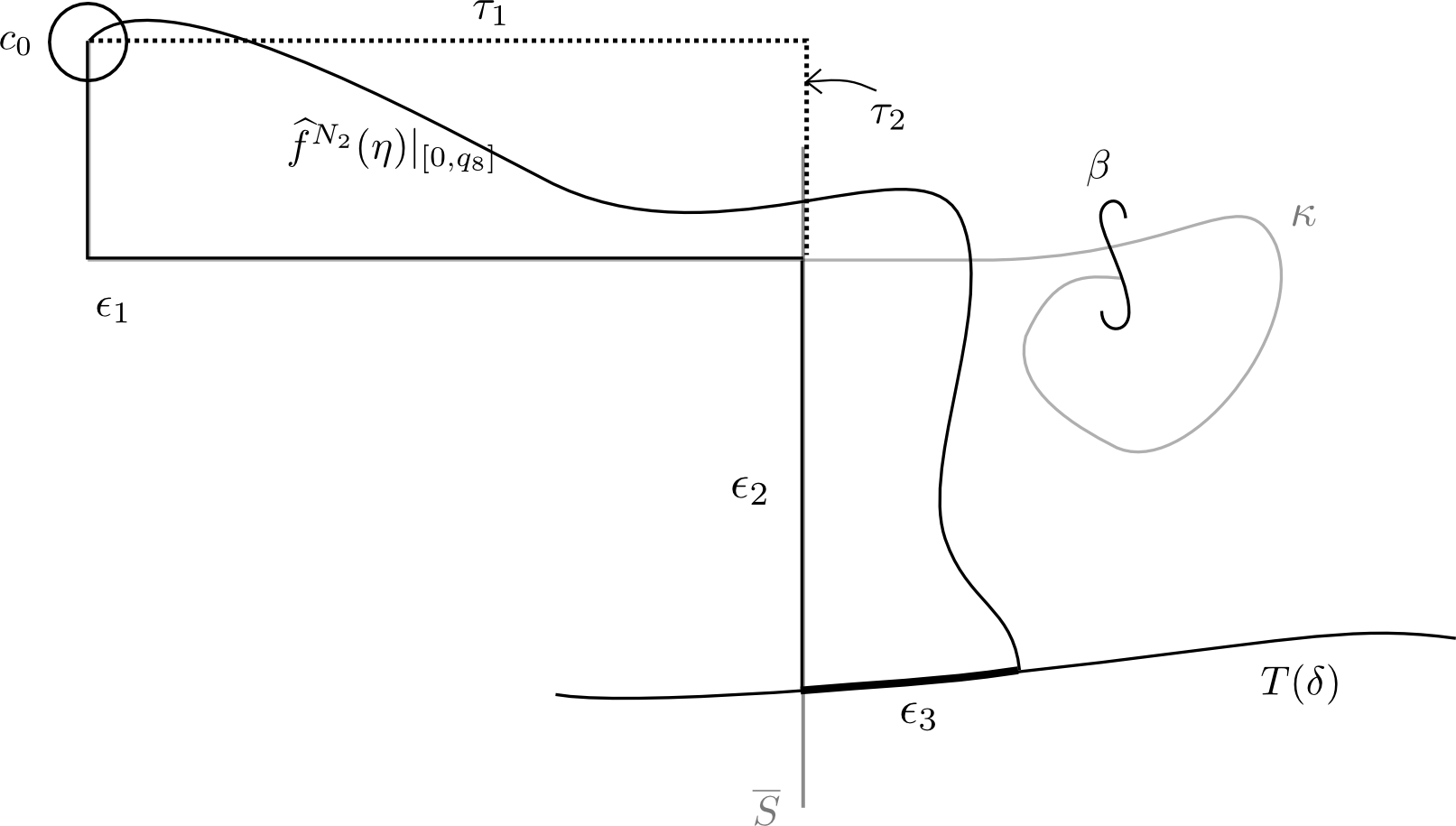}
\caption{The arcs $\e_i$ and $\tau_i$ from the proof of Claim \ref{cl.bp-ib1}.}
\label{fig.cl1}
\end{center}  
\end{figure}

\begin{prueba}[Proof of Claim \ref{cl.bp-ib1}.]
Without loss of generality, we may make the following assumption on the covering map $\Pi_0:\D\to \R^2\minus E$. If a sequence of loops $\sigma_i\subset \txt{int}(c_0)$ is such that $\eta(0)\in\txt{int}(\sigma_i)$ for all $i$ and diam$(\sigma_i)\to 0$ as $i\to\infty$, and if $\ol{\sigma}_i$ denotes the lift of $\sigma_i$ contained in $L(\ol{c}_0)$, then diam$(\ol{\sigma}_i)\to 0$, as $i\to\infty$, where in both cases diam$(X)$ denotes the euclidean diameter of a set $X$. 

By this assumption, the limit $\lim_{t\to 0}\ol{\eta}(t)=:z_0$ exists, $z_0\in\pr\D\subset\C$, and we may extend the curve $\ol{\eta}:(0,1]\ra\D$ continuously to a curve $\ol{\eta}:[0,1]\to \D\cup\{z_0\}$ as $\ol{\eta}(0) = z_0$. Analogously, for any $t$ extend $\ol{f}_{t}\ol{\eta}$ continuously to a curve $\ol{f}_t\ol{\eta}:[0,1]\to\D\cup\{z_0\}$. 

Define 
$$q_8= \min\{ t  \, : \, \ol{f}_{N_2}\ol{\eta}(t) \in \ol{\delta} \}$$ 
(such minimum exists, as by item \ref{itb4} from Lemma \ref{lema.fam-Pit2} $\ol{\delta}\subset R(\ol{c}_0)$). As $\D\cup\{z_0\}$ is simply connected, the arc $\ol{f}_{N_2}\ol{\eta}|_{[0,q_8]}$ is homotopic wfe to the arc $\ol{\e}= \ol{\e}_1\cdot\ol{\e}_2\cdot \ol{\e}_3$, where the arcs $\ol{\e}_i$ are such that:
\begin{itemize}
\item $\ol{\e}_1\subset \ol{\kappa}$ going from $\ol{\e}_1(0)=\ol{f}_{N_2}\ol{\eta}(0)=\ol{\eta}(0)$ to $\ol{\e}_1(1) = \ol{\kappa}(q_4)$, with $q_4$ as in (\ref{eq.10a});
\item $\ol{\e}_2 \subset \ol{S}$ going from $\ol{\e}_1(1)$ to $\ol{\e}_2(1)= \ol{S}(q_9)$, where $q_9= \min\{t \in [0,1] \, : \, \ol{S}(q_9) \in\ol{\delta}\}$;
\item $\ol{\e}_3\subset  \ol{\delta}$ going from $\ol{\e}_2(1)$ to $\ol{\e}_3(1)=\ol{f}_{N_2}\ol{\eta}(q_8)$.
\end{itemize}

By item \ref{itb2} from Lemma \ref{lema.fam-Pit2} we may define 
$$q_7 = \min\{ t\, : \, \ol{f}_{N_2}\ol{\eta}(t)\in \ol{c}_0 \}.$$

Define the arc $\e_1\subset\R^2$ by $\e_1(0)=\eta(0)$ and $\e_1(t)=\Pi_{N_2}(\ol{\e}_1(t))$ for $t\in(0,1]$. For $i=2,3$, let $\e_i=\Pi_{N_2}(\ol{\e}_i)$, and define $\e=\e_1\cdot\e_2\cdot\e_3$.

The homotopy wfe from  $\ol{f}_{N_2}\ol{\eta}|_{[0,q_8]}$ to $\ol{\e}$ and the projection $\Pi_{N_2}$ induce a homotopy wfe Rel$(E^{N_2})$ from $\wh{f}^{N_2}\eta|_{[0,q_8]}$ to $\e$ (recall that $E^t= \wh{f}_t(E)$). By Item \ref{itb1} from Lemma \ref{lema.fam-Pit2}
$$ \Pi_0|_{\ol{\kappa}\cup\ol{S}\cup\ol{\delta}} = \Pi_t|_{\ol{\kappa}\cup\ol{S}\cup \ol{\delta}} \ \ \ \forall\, t\in [0,N_2],$$
and then $\epsilon_1\subset\kappa$, $\epsilon_2\subset S$, and $\epsilon_3\subset T(\delta)$. Also, by definition of $\epsilon_1$ and $E^{N_2}$, the arc $\epsilon_1$ is homotopic wfe Rel$(E^{N_2})$ to the arc $\tau =\tau_1\cdot\tau_2$, where $\tau_1$ is straight horizontal and $\tau_2$ is straight vertical (see Fig. \ref{fig.cl1}). Then, 
$$ \e  \ \  \txt{is homotopic wfe Rel$(E^{N_2})$ to} \ \  \tau_1\cdot \tau_2\cdot \epsilon_2 \cdot \epsilon_3.$$
As $\wh{f}^{N_2}\eta|_{[0,q_8]}$ is homotopic wfe Rel$(E^{N_2})$ to $\e$, we obtain that 
\begin{equation}     \label{eq.25a}
\wh{f}^{N_2}\eta|_{[0,q_8]}   \ \ \ \txt{is homotopic wfe Rel$(E^{N_2})$ to} \ \ \ \tau_1\cdot \tau_2\cdot \epsilon_2 \cdot \epsilon_3.
\end{equation}
Note that, as $E^{N_2}  \supset\wh{f}^{-N_1}(F_{1})$, the homotopy in (\ref{eq.25a}) holds in particular Rel$(\wh{f}^{-N_1}(F_{1}))$. As
\begin{itemize}
\item $\tau_1$ is straight horizontal,
\item $\tau_2\cdot \epsilon_2$ is straight vertical, and
\item $\epsilon_3\subset T(\delta)$,
\end{itemize} 
we have that $\wh{f}^{N_2+n}\wt{\gamma}\supset \wh{f}^{N_2}\eta$ satisfies the definition of good intersection with $\wh{f}^{-N_1}(\beta_{1})$, with $\mu=\wh{f}^{N_2}\eta|_{[0,q_8]}$, $\nu_1=\tau_1$, $\nu_2=\tau_2\cdot\e_2$ and $\nu_3=\e_3$. This proves Claim \ref{cl.bp-ib1}.
\end{prueba}

\begin{figure}[h]        
\begin{center} 
\includegraphics{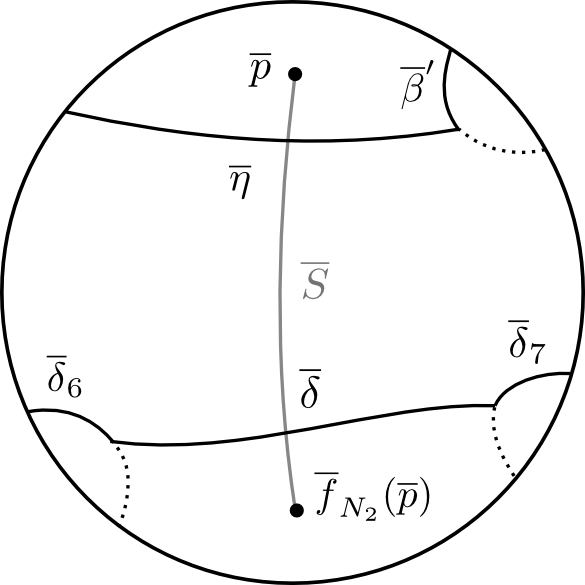}
\caption{The arcs $\ol{\beta}'$, $\ol{\delta}_6$ and $\ol{\delta}_7$ (in the case $\ol{\delta}_7\neq\empt$). From Lemma \ref{lema.cl.bp-ib2}, $\ol{p}\in L(\ol{\eta}\cdot\ol{\beta}')$ and $\ol{f}_{N_2}(\ol{p}) \in R(\Delta)$.}
\label{fig.delta67}
\end{center}  
\end{figure}

Our objective now is to prove Claim \ref{cl.bp-ib2}. To this end, we first make some definitions and prove a lemma related to them. Let 
$$\ol{\beta}' = \ol{\beta}|_{[q_{10},1)},$$ 
where $q_{10}$ is such that $\ol{\beta}(q_{10}) = \ol{\eta}(1)$ (see Fig. \ref{fig.delta67}). Recall our definitions of $\ol{\delta}_1$,$\ol{\delta}_2$ from (\ref{eq.d111}), (\ref{eq.c1}). Let 
$$\ol{\delta}_6 = \ol{\delta}_1|_{[q_{11},1)},$$ 
where $q_{11}$ is such that $ \ol{\delta}_1(q_{11}) = \ol{\delta}(0)$. If $\wt{\delta}_2\neq\empt$, define
$$ \ol{\delta}_7 = \ol{\delta}_2|_{(0,q_{12}]},  $$
where $q_{12}$ satisfies $ \ol{\delta}_2(q_{12}) = \ol{\delta}(1)$, and in case $\wt{\delta}_2=\empt$ define $\ol{\delta}_7=\empt$.

Define $\ol{p}\in\D$ as the lift of $p_0$ by $\Pi_0$ given by 
$$   \ol{p} = \ol{S}(0). $$

Let $\theta\subset\R^2$ be the straight vertical line passing through $p_0$, oriented downwards. Then, $c_0 \subset R(\theta)$ and $\beta\subset L(\theta)$. Let $\ol{\theta}\subset\D$ be a lift of $\theta$ by $\Pi_0$ passing through $\ol{p}$. It is easy to see that
\begin{equation}   \label{eq.c2}
\ol{c}_0 \subset R(\ol{\theta}) \ \ \txt{and} \ \ \ol{\beta}\subset L(\ol{\theta}).
\end{equation} 
Define
$$\Delta = \ol{\delta}_6^{-1}\cdot \ol{\delta} \cdot \ol{\delta}_7^{-1}.$$

Recall that IN$(\Gamma_1,\Gamma_2)$ denotes the intersection number of two curves $\Gamma_1$ and $\Gamma_2$ with disjoint ends (cf. Section \ref{sec.in}). Without loss of generality we may assume that the curves $\Delta$, $\ol{\eta}\cdot\ol{\beta}'$, $\ol{\theta}$ have disjoint ends, so their intersection number is well defined.

\begin{lema}       \label{lema.cl.bp-ib2}
The following hold:
\begin{enumerate}
\item  $\txt{IN}(\ol{\theta}, \ol{f}_{N_2}(\ol{\eta}\cdot\ol{\beta}')) = 1$,     \label{it.1a}
\item  $\ol{p} \in L(\ol{\eta} \cdot \ol{\beta}')$,      \label{it.2a}
\item  $\ol{f}_{N_2}(\ol{p}) \in R(\Delta)$, and       \label{it.3a}
\item  $\ol{L}(\ol{c}_0) \cup \ol{R}(\ol{\beta})  \subset L(\Delta)$.       \label{it.4a}
\end{enumerate}
\end{lema}
\begin{prueba}
By item \ref{itb2} from Lemma \ref{lema.fam-Pit2}, $\ol{f}_{N_2}\ol{\eta}(t)\in L(\ol{c}_0)$ for $t$ small. By item 3 from Lemma \ref{lema.fam-Pit}, $\ol{\beta}$ is a Brouwer line for $\ol{f}_{N_2}$, and then $\ol{f}_{N_2}(\ol{\eta}\cdot\ol{\beta}') \subset R(\ol{\beta})$. By this and by (\ref{eq.c2}), IN$(\ol{\theta}, \ol{f}_{N_2}(\ol{\eta}\cdot\ol{\beta}' ) = 1$, which proves Item \ref{it.1a}.

We now prove Item \ref{it.2a}. Recall that, by definition, $\eta$ is homotopic wfe Rel$(F_0)$ to an arc $\kappa=\kappa_1\cdot\kappa_2$, where $\kappa_1\subset L(\cl{C}_{-1})$ is straight vertical and $\kappa_2$ is below $p_0$. Analogously as in Item \ref{it.1a}, it follows 
\begin{equation}      \label{eq.30a}
\txt{IN}(\ol{\theta}, \ol{\kappa}\cdot\ol{\beta}') = 1.
\end{equation}
By definition, $p_0$ is above $\kappa_2\cup\beta$, and then if $q_{13}\in\R$ is such that $\ol{\theta}(q_{13}) = \ol{p}$, we have that
\begin{equation}      \label{eq.31a}
q_{13} < \min \{ t\in\R \, : \, \ol{\theta}(t) \in \ol{\kappa}\cdot\ol{\beta}' \}.
\end{equation} 
By (\ref{eq.30a}), (\ref{eq.31a}) and by Lemma \ref{lema.in} we have
$$ \ol{p}  \in  L(\ol{\kappa}\cdot \ol{\beta}' ).$$
As $\kappa$ is homotopic wfe Rel$(F_0)$ to $\eta$, and as $p_0\in F_0$, we conclude that
$$ \ol{p} \in L(\ol{\eta}\cdot \ol{\beta}'),$$
which proves item \ref{it.2a}.

We now prove Item \ref{it.3a}. Analogously as in Item \ref{it.1a}, we have
\begin{equation}     \label{eq.32a}
\txt{IN}(\ol{\theta}, \Delta)  = 1.
\end{equation}
Also, as $\wh{f}^{N_2}(p_0)$ is below $T(\delta)$, we have that 
$$ \max \{ t\in\R \, : \, \theta(t) \in T(\delta) \} < q_{14},$$
where $q_{14}$ is such that $\theta(q_{14}) = \wh{f}^{N_2}(p_0).$ This implies that
\begin{equation}     \label{eq.33a}
\max \{ t\in\R \, : \, \ol{\theta}(t) \in \Delta \} = \max \{ t\in\R \, : \, \ol{\theta}(t) \in \ol{\delta} \} < q_{14}.
\end{equation}
By (\ref{eq.32a}), (\ref{eq.33a}) and by Lemma \ref{lema.in} we have that
$$\ol{f}_{N_2}(\ol{p}) \in R( \Delta),$$
which proves item \ref{it.3a}.

Finally, we prove Item \ref{it.4a}. As $T(\delta)$ is below $\kappa$, we have the inequality
\begin{align*}
 \max \{ t\in\R \, : \, \ol{\theta}(t) \in \ol{\kappa} \cup\ol{L}(\ol{c}_0) \cup \ol{R}(\ol{\beta}) \} &= \max \{ t\in\R\, : \, \ol{\theta}(t) \in \ol{\kappa} \} \\
                     &< \min \{ t\in\R \, : \, \ol{\theta}(t) \in \ol{\delta} \}  \\
										 &= \min \{ t\in\R \, : \, \ol{\theta}(t) \in \Delta \}.
\end{align*}
By this, as $\txt{IN}(\ol{\theta}, \Delta) = 1$ (by (\ref{eq.32a})) and by Lemma \ref{lema.in}, we get that the connected set $\ol{\kappa} \cup \ol{L}(\ol{c}_0) \cup \ol{R}(\ol{\beta})$ is contained in $L(\Delta)$. This implies Item \ref{it.4a}. 
\end{prueba}

We may now proceed to the proof of Claim \ref{cl.bp-ib2}.

\begin{prueba}[Proof of Claim \ref{cl.bp-ib2}]
We will first show that 
\begin{equation}     \label{eq.36a}
\ol{f}_{N_2}\ol{\eta} \cap \Delta \neq\empt. 
\end{equation}
Recall that by Item \ref{itb2} from Lemma \ref{lema.fam-Pit2}, there is $a>0$ such that 
\begin{equation}   \nonumber
\ol{f}_{N_2}\ol{\eta}(t) \in L(\ol{c}_1), \ \ \ \ \txt{ for $t< a$}.
\end{equation}
 By Item \ref{ita3} of Lemma \ref{lema.fam-Pit} $\ol{\beta}$ is a Brouwer curve for $\ol{f}_{N_2}$, and therefore 
\begin{equation}    \nonumber 
\ol{f}_{N_2}(\ol{\eta}(1))\in\ol{f}_{N_2}(\ol{\beta})\subset R(\ol{\beta}).
\end{equation}
Also, by Item \ref{it.4a} from Lemma \ref{lema.cl.bp-ib2}, we know that 
\begin{equation}     \label{eq.37a}
\ol{L}(\ol{c}_0) \cup \ol{R}(\ol{\beta}) \subset L(\Delta),
\end{equation}
and therefore, to prove (\ref{eq.36a}) it suffices to show that 
\begin{equation}      \label{eq.34a}
\ol{f}_{N_2}\ol{\eta} \cap \ol{R}(\Delta ) \neq\empt.
\end{equation}
We suppose that (\ref{eq.34a}) does not hold, and we will find a contradiction. By Item \ref{it.2a} from Lemma \ref{lema.cl.bp-ib2}, we know that $\ol{f}_{N_2}(\ol{p}) \in \ol{f}_{N_2}(L(\ol{\eta}\cdot \ol{\beta}')) = L(\ol{f}_{N_2}(\ol{\eta}\cdot \ol{\beta}'))$. By this and by item 3 from Lemma \ref{lema.cl.bp-ib2} we have that
$$L(\ol{f}_{N_2}(\ol{\eta}\cdot \ol{\beta}')) \cap R(\Delta) \neq\empt.$$
As $\ol{\beta}'\cap \Delta =\empt$, we obtain that, if (\ref{eq.34a}) does not hold, then
\begin{equation}     \label{eq.35a}
L(\ol{f}_{N_2}(\ol{\eta} \cdot \ol{\beta}')) \supset R(\Delta).
\end{equation}
On the other hand, as $\ol{f}_{N_2}\ol{\eta}$ intersects $L(\ol{c}_0)$ and $R(\ol{\beta})$, and by (\ref{eq.37a}), we have 
$$\ol{f}_{N_2}\ol{\eta} \subset L(\Delta).$$
Then, we must have 
$$\max \{ t\in\R \, : \, \ol{\theta}(t)\in \ol{f}_{N_2}\ol{\eta} \} < \min \{ t\in\R \, : \, \ol{\theta}(t)\in \Delta \}.$$
As IN$(\ol{\theta}, \ol{f}_{N_2}\ol{\eta}\cdot\ol{\beta}') = 1$ (by item \ref{it.1a} from Lemma \ref{lema.cl.bp-ib2}), by Lemma \ref{lema.in} it holds
$$ R(\Delta) \subset R(\ol{f}_{N_2}(\ol{\eta} \cdot \ol{\beta}')),$$
which contradicts (\ref{eq.35a}). This contradiction proves (\ref{eq.34a}), and therefore (\ref{eq.36a}).

We will now see that 
\begin{equation}      \label{eq.38a}
\ol{f}_{N_2}\ol{\eta}\cap\ol{\delta}_6= \ol{f}_{N_2}\ol{\eta}\cap \ol{\delta}_7 =\empt.
\end{equation}
This, together with (\ref{eq.36a}), will imply that $\ol{f}_{N_2}\ol{\eta}$ intersects $\ol{\delta}$, which proves Claim \ref{cl.bp-ib2}. 

Note that, by item \ref{ita3} from Lemma \ref{lema.fam-Pit}, $\ol{\delta}_1$ and $\ol{\delta}_2$ are Brouwer curves for $\ol{f}_1$. If we had that $\ol{f}_{N_2}\ol{\eta}\cap \ol{\delta}_6\neq\empt$, as $\ol{\delta}_6\subset\ol{\delta}_1$ this would give 
\begin{equation}      \label{eq.39a}
(\ol{f}_1)^{-m} \ol{f}_{N_2}\ol{\eta} \cap L(\ol{\delta}_1) = (\ol{f}_1)^{-m+N_2} \ol{\eta} \cap L(\ol{\delta}_1) \neq\empt \ \ \ \ \forall m\geq 0.
\end{equation}
On the other hand, by definition $\eta$ is contained in some iterate $\wh{f}^{n} \tl{\gamma}$, $n>0$, of the curve $\tl{\gamma}$, and $\tl{\gamma}$ is disjoint from $\delta_1\cup\delta_2$, which implies
\begin{equation}   \label{eq.d112}
(\ol{f}_1)^{-n}\ol{\eta}\cap\ol{\delta}_1=\empt.
\end{equation}
By (\ref{eq.37a}), $L(\ol{c}_0)\subset L(\Delta) \subset R(\ol{\delta}_1)$, and then, by item \ref{itb2} from Lemma \ref{lema.fam-Pit2} we have that $(\ol{f}_1)^{-n} \ol{\eta}(t)$ belongs to $L(\ol{c}_0)\subset R(\ol{\delta}_1)$ for $t$ small. By this and by (\ref{eq.39a}) we have that $(\ol{f}_1)^{-n}\ol{\eta}$ intersects $\ol{\delta}_1$, which contradicts (\ref{eq.d112}).
Therefore, we must have that $\ol{f}_{N_2}\ol{\eta}\cap \ol{\delta}_6 = \empt$. 

The same argument shows that $\ol{f}_{N_2}(\ol{\eta})\cap \ol{\delta}_7 = \empt$. This shows (\ref{eq.38a}) and finishes the proof of the claim. 
\end{prueba}

\subsubsection{Proof of Lemma \ref{ib-bp}}       \label{sec.ib-bp}

By definition, if $\wh{f}^n(\wt{\gamma})$ has good intersection with $\wh{f}^{-N_1}(\beta_i)$ we have that $\wh{f}^n(\tl{\gamma})$ contains a subarc $\mu$ which goes from $\wh{f}^n(\tl{\gamma}(0))=\tl{\gamma}(0)$ to $\wh{f}^{-N_1}(\beta_i)$, which is homotopic wfe Rel$(\wh{f}^{-N_1}(F_i))$ to an arc $\nu$ of the form $\nu=\nu_1\cdot\nu_2\cdot\nu_3$, where $\nu_1$ is horizontal, $\nu_2$ is vertical, and $\nu_3$ is contained $T^i(\delta)$ (cf. Definition \ref{def.ib} and Fig. \ref{fig.ib}). Clearly, $\nu_1\cdot\nu_2$ is homotopic wfe Rel$(\wh{f}^{-N_1}(F_i))$ to the arc $\nu'_1\cdot\nu'_2$, where $\nu'_1$ is straight vertical and $\nu'_2$ is straight horizontal. If $\nu'=\nu'_1\cdot\nu'_2\cdot\nu_3$, we then have that $\mu$ is homotopic wfe Rel$(\wh{f}^{-N_1}(F_i))$ to $\nu'$.

Let $\e$ be the arc $\e(t) = \wh{f}_t (\mu(1))$, $0\leq t \leq N_1$, and reparametrize $\e$ so that it is defined in $[0,1]$. By definition of $p_0$ in $\S$\ref{sec.plg}, we have that $\{\wh{f}_t(p_0)\, : \, t\in[-N_1,0]\}$ is above $\{\wh{f}_t(\beta_0)\, : \, t\in[-N_1,0]\}$, and by item 2 of Claim \ref{claim.N2} we have that $\wh{f}^{-j(N_1+N_2)}(p_j)$ is above $p_0$ for all $j>0$. Therefore, by periodicity 
\begin{equation}   \label{eq.cc5}
\{\wh{f}_t(\wh{f}^{-j(N_1+N_2)}(p_j))\, : \, t\in[-N_1,0]\} \ \ \txt{is above} \ \ \{\wh{f}_t(\beta_i)\, : \, t\in[-N_1,0]\} \ \ \txt{for all $j>i$.}
\end{equation}
Also, the set $\{\wh{f}_t(\beta_i)\, : \, t\in[-N_1,0]\}$ is disjoint from the sets sing$(\beta_j),\txt{sing}(T^j(\delta_1)),\txt{sing}(T^j(\wt{\delta}_2))$ for any $j$. By, by (\ref{eq.cc5}) and as $\e\subset \{\wh{f}_t(\beta_i)\, : \, t\in[-N_1,0]\}$ we conclude that 
\begin{equation}  \label{eq.cc3}
\{ \wh{f}_t(F_i) \, : \, t\in[-N_1,0]\} \cap \e = \empt. 
\end{equation}
Similarly, we have 
\begin{equation}  \label{eq.cc4}
\{ \wh{f}_t(F_i) \, : \, t\in[-N_1,0]\} \cap \nu' = \empt. 
\end{equation}
By (\ref{eq.cc3}) and (\ref{eq.cc4}), Lemma \ref{isot.lema2} gives us that $\nu'\cdot\e$ is homotopic wfe Rel$(\wh{f}^{N_1} (\wh{f}^{-N_1}(F_i))) = \txt{Rel}(F_i)$ to $\wh{f}^{N_1}\mu$.

Let $\ol{\nu}_3$ be a lift of $\nu_3$ to the universal cover $\D$ of $\R^2\minus (F_i \cap \txt{sing}(\wh{\cl{F}}) )$, let $\ol{\e}$ be the lift of $\e$ to $\D$ with $\ol{\e}(0) = \ol{\nu}_3(1)$, and let $\ol{\beta}_i$ be a lift of $\beta_i$ such that $\ol{\e}(1) \in \ol{\beta}_i$. Observe that, by Theorem \ref{teo.jau}, $\e$ is homotopic wfe Rel$(\txt{sing}(\wh{\cl{F}}))$ to an arc $\e'$ which is positively transverse to $\wh{\cl{F}}$, and then $\ol{\e}(0)\in L(\ol{\beta}_i)$. As $\nu_3\cap\beta_i=\empt$, we have also $\ol{\nu}_3\cdot\ol{\e}(0)=\ol{\nu}_3(0)\in L(\ol{\beta}_i)$. Note also that the arcs $\e$ and $\nu_3$ are contained in the strip
$$ B_i=\{ x \in\R^2 \, : \, x \txt{ is below $p_i$ and above $\beta_i\cup T^i(\delta_1\cup\wt{\delta}_2)$} \},$$ 
and that if $\ol{\nu}'_2$ is the lift of $\nu'_2$ with $\ol{\nu}'_2(1)=\ol{\nu}_3(0)$, we have $\ol{\nu}'_2(0)\in L(\ol{\beta}_i)$ (as $\nu'_2\cap\beta_i=\empt$). 

We conclude that $\wh{f}^{N_1}\mu$ is homotopic wfe Rel$(F_i)$ to the arc $\nu' \cdot \e = \nu'_1\cdot \nu'_2 \cdot \nu'_3 \cdot \e$, where:
\begin{itemize}
\item $\nu'_1$ is straight vertical, and
\item $\nu'_2\cdot\nu_3\cdot\e$ is contained in $B_i\minus F_i$ and lifts to the arc $\ol{\nu}'_2\cdot\ol{\nu}_3\cdot\ol{\e}$, with $\ol{\nu}'_2\cdot\ol{\nu}_3\cdot\ol{\e}(0)\in L(\ol{\beta}_i)$.
\end{itemize}
This means that $\wh{f}^{n+N_1}(\wt{\gamma}) \supset \wh{f}^{N_1}\mu$ satisfies the definition of being well positioned with respect to $p_i$, with $\eta=\wh{f}^{N_1}\mu$, $\kappa_1=\nu'_1$ and $\kappa_2=\nu'_2\cdot\nu_3\cdot\e$, which proves Lemma \ref{ib-bp}.

\subsection{Construction of the arc $\delta$; proof of Lemma \ref{lema.delta1}.}     \label{sec.delta}

We begin by recalling the statement of the lemma. 

\begin{lema.sn}         
There exists an arc $\delta$ such that:
\begin{enumerate}
\item $\txt{int}(\delta) \subset \wh{f}^{-N_1}\beta$,
\item $\delta(0)\in\cl{C}_{-2}$, $\delta(1)\in\cl{C}_{0}$, $\txt{int}({\delta})\subset R(\cl{C}_{-2})\cap L(\cl{C}_{0})$,
\item $\delta$ leaves a leaf $\delta_1\subset \cl{C}_{-2}$ on $t=0$ by the right (cf. Definition \ref{def.btl}), 
\item one of the following holds:
	\begin{enumerate}
	\item $\delta(1) \in \txt{sing}(\wh{\cl{F}})\cap\cl{C}_0$. 
	\item $\delta$ arrives in a leaf $\delta_2\subset \cl{C}_{0}$ on $t=1$ by the right (see Fig. \ref{fig.delta}).
	\end{enumerate}  
\end{enumerate}
\end{lema.sn}

We also recall from Notation \ref{notacion1} our definition of $\wt{\delta}_2=\empt$ in case item $4(a)$ holds, and $\wt{\delta}_2=\delta_2$ in case item $4(b)$ holds.

\subsubsection{Case that either $\alpha(\gamma)$ or $\omega(\gamma)$ is a singularity of $\wh{\cl{F}}$, and either $\alpha(\beta)$ or $\omega(\beta)$ is a singularity of $\wh{\cl{F}}$.}           \label{sec.sing}

We will assume that $\omega(\gamma)$ and $\omega(\beta)$ are singularities, $\omega(\gamma)=\{s_1\}$ and $\omega(\beta) = \{s_2\}$, for some $s_1,s_2\in \txt{sing}(\wh{\cl{F}})$. The complementary case that (at least) one of the sets $\alpha(\gamma)$ and $\alpha(\beta)$ is a singularity, and (at least) one of the sets $\omega(\gamma)$ and $\omega(\beta)$ is not a singularity will follow from the same argument.

\begin{definicion}        \label{def.removable}
We say that an arc $\eta$ has a \textit{removable intersection} in $t\in (0,1)$ with $\cl{C}_{-2}$ if there exists a leaf $\e$ of $\wh{\cl{F}}$ contained in $\cl{C}_{-2}$ and $t'\in(t,1)$ such that:
\begin{itemize}
\item $\eta(t),\eta(t')\in\e$, and
\item $\eta|_{[t,t']}$ is homotopic wfe Rel$(\txt{sing}(\wh{\cl{F}}))$ to the subarc of $\e$ going from $\eta(t)$ to $\eta(t')$. 
\end{itemize} 
We call the interval $[t,t']$ a \textit{removable interval} for $\eta$. Also, we denote by $U_{\eta}\subset(0,1)$ the union of the interiors of every removable interval for $\eta$.
\end{definicion}

\begin{figure}[h]        
\begin{center} 
\includegraphics{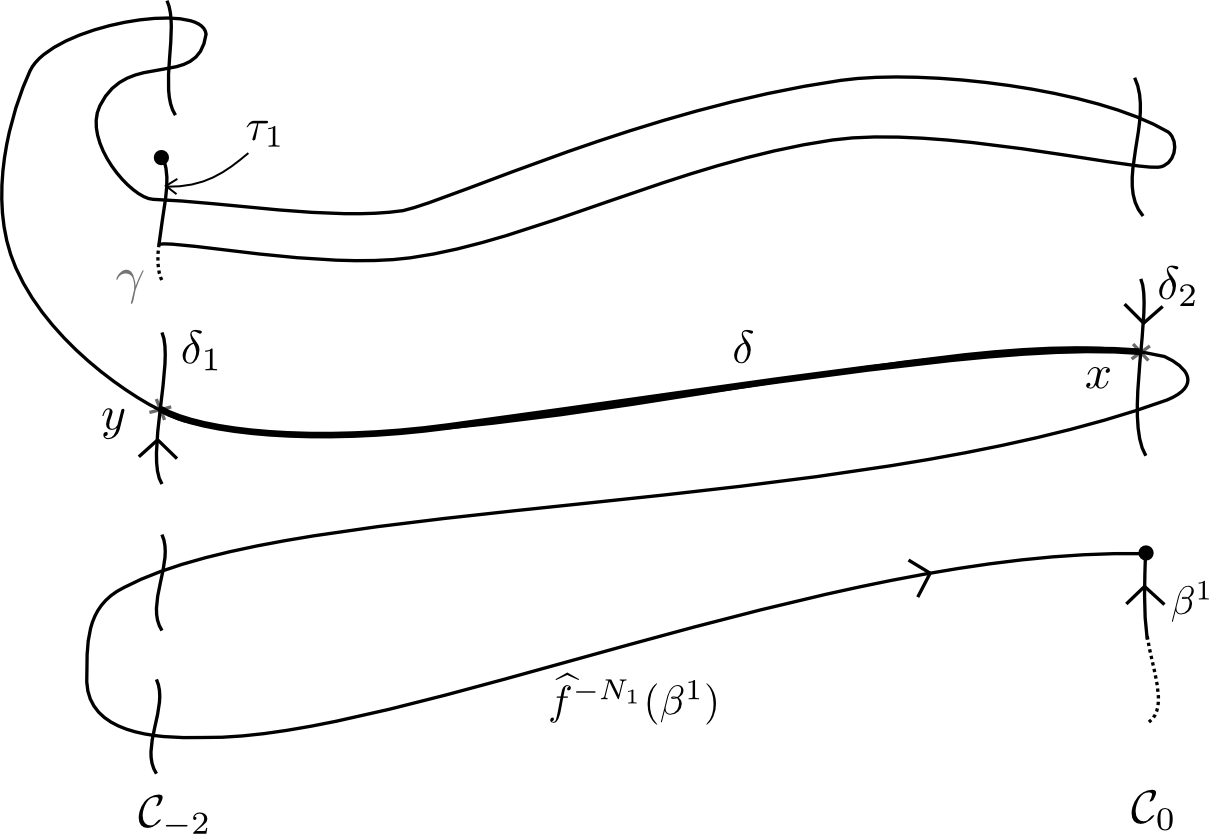}
\caption{An illustration of the construction of $\delta$, in the case that $\omega(\gamma)$ and $\omega(\beta)$ are singularities.}
\label{fig.delta1}
\end{center}  
\end{figure}

Extend the leaves $\gamma,\beta:(0,1)\to\R^2$ to $(0,1]$, so that $\gamma(1)=s_1$ and $\beta(1)=s_2$. Recall that by definition, $\beta$ is the arc that contains the last point of intersection of $\wh{f}^{N_1}\gamma$ with $\cl{C}_0$. That is, $\beta$ contains the point $z=\wh{f}^{N_1}\gamma(t_1)$, where
$$t_1= \max \{ t\in(0,1) \, : \,  \wh{f}^{N_1}\gamma(t) = \cl{C}_0 \}.$$
Let $\beta^1\subset\beta$ be the subarc of $\beta$ going from $z$ to $\beta(1)=s_2$ (see Fig. \ref{fig.delta1}). Consider the set $U_{\wh{f}^{-N_1}\beta^1}\subset(0,1)$ from Definition \ref{def.removable}, which is the union of the interiors of the removable intervals for $\wh{f}^{-N_1}\beta^1$. There exist two possibilities:
\vspace{2mm}\\
\textit{Case 1}: For every $t\in(0,1)$ such that $\wh{f}^{-N_1}\beta^1(t)\in \cl{C}_0$, it holds $t\in U_{\wh{f}^{-N_1}\beta^1}$. \\
In this case, let
$$ t_3 = \max \{ t\in (0,1) \, : \, \wh{f}^{-N_1}\beta^1(t) \in \cl{C}_{-2} \},$$
and define $\delta_1\subset\cl{C}_{-2}$ as the leaf of $\wh{\cl{F}}$ containing $\wh{f}^{-N_1}\beta(t_3)$. Then, $\wh{f}^{-N_1}\beta^1$ must leave $\delta_1$ in $t_3$ by the right; otherwise, it would leave by the left, and by Lemma \ref{lemabrouwer2} we would have 
$$\beta^1\cap\cl{C}_{-2} = \wh{f}^{N_1}(\wh{f}^{-N_1}\beta^1)\cap \cl{C}_{-2} \supset \wh{f}^{N_1}(\wh{f}^{-N_1}\beta^1)\cap\delta_1  \neq\empt,$$ 
which contradicts the fact that $\beta^1\subset\cl{C}_0$. 

Also, $\wh{f}^{-N_1}\beta^1|_{(t_3,1)}\subset R(\cl{C}_{-2})\cap L(\cl{C}_0)$. To see this, note that, by definition of $t_3$, if $t\in(t_3,1)$ is such that $\wh{f}^{-N_1}\beta^1(t)\in\cl{C}_0$, then $t$ must be contained in $(0,1)\minus U_{\wh{f}^{-N_1}}\beta^1$, which is excluded from Case 1. 

Define then 
$$\delta=\wh{f}^{-N_1}\beta^1|_{[t_3,1]}$$ 
and reparametrize $\delta$ to be defined in $[0,1]$. By the above, $\delta$ satisfies items 1 to 3 and 4(a) of the lemma. 
\vspace{2mm} \\
\textit{Case 2}: There exists $t\in(0,1)$ such that $\wh{f}^{-N_1}\beta^1(t) \in\cl{C}_0$, and $t\notin U_{\wh{f}^{-N_1}\beta^1}.$ \\
Let $x= \wh{f}^{-N_1}\beta^1(t_4)$, where
$$t_4= \min \{ t\in (0,1) \minus U_{\wh{f}^{-N_1}\beta^1}  \, : \, \wh{f}^{-N_1}\beta^1(t) \in \cl{C}_0 \}.$$
Define 
\begin{equation}   \label{eq.c12}
\beta^2= \wh{f}^{-N_1}\beta^1|_{[0,t_4]}
\end{equation}
(reparametrized to $[0,1]$), and let $\delta_2$ be the leaf of $\wh{\cl{F}}$ contained in $\cl{C}_0$ and containing $x= \beta^2(1)$. By definition, $\beta^2$ is homotopic wfe Rel$(\txt{sing}(\wh{\cl{F}}))$ to an arc whose first intersection with $\cl{C}_0$ is the endpoint $x=\beta^2(1)\in\delta_2$.

\begin{claim}    \label{cl.s1}
$\beta^2$ arrives in $\delta_2$ in $t=1$ by the right.
\end{claim}

\begin{prueba}
Let $\tau_1\subset\gamma$ be the arc going from $\gamma(1)=s_1$ to $\beta^2(0)=\wh{f}^{-N_1}\beta^1(0)$. Suppose that $\beta^2$ does not arrive in $\delta_2$ in $1$ by the right. Then it arrives by the left, and then by Lemma \ref{lemabrouwer2} we have that 
\begin{equation}     \label{eq.40a}
\wh{f}^{N_1}(\tau_1\cdot\beta^2)\cap\delta_2\neq\empt. 
\end{equation}
By construction of $\beta^1$ and $\beta^2$, we have that 
\begin{equation}    \label{eq.41a}
\wh{f}^{N_1}(\tau_1)\cap\cl{C}_0=\{ \wh{f}^{N_1}\tau_1(1) \} = \{ \wh{f}^{N_1}\beta^2(0) \} \subset \wh{f}^{N_1}\beta^2 \cap \cl{C}_0.
\end{equation}
By this, and as $\wh{f}^{N_1}(\beta^2)\subset\beta$, (\ref{eq.40a}) implies that 
$$\delta_2=\beta.$$ 
We will see that this yields a contradiction, and therefore $\beta^2$ must arrive in $\delta_2$ in 1 by the right.

Let $\ol{\beta}$ be a lift of $\beta$ to the universal cover of $\R^2\minus \txt{sing}(\beta)$. Let $\ol{\beta}^2$ be the lift of $\beta^2$ such that $\ol{\beta}^2(1)\in\ol{\beta}$, and let $\ol{\tau}_1$ be the lift of $\tau_1$ such that $\ol{\tau}_1(1)=\ol{\beta}^2(0)$. By construction, $\beta^2$ is homotopic wfe Rel$(\txt{sing}(\wh{\cl{F}}))$ to an arc whose only intersection with $\cl{C}_0$ is in $\beta^2(1)\in\beta$. Therefore, by our assumption that $\beta^2$ arrives in $\delta_2=\beta$ in $1$ by the left, we have that 
$$\ol{\tau}_1\cdot\ol{\beta}^2(0)\in L(\ol{\beta}).$$
If $\ol{f}$ is the canonical lift of $\wh{f}$, $\ol{\beta}$ is a Brouwer curve for $\ol{f}$, and therefore by Lemma \ref{lemabrouwer2} $\ol{f}^{N_1}(\ol{\tau}_1\cdot\ol{\beta}^2) \cap \ol{\beta}\neq\empt$. By (\ref{eq.41a}), this yields 
$$\ol{f}^{N_1}\ol{\beta}^2\cap\ol{\beta} \supset \ol{f}^{N_1}(\ol{\tau}_1\cdot\ol{\beta}^2) \cap \ol{\beta} \supset \{ \wh{f}^{N_1}\ol{\beta}^2(0) \}  \neq\empt,$$
and as $\beta^2\subset \wh{f}^{-N_1}\beta$, we must have
$$\ol{f}^{N_1}\ol{\beta}^2 \subset\ol{\beta}.$$
As $\ol{\beta}^2(1)\in \ol{\beta}$, we have that $\ol{f}^{N_1}\ol{\beta}\cap \ol{\beta}\neq\empt$, which contradicts the fact that $\ol{\beta}$ is a Brouwer curve for $\ol{f}$. This is the sought contradiction, and this finishes the proof of the claim. 
\end{prueba}

Let $y = \beta^2(t_5)$, where
$$ t_5 = \max \{ t\in (0,1) \, : \, \beta^2(t)\in \cl{C}_{-2} \}, $$
and let $t_6\in[0,t_4]$ be such that $\wh{f}^{-N_1}\beta^1(t_6)=\beta^2(t_5)$.

Define $\delta_1$ as the leaf of $\wh{\cl{F}}$ contained in $\cl{C}_{-2}$ which contains $y$.

\begin{claim}     \label{cl.s2}
$\wh{f}^{-N_1}\beta^1$ leaves $\delta_1$ in $t_6$ by the right.
\end{claim}

Before proving this claim we see how it implies Lemma \ref{lema.delta1} in Case 2. Define $\delta$ as the subarc of $\wh{f}^{-N_1}\beta$ going from $y$ to $x$. As $\delta\subset \wh{f}^{-N_1}\beta$, $\delta$ satisfies item 1 from the lemma. By definition of $x$ and $y$, int$(\delta)\subset R(\cl{C}_{-2})\cap L(\cl{C}_0)$, and then $\delta$ satisfies item 2. By Claim \ref{cl.s1}, $\delta$ arrives in the leaf $\delta_2\subset\cl{C}_0$ in $t=1$ by the right, and by Claim \ref{cl.s2} $\delta$ leaves the leaf $\delta_1\subset\cl{C}_{-2}$ in $t=0$ by the right, and then $\delta$ satisfies items 3 and 4(b) of the lemma. This concludes the proof. 

We now give the proof of Claim \ref{cl.s2}.

\begin{prueba}[Proof of Claim \ref{cl.s2}]
If $\wh{f}^{-N_1}\beta^1$ does not leave $\delta_1$ in $t_6$ by the right, it leaves by the left. We will show that this implies a contradiction.

Consider the set $U_{\wh{f}^{-N_1}\beta^1}\subset[0,1]$, which is a (possibly empty) union of open intervals $I_i$ (cf. Definition \ref{def.removable}). We know that, for each $i$, $\wh{f}^{-N_1}\beta^1|_{I_i}$ is homotopic wfe Rel$(\txt{sing}(\wh{\cl{F}}))$ to an arc $\e_i\subset\cl{C}_{-2}$ contained in a leaf of $\wh{\cl{F}}$. 

Consider the canonical $\ol{f}$ of $\wh{f}$ to the universal cover of $\R^2\minus \txt{sing}(\delta_1)$, let $\ol{\delta}_1$ be any lift of $\delta_1$, and let $\ol{\beta}$ be a lift of $\wh{f}^{-N_1}\beta^1$ such that $\ol{\beta}(t_6) \in \ol{\delta}_1$. By definition, $t_6$ is not a removable intersection for $\wh{f}^{-N_1}\beta^1$. By this, and as we are assuming that $\wh{f}^{-N_1}\beta^1$ leaves $\delta_1$ in $t_6$ by the left, we have that $\ol{\beta}$ leaves $\ol{\delta}_1$ in $t_6$ by the left and
\begin{equation}    \label{eq.c11}
\ol{\beta}(1) \in L(\ol{\delta}_1).
\end{equation}
By Lemma \ref{lemabrouwer1}, this yields $\ol{f}^{N_1}\ol{\beta}\cap\ol{\delta}_1\neq\empt$, and therefore
$$ \beta^1\cap\cl{C}_{-2} = \wh{f}^{N_1}(\wh{f}^{-N_1}\beta^1) \cap \cl{C}_{-2} \supset \wh{f}^{N_1}(\wh{f}^{-N_1}\beta^1) \cap \delta_1 \neq\empt,$$
which is the desired contradiction, as $\beta^1\subset \beta \subset \cl{C}_0$. 
\end{prueba}

\subsubsection{Case that none of the sets $\alpha(\gamma)$, $\omega(\gamma)$, $\alpha(\beta)$, and $\omega(\beta)$ consist of a singularity.}         \label{sec.nosing}

We recall the definition of $\beta$ from Section \ref{sec.plg}, for the case that neither $\omega(\gamma)$ nor $\omega(\beta)$ consist of a singularity.

We first defined a straight horizontal arc $\lambda_1$ going leftwards from a singularity $s\in\txt{Fill}(\omega(\gamma))$ to a point of $\wh{f}^{N_1}(\gamma)$ sufficiently close to $\wh{f}^{N_1}(\omega(\gamma))$ so that $\wh{f}^{-N_1}(\lambda_1)\cap\cl{C}_0=\empt$. Then, we defined $\lambda_2$ as a subarc of $\wh{f}^{N_1}\gamma$ going from $\lambda_1(1)$ to a point $z\in\wh{f}^{N_1}\gamma\cap\cl{C}_0$ such that $\lambda_2|_{[0,1)}\subset L(\cl{C}_0)$. Finally, we set $\lambda=\lambda_1\cdot\lambda_2$, and defined $\beta$ as the leaf of $\wh{\cl{F}}$ that contains $\lambda(1)$.

\begin{claim}
The arc $\lambda=\lambda_1\cdot\lambda_2$ arrives in $\beta$ in $t=1$ by the left.
\end{claim}

\begin{prueba}
Suppose that $\lambda$ does not arrive in $\beta$ in 1 by the left. Then $\lambda$ must arrive by the right, and by Lemma \ref{lemabrouwer2} we have that 
$$ \wh{f}^{-N_1}\lambda\cap \cl{C}_0 \supset \wh{f}^{-N_1}\lambda\cap \beta \neq\empt,$$ 
which contradicts the fact that $\wh{f}^{-N_1}\lambda_1 \cap \cl{C}_0=\empt$. Therefore, we must have that $\lambda$ arrives in $\beta$ by the left.
\end{prueba}

\begin{lema}     \label{lema.delta2}
There exist arcs $\beta^1$ and $\beta^2$ such that:
\begin{enumerate}
\item $\beta^1(0) \in \txt{sing}(\wh{\cl{F}})\cap \txt{Fill}(\omega(\beta))$,
\item $(\beta^1\cup\wh{f}^{N_1}\beta^1) \cap \cl{C}_{-2} = \empt$, 
\item $\beta^2(0)=\beta^1(1)$, $\beta^2 \subset \wh{f}^{-N_1}\beta$ and $\beta^2(1)= \wh{f}^{-N_1}\lambda(1)$, 
\item there exists $t_1\in[0,1)$ such that $\beta^2(t_1)\in\beta$ and $\beta^2(t) \in R(\cl{C}_{-2})$ for all $t\in [0,t_1]$. 
\end{enumerate}
\end{lema}

\begin{figure}[h]        
\begin{center} 
\includegraphics{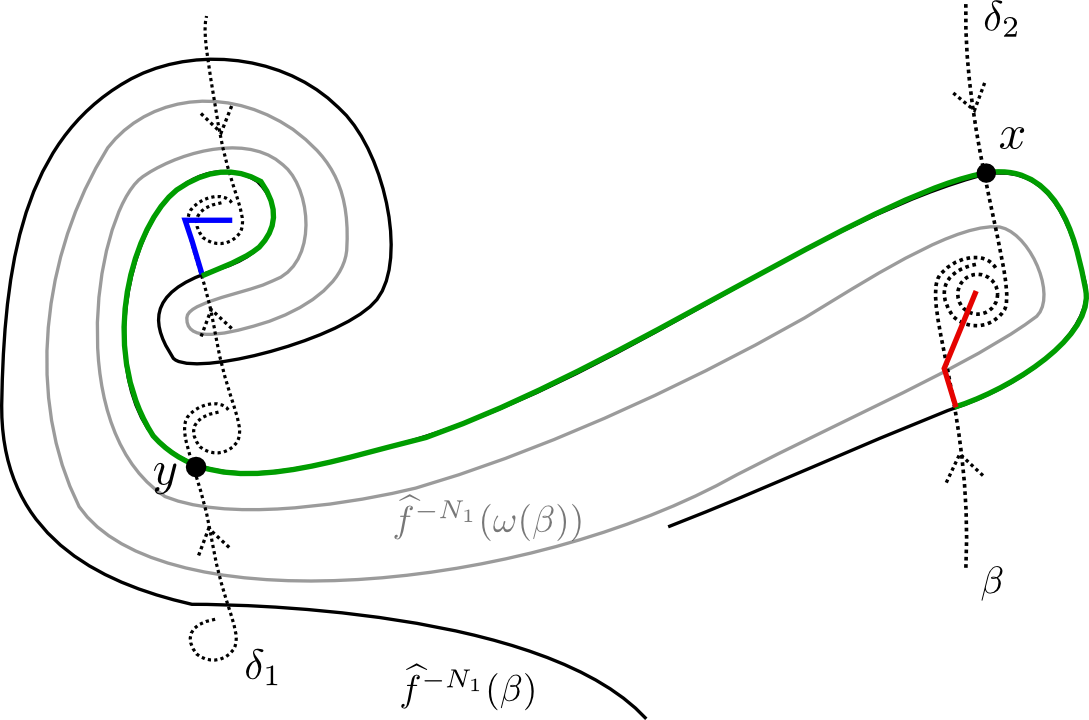}
\caption{An illustration of the arc $\Gamma$. In blue the arc $\wh{f}^{-N_1}\lambda$, in green $\beta^2$ and in red $\beta^1$.}
\label{fig.Gamma}
\end{center}  
\end{figure}

We postpone the proof of this lemma to the end of this section. Define
$$\Gamma = \wh{f}^{-N_1}\lambda \cdot(\beta^2)^{-1}\cdot(\beta^1)^{-1},$$
(cf. Fig. \ref{fig.Gamma}). Consider the set $U_{\Gamma}$, which is the (possibly empty) union of removable intervals for $\Gamma$ (cf. Definition \ref{def.removable}). Let $t_1$ be as in item 4 from Lemma \ref{lema.delta2}, and let $t_1'$ be such that $\Gamma(t_1') = \beta^2(t_1)\in\beta$. As $\beta^1\cap\cl{C}_{-2}=\empt$ (by item 2 of such lemma), it follows that $t'_1$ is not contained in $U_{\Gamma}$.

We may then define $x = \Gamma(t_2)$, where
$$ t_2= \min \{ t\in (0,1)\minus U_{\Gamma} \, : \, \Gamma(t) \in \cl{C}_0 \}. $$
Observe that as $\Gamma(t'_1)=\beta^2(t_1)\in\beta$,
\begin{equation}   \label{eq.c13}
\Gamma|_{[0,t_2]} \subset (\wh{f}^{-N_1}\lambda)\cdot (\beta^2)^{-1}.
\end{equation}
Let $\delta_2 \subset \cl{C}_0$ be the leaf of $\wh{\cl{F}}$ that contains $x$.

\begin{claim}    \label{cl.11}
$\Gamma$ arrives in $\delta_2$ in $t_2$ by the right.
\end{claim}

Before proving this claim, we use it to construct the arc $\delta$. Let $y = \Gamma(t_3)$, where
$$t_3= \max \{ t\in (0,t_2) \, : \, \Gamma(t) \in \cl{C}_{-2} \},$$
and define 
$$\delta = \Gamma_{[t_3,t_2]},$$
(reparametrized to $[0,1]$). Let $\delta_1$ be the leaf of $\wh{\cl{F}}$ that contains $y$. 

By definition of $t_3$ we have that $\txt{int}(\delta) \subset R(\cl{C}_{-2}) \cap L(\cl{C}_0)$ and $\delta \subset \beta^2 \subset \wh{f}^{-N_1}\beta$, and then $\delta$ satisfies items 1 and 2 from the Lemma \ref{lema.delta1}. By Claim \ref{cl.11} $\delta$ arrives in $\delta_2$ in 1 by the right and then $\delta$ satisfies item 4(b). Finally, the claim below will give us that $\delta$ satisfies item 3.

\begin{claim}    \label{cl.12}
$\delta $ leaves $\delta_1$ in $t=0$ by the right.
\end{claim}  

We proceed to the proof of Claims \ref{cl.11} and \ref{cl.12}.

\begin{prueba}[Proof of Claim \ref{cl.11}]
If the claim does not hold, then $\Gamma$ arrives in $\delta_2$ in $t_2$ by the left. By definition of $t_2$, $\Gamma|_{[0,t_2]}$ is homotopic wfe Rel sing$(\wh{\cl{F}})$ to an arc whose first and only intersection with $\cl{C}_0$ is in $x$, and such intersection is by the left. Then, by Lemma \ref{lemabrouwer2} we have that there is $t\in[0,t_2]$ such that $\wh{f}^{N_1}\Gamma(t)\in\delta_2$. By (\ref{eq.c13}) and by definition of $\lambda$, we have that 
$$\wh{f}^{N_1}(\Gamma|_{[0,t_2]})\cap \cl{C}_{0}  \subset \wh{f}^{N_1}(\wh{f}^{-N_1}(\lambda) \cdot (\beta^2)^{-1} ) \cap \cl{C}_0 \subset\beta.$$ 
This yields $\delta_2=\beta$. We will see that this leads to a contradiction.

Consider the canonical lift $\ol{f}$ of $\wh{f}$ to the universal cover of $\R^2\minus\txt{sing}(\beta)$. Fix a lift $\ol{\beta}$ of $\beta$, and consider the lift $\ol{\Gamma}$, of $\Gamma|_{[0,1)}$ such that $\ol{\Gamma}(t_2) \in \ol{\beta}$. As $\Gamma|_{[0,t_2]}$ is homotopic wfe Rel$(\txt{sing}(\wh{\cl{F}}))$ to an arc that intersects $\delta_2$ only in $x=\Gamma(t_2)$, and by our assumption that $\Gamma$ arrives in $\delta_2$ in $x$ by the left, we have
$$ \ol{\Gamma}(0) \in L(\ol{\beta}),$$
and then by Lemma \ref{lemabrouwer1} we have that $\ol{f}^{N_1}\ol{\Gamma}|_{[0,t_2]} \cap \ol{\beta} \neq\empt$. Also, by (\ref{eq.c13}) and by definition of $\lambda$, $\wh{f}^{N_1}\Gamma|_{[0,t_2]} \cap \beta \subset \wh{f}^{N_1}((\beta^2)^{-1}) \subset \beta$. Then, we must have that, if $\ol{\beta}^2$ is the lift of $\beta^2$ contained in $\ol{\Gamma}$,
$$ \ol{f}^{N_1}\ol{\beta}^2 \subset \ol{\beta}.$$
As $\ol{\Gamma}(t_2)\in \ol{\beta}^2$ and $\ol{\Gamma}(t_2)\in\ol{\beta}$, this contradicts the fact that $\ol{\beta}$ is a Brouwer curve for $\ol{f}$. This contradiction proves the claim.
\end{prueba}

\begin{prueba}[Proof of Claim \ref{cl.12}]
We suppose on the contrary that $\delta$ leaves $\delta_1$ in $t=0$ by the left, and we will find a contradiction. 
Consider the canonical lift $\ol{f}$ of $\wh{f}$ to the universal cover of $\R^2\minus\txt{sing}(\delta_1)$. Fix a lift $\ol{\delta}_1$ of $\delta_1$, and let $\ol{\Gamma}$ be the lift of $\Gamma|_{(0,1]}$ such that $\ol{\Gamma}(t_3) \in \ol{\delta}_1$. Note that, as $t_2 \notin U_{\Gamma}$, $t_3$ is not a removable intersection for $\Gamma$. By this and as $\Gamma$ leaves $\delta_1$ in $t_3$ by the left, we must have that 
$$\ol{\Gamma}(1) \in L(\ol{\delta}_1).$$
By Lemma \ref{lemabrouwer1} we have that $\ol{f}^{N_1}\ol{\Gamma}|_{[t_3,1]} \cap \ol{\delta}_1 \neq\empt$, and by item 2 from Lemma \ref{lema.delta2} $\wh{f}^{N_1}\beta^1\cap\cl{C}_{-2}=\empt$.  Thus 
$$ \wh{f}^{N_1}\beta_2\cap\cl{C}_{-2} \supset  \wh{f}^{N_1}(\beta^2\cdot\beta^1) \cap \cl{C}_{-2} \supset \wh{f}^{N_1} (\Gamma|_{[t_3,1]})\cap\delta_1 \neq\empt,$$
which contradicts the fact that $\wh{f}^{N_1}\beta^2\subset\beta\subset\cl{C}_0$. This is the sought contradiction, and this proves the claim. 
\end{prueba}

\begin{prueba}[Proof of Lemma \ref{lema.delta2}]
We have four possibilities. \\
\\
\textbf{Case 1:} $\omega(\beta)$ contains no singularities, and $\wh{f}^{-1}\omega(\beta) \subset \txt{Fill}(\omega(\beta))$.\\
Let $\beta_1^1$ be an arc contained in Fill$(\wh{f}^{-N_1}\omega(\beta))$, going from a singularity up to a point $p\in\wh{f}^{-1}\omega(\beta)$. 
Define $\beta_2^1$ as an arc going from $p$ to a point of $\wh{f}^{-N_1}\beta$, positively transverse to $\wh{f}^{-N_1}\beta$ and suficiently small in order that $\beta_2^1\subset \txt{Fill}(\omega(\beta))$ and $\beta_2^1\cup \wh{f}^{N_1}\beta_2^1 \cap \cl{C}_{-2}=\empt$ (see Fig. \ref{fig.delta2}). 

Define then $\beta^1=\beta_1^1\cdot\beta_2^1$, and note that
$$( \beta^1\cup \wh{f}^{N_1}\beta^1 ) \cap \cl{C}_{-2} = \empt,$$
and then $\beta^1$ satisfies items 1 and 2.

Define now $\beta^2$ as the subarc of $\wh{f}^{-N_1}(\beta)$ going from $\beta^1(1)$ to $\wh{f}^{-N_1}\lambda(1)\in\cl{C}_{-2}$. The arc $\beta^2$ satisfies then item (3) from the lemma. By last, we see that $\beta^2$ satisfies item (4). 

Let 
$$t_0 = \max \{ t \in[0,1] \, : \, \beta^2(t)\in \omega(\beta) \txt{ and } \beta^2(s)\in \txt{Fill}(\omega(\beta)) \ \txt{for all} \  s\in[0,t] \}.$$
Then, there exist points $t\in(t_0,1)$ arbitrarily close to $t_0$ such that $\beta^2(t) \notin \txt{Fill}(\omega(\beta))$, and such that $\beta^2(t)\in \beta$ (see Fig. \ref{fig.delta2}). Choose then $t_1\in(t_0,1)$ close enough to $t_0$ so 
$$\beta^2|_{[0,t_1]} \subset R(\cl{C}_{-2}).$$
The point $t_1$ and the arc $\beta^2$ satisfy then item (4), as we wanted. \\

\begin{figure}[h]        
\begin{center} 
\includegraphics{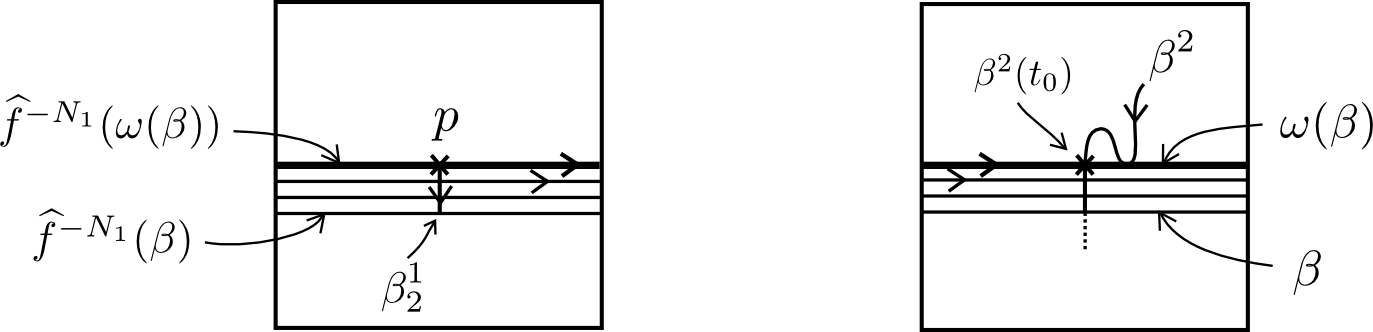}
\caption{Illustration of Case 1. Left: the arc $\beta^1_2$. Right: the arc $\beta^2$.}
\label{fig.delta2}
\end{center}  
\end{figure}

\noindent \textbf{Case 2:} $\omega(\beta)$ contains singularities and $\wh{f}^{-1}\omega(\beta) \subset \txt{Fill} (\omega(\beta))$. \\
By Theorem \ref{poinc.ben} we know that in this case, $\omega(\beta)$ is a generalized cycle of connections of $\wh{\cl{F}}$. Let $\epsilon$ be a leaf of $\wh{\cl{F}}$ contained in $\wh{f}^{-N_1}(\omega(\beta))$, and think of $\epsilon$ as an arc whose endpoints are singularities (both endpoints may coincide, and in this case $\e$ is a loop containing one singularity). Let $p\in\txt{int}(\epsilon)$, and let $\beta_1^1$ be the subarc of $\epsilon$ going from $\epsilon(0)$ to $p$ (see Fig. \ref{fig.delta3}). If $p$ is close enough to an endpoint of $\e$, then there is an arc $\beta_2^1$ going from $p$ to a point of $\wh{f}^{-N_1}\beta \cap \txt{Fill}(\omega(\beta))$, with $\beta_2^1$ sufficiently small in order that 
$$\wh{f}^{N_1}\beta_2^1\cap\cl{C}_{-2}=\empt.$$
Let $\beta^1=\beta_1^1\cdot\beta_2^1$. The arc $\beta^1$ satisies items 1 and 2. 

Define $\beta^2$ as the subarc of $\wh{f}^{-N_1}\beta$ going from $\beta^1(1)$ to $\wh{f}^{-N_1}\lambda(1)\in\cl{C}_{-2}$. The arc $\beta^2$ therefore satisfies item (3) from the lemma. By the same argument from Case 1 one can show that $\beta^2$ satisfies item (4). 
\vspace{1mm}
\begin{figure}[h]        
\begin{center} 
\includegraphics{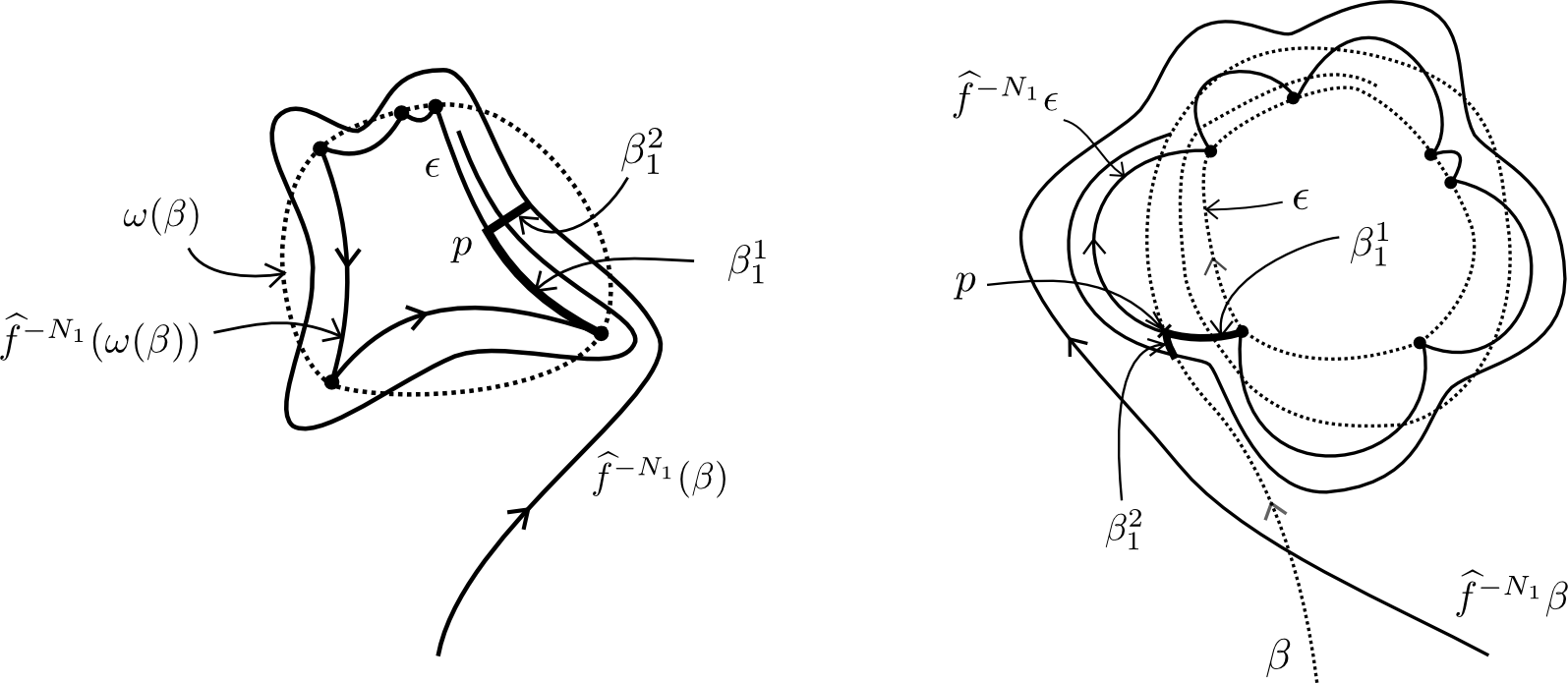}
\caption{Cases where $\omega(\beta)$ contains singularities. Left: Case 2, $\wh{f}^{-1}(\omega(\beta))\subset \txt{Fill}(\omega(\beta))$. Right: Case 3, $\wh{f}^{-1}(\omega(\beta))\supset \txt{Fill}(\omega(\beta))$.}
\label{fig.delta3}
\end{center}  
\end{figure}
\vspace{1mm}\\
\textbf{Case 3:} $\omega(\beta)$ contains singularities and $\txt{Fill}(\wh{f}^{-1}(\omega(\beta)) \supset \omega(\beta)$.\\
Let $\e$ be a leaf of $\wh{\cl{F}}$ contained in $\omega(\beta)$. Think of $\e$ as an arc whose endpoints $\e(0)$, $\e(1)$ are singularities (both endpoints may coincide). There exist points $p$ of $\beta$ arbitrarily close to $\e(0)$ and such that $\beta$ arrives in $\wh{f}^{-N_1}\e$ in $p$ by the left (see Fig. \ref{fig.delta3}). Let $\beta_1^1$ be a subarc of $\wh{f}^{-N_1}\e$ going from $\wh{f}^{-N_1}\epsilon(0)$ to one of such points $p$ sufficiently close to $\e(0)$ such that $\beta_1^1\cap\cl{C}_{-2}=\empt$. Let $\beta_2^1$ be a subarc of $\beta$ going from $p$ to a point in $\wh{f}^{-N_1}\beta$, with $\beta_2^1$ sufficiently small so that 
$$ \beta_2^1 \cup \wh{f}^{N_1}\beta_2^1\cap\cl{C}_{-2} = \empt. $$
Let $\beta^1=\beta_1^1\cdot\beta_2^1$. The arc $\beta^1$ satisfies then items 1 and 2. 

Define $\beta^2$ be the subarc of $\wh{f}^{-N_1}(\beta)$ going from $\beta^1(1)$ to $\wh{f}^{-N_1}\lambda(1)\in\cl{C}_{-2}$. The arc $\beta^2$ satisfies then item (3) from the lemma. Letting $t_1=0$, we clearly have that $t_1$ and $\beta^2$ satisfy item (4).\\
\\
\textbf{Case 4:} $\omega(\beta)$ has no singularities and $\txt{Fill}(\wh{f}^{-1}(\omega(\beta)))\supset \omega(\beta)$. \\
Let $\ell$ be a straight vertical line such that $\cl{C}_{-2}\subset L(\ell)$ and $\cl{C}_0\subset R(\ell)$. Let $D$ be the connected component of $\wh{f}^{-N_1}(\txt{Fill}(\omega(\beta)))\cap R(\ell)$ that contains $\omega(\beta)$ (see Fig. \ref{fig.delta4}).  
\begin{claim}  \label{claim.delta1}
$\wh{f}^{-N_1}(\txt{Fill}(\omega(\beta)))$ does not contain $\beta$. 
\end{claim}
With this claim we have that there exists $t_0\in[0,1]$ such that $\beta(t)\in D$ for all $t\in[t_0,1]$. Let $p=\beta(t_0)\in \pr D \cap \wh{f}^{-N_1}(\omega(\beta))$. Then, $\beta$ arrives in $\wh{f}^{-N_1}(\omega(\beta))$ in $p$ by the left. Let $\beta_1^1$ be an arc contained in $D$ going from a singularity of $\wh{\cl{F}}$ in the interior of $D$ to the point $p$. Let $\beta_2^1$ be a subarc of $\beta$ going from $p$ to a point of $\wh{f}^{-N_1}\beta$, with $\beta_2^1$ sufficiently small such that 
$$\beta_2^1\cup\wh{f}^{N_1}\beta_2^1\cap\cl{C}_{-2}=\empt.$$ 
Let $\beta^1=\beta_1^1\cdot\beta_2^1$. The arc $\beta^1$ satisfies items 1 and 2. 

Define $\beta^2$ as the subarc of $\wh{f}^{-N_1}\beta$ going from $\beta^1(1)$ to $\wh{f}^{-N_1}\lambda(1)\in\cl{C}_{-2}$. The arc $\beta^2$ then satisfies item (3) from the lemma. Letting $t_1=0$ we obtain that $t_1$ and $\beta^2$ satisfy item (4). 
\vspace{1mm}
\begin{figure}[h]        
\begin{center} 
\includegraphics{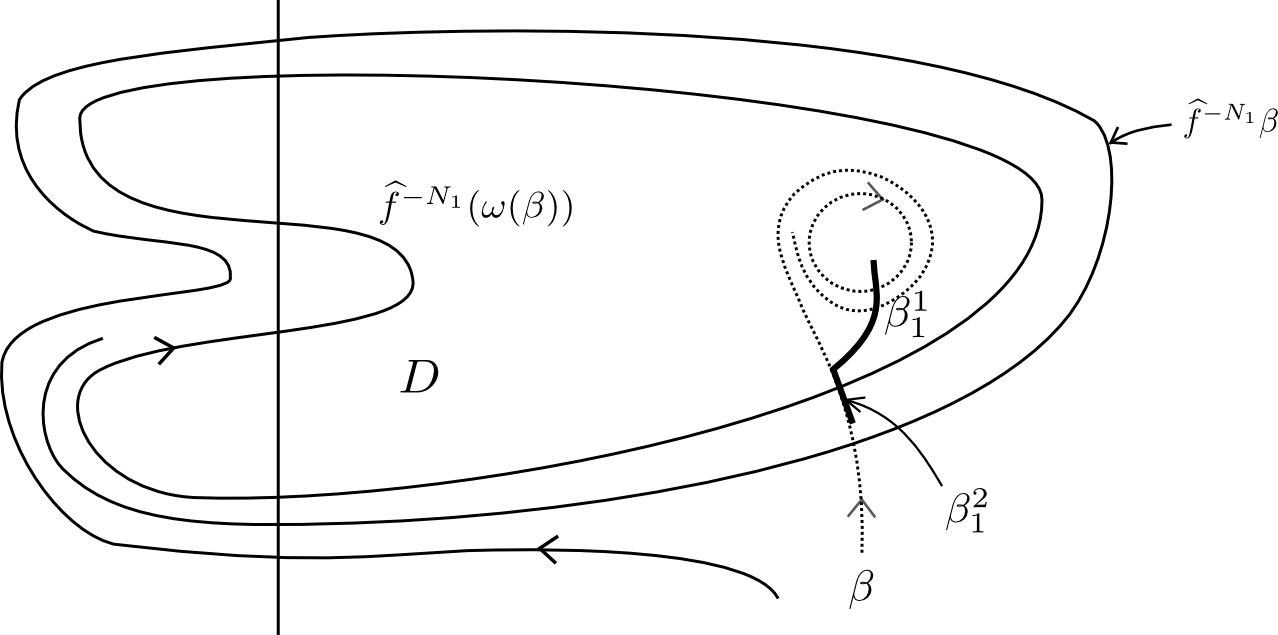}
\caption{Illustration of Case 4.}
\label{fig.delta4}
\end{center}  
\end{figure}
\begin{prueba}[Proof of Claim \ref{claim.delta1}]
As $\omega(\beta)$ contains no singularities, we must have $\alpha(\beta)\cap\omega(\beta)=\empt$. Suppose first that $\alpha(\beta)$ contains singularities, and denote the set of such singularities by $S$. As $S\subset\R^2\minus\txt{Fill}(\omega(\beta))$, and as $\wh{f}$ is isotopic to the identity, $S\subset\R^2\minus\wh{f}^{-N_1}(\txt{Fill}(\omega(\beta)))$, which implies that $\alpha(\beta) \nsubseteq \wh{f}^{-N_1}(\txt{Fill}(\omega(\beta)))$, and then $\beta \nsubseteq \wh{f}^{-N_1} \txt{Fill}(\omega(\beta))$.

Now suppose that $\alpha(\beta)$ does not contain singularities. In this case there is a non-empty set $S$ of singularities contained in $\txt{int}(\txt{Fill}(\alpha(\beta)))$. Also, $S\subset \R^2\minus \wh{f}^{-N_1}(\txt{Fill}(\omega(\beta)))$. If $\wh{f}^{-N_1}(\txt{Fill}(\omega(\beta)))$ contained $\beta$, it would also contain the loop $\alpha(\beta)$, and as $S \subset \txt{int}(\txt{Fill}(\alpha(\beta)))$, we would have that $\wh{f}^{-N_1}(\txt{Fill}(\omega(\beta)))$ is not simply connected, which is a contradiction. 
\end{prueba}
This finishes the proof of Lemma \ref{lema.delta2}.
\end{prueba}

\subsubsection{Complementary cases.}

\paragraph{ None of the sets $\alpha(\gamma)$ and $\omega(\gamma)$ is a singularity and at least one of the sets $\alpha(\beta)$ and $\omega(\beta)$ is a singularity.}\mbox{} \\

\noindent Recall that in this case, the arc $\lambda$ was defined in $\S$\ref{sec.plg} as $\lambda_1\cdot\lambda_2$, where $\lambda_1$ is straight horizontal, going leftwards from a singularity contained in $\txt{Fill}(\omega(\gamma))$ to a point of $\wh{f}^{N_1}(\gamma)$, and $\lambda_2$ is a subarc of $\wh{f}^{N_1}\gamma$ going from $\lambda_1(1)$ to a point $z\in\wh{f}^{N_1}\gamma\cap\cl{C}_0$ such that $\lambda_2|_{[0,1)}\subset L(\cl{C}_0)$. Also, $\beta\subset\cl{C}_0$ was defined as the leaf of $\wh{\cl{F}}$ such that $\lambda(1)\in\beta$.

Without loss of generality, we assume that $\omega(\beta)$ consists of a singularity $s$. Extend $\beta$ to $(0,1]$ as $\beta(1)=s$, and let $\beta^1\subset\beta$ be the subarc of $\beta$ going from $\lambda(1)$ to $\beta(1) = s$. Let 
$$ t_2 = \min \{ t\in (0,1] \minus U_{\wh{f}^{-N_1}\beta^1} \, : \, \wh{f}^{-N_1}\beta^1(t) \in \cl{C}_0 \}, $$
(cf. Definition \ref{def.removable}) and
$$ t_3 = \max \{ t\in [0,t_2) \, : \, \wh{f}^{-N_1}\beta^1(t) \in \cl{C}_{-2} \}.$$
With the same arguments from Section \ref{sec.sing}, one can prove that if
$$\delta= \wh{f}^{-N_1}\beta^1|_{[t_3,t_2]}$$
(reparametrized to $[0,1]$), then $\delta$ leaves a leaf $\delta_1\subset\cl{C}_{-2}$ of $\wh{\cl{F}}$ in $t=0$ by the right, int$(\delta)\subset R(\cl{C}_{-2})\cap L(\cl{C}_0)$, and either $\delta$ arrives in a leaf $\delta_2\subset\cl{C}_0$ of $\wh{\cl{F}}$ in $t?1$ by the right, or $\delta(1) = s \in\txt{sing}(\wh{\cl{F}})$, and then $\delta$ satisfies items 1 to 4 from Lemma \ref{lema.delta1}.

\paragraph{ None of the sets $\alpha(\beta)$ and $\omega(\beta)$ is a singularity and at least one of the sets $\alpha(\gamma)$ and $\omega(\gamma)$ is a singularity.} \mbox{} \\

\noindent Without loss of generality, assume that $\omega(\gamma)$ consists of a singularity $s$. Extend $\gamma$ to $(0,1]$ as $\gamma(1)=s$. In this case, the arc $\beta\subset\cl{C}_0$ was defined in $\S$\ref{sec.plg} as the leaf of $\wh{\cl{F}}$ that contains $\gamma(t_*)$, where 
$$t_*= \max \{ t\in (0,1) \, : \, \wh{f}^{N_1}\gamma(t) \in \cl{C}_0 \},$$
and the arc $\lambda$ was defined as $(\wh{f}^{N_1}\gamma|_{[t_*,1]})^{-1}$.  

The proof of Lemma \ref{lema.delta2}, unmodified, gives us two arcs $\beta^1$ and $\beta^2$ satisfying the conclusions of that lemma. Let $\Gamma=\wh{f}^{-N_1}\lambda\cdot (\beta^2)^{-1}\cdot (\beta^1)^{-1}$ and define
$$ t_2 = \min \{ t\in (0,1] \minus U_{\Gamma} \, : \, \Gamma(t) \in \cl{C}_0 \}, $$
(cf. Definition \ref{def.removable}), and
$$ t_3 = \max \{ t\in [0,t_2) \, : \, \Gamma(t) \in \cl{C}_{-2} \}.$$
The same arguments from Section \ref{sec.nosing} give us that if
$$\delta= \Gamma|_{[t_3,t_2]}$$
(reparametrized to $[0,1]$), then the arc $\delta$ satisfies the conclusions of Lemma \ref{lema.delta1}.

\addcontentsline{toc}{section}{References}

\bibliographystyle{amsalpha}
\bibliography{bib-sublinearA}

\providecommand{\bysame}{\leavevmode\hbox to3em{\hrulefill}\thinspace}
\providecommand{\MR}{\relax\ifhmode\unskip\space\fi MR }
\providecommand{\MRhref}[2]{%
  \href{http://www.ams.org/mathscinet-getitem?mr=#1}{#2}
}
\providecommand{\href}[2]{#2}
\begin{thebibliography}{BCJR09}

\bibitem[Atk76]{atk}
G.~Atkinson, \emph{Recurrence of cocycles and random walks}, Journal of the
  London Mathematical Society \textbf{13} (1976), no.~2, 486--488.

\bibitem[BCJR09]{bcjr}
F.~B\'eguin, S.~Crovisier, T.~J{\"a}ger, and F.~Le Roux, \emph{Denjoy
  constructions for fibered homeomorphisms of the torus}, Trans. Amer. Math.
  Soc. \textbf{361} (2009), no.~11, 5851--5883.

\bibitem[Cai51]{cairns}
S.~Cairns, \emph{An elementary proof of the {J}ordan-{S}choenflies theorem},
  Proc. of the Amer. Math. Soc. \textbf{91} (1951), no.~2, 860--867.

\bibitem[Cal05]{lc2}
P.~Le Calvez, \emph{Une version feuillet\'ee \'equivariante du th\'eor\`eme de
  translation de {B}rouwer}, Publications Math\'ematiques de l'IH\'ES
  \textbf{102} (2005), no.~1, 1--98.

\bibitem[Dav86]{daver}
R.~J. Daverman, \emph{Decompositions of manifolds}, vol. 124, Academic Press
  Inc., Orlando, FL, 1986.

\bibitem[D{\'a}v13]{dav}
P.~D{\'a}valos, \emph{On torus homeomorphisms whose rotation set is an
  interval}, Mathematische Zeitschrift (2013), DOI
  10.1007/s00209--013--1168--3.

\bibitem[Fay02]{fay1}
B.~Fayad, \emph{Weak mixing for reparametrized linear flows on the torus},
  Ergodic Theory and Dynamical Systems \textbf{22} (2002), 187--201.

\bibitem[FM90]{fm}
J.~Franks and M.~Misiurewicz, \emph{Rotation sets of toral flows}, Proc. Amer.
  Math. Soc. \textbf{109} (1990), no.~1, 243--249.

\bibitem[Fra89]{f2}
J.~Franks, \emph{Realizing rotation vectors for torus homeomorphisms},
  Transactions of the American Mathematical Society \textbf{311} (1989), no.~1,
  107--115.

\bibitem[Fra95]{f3}
\bysame, \emph{The rotation set and periodic points for torus homeomorphisms},
  Dynamical Systems \& Chaos (Aoki, Shiraiwa, and Takahashi, eds.), World
  Scientific, Singapore (1995), 41--48.

\bibitem[Fur61]{fur}
H.~Furstenberg, \emph{Strict ergodicity and transformation of the torus}, Amer.
  J. of Math. \textbf{83} (1961), 573--601.

\bibitem[GKT12]{gkt}
N.~Guelmanm, A.~Koropecki, and F.~A. Tal, \emph{A caracterization of annularity
  for area preserving toral homeomorphisms}, eprint arXiv:1211.5044v1 (2012).

\bibitem[Gut79]{Gut}
C.~Guti\'errez, \emph{Smoothing continuous flows and the converse of
  {D}enjoy-{S}chwartz theorem}, An. Acad. Brasil. Ci\^enc. \textbf{51} (1979),
  no.~4, 581--589.

\bibitem[Han89]{handel}
M.~Handel, \emph{Periodic point free homeomorphisms of $\mathbb{T}^2$}, Proc.
  of the Amer. Math. Soc. \textbf{107} (1989), no.~2, 511--515.

\bibitem[Han90]{han}
\bysame, \emph{The rotation set of a homeomorphism of the annulus is closed},
  Communications in Mathematical Physics \textbf{127} (1990), 339--349.

\bibitem[J{\"a}g09a]{j2}
T.~J{\"a}ger, \emph{The concept of bounded mean motion for toral
  homeomorphisms}, Dyn. Sys. \textbf{24} (2009), no.~3, 277--297.

\bibitem[J{\"a}g09b]{j1}
\bysame, \emph{Linearisation of conservative toral homeomorphisms}, Invent.
  Math. \textbf{176} (2009), no.~3, 601--616.

\bibitem[Jau13]{jau}
O.~Jaulent, \emph{Existence d'un feuilletage positivement transverse \`a un
  homeomorphism de surface}, arXiv:1206.0213 (2013).

\bibitem[KK08]{kk}
A.~Kocksard and A.~Koropecki, \emph{Free curves and periodic points for torus
  homeomorphisms}, Ergodic Theory \& Dynamical Systems \textbf{28} (2008),
  1895--1915.

\bibitem[KK09]{kk2}
\bysame, \emph{A mixing-like property and inexistence of invariant foliations
  for minimal diffeomorphisms of the 2-torus}, Proc. of the Amer. Math. Soc.
  \textbf{137} (2009), no.~10, 3379--3386.

\bibitem[KT]{st}
A.~Koropecki and F.~A. Tal, \emph{Strictly toral dynamics}, To appear in
  Invent. Math.

\bibitem[KT12a]{kt1}
\bysame, \emph{Area preserving irrotational diffeomorphisms of the torus with
  sublinear diffusion}, To appear in Proc. Amer. Math. Soc. (2012).

\bibitem[KT12b]{kt2}
\bysame, \emph{Bounded and unbounded behaviour for area-preserving rational
  pseudorotations}, Preprint, arXiv:1207.5573v3 (2012).

\bibitem[Kwa92]{kw2}
J.~Kwapisz, \emph{Every convex polygon with rational vertices is a rotation
  set}, Erg. Theory and Dyn. Sys. \textbf{12} (1992), 333--339.

\bibitem[Kwa95]{kw}
\bysame, \emph{A toral diffeomorphism with a nonpolygonal rotation set},
  Nonlinearity \textbf{8} (1995), no.~4, 461--476.

\bibitem[LM91]{lm}
J.~Llibre and R.~S. Mackay, \emph{Rotation vectors and entropy for torus
  homeomorphisms isotopic to the identity}, Ergodic Theory \& Dynamical Sistems
  \textbf{11} (1991), 115--128.

\bibitem[Mat97]{mat}
S.~Matsumoto, \emph{Rotation sets of surface homeomorphisms}, Bol. Soc. Bras.
  Mat. \textbf{28} (1997), no.~1, 89--101.

\bibitem[Mos65]{moser}
J.~Moser, \emph{On the volume elements on a manifold}, Trans. of the Amer.
  Math. Soc. \textbf{120} (1965), no.~2, 286--294.

\bibitem[MZ89]{mz}
M.~Misiurewicz and K.~Ziemian, \emph{Rotation set for maps of tori}, Journal of
  the London Mathematical Society \textbf{40} (1989), no.~2, 490--506.

\bibitem[NS89]{NemStep}
V.~V. Nemitskii and V.~V. Stepanov, \emph{Qualitative theory of differential
  equations}, Courier Dover Publications, 1989.

\bibitem[Poi52]{po}
H.~Poincar\'e, \emph{Oeuvres compl\`etes}, {G}authier-{V}illars, Paris, 1952.

\bibitem[Pol92]{pol}
M.~Pollicott, \emph{Rotation sets for homeomorphisms and homology}, Trans.
  Amer. Math. Soc. \textbf{331} (1992), no.~2, 881--894.

\bibitem[Sch57]{schw}
S.~Schwartzmann, \emph{Assymptotic cycles}, Ann. of Math. \textbf{66} (1957),
  no.~2, 270--284.

\bibitem[Sol45]{soln}
G.~Solntzev, \emph{On the asymptotic behaviour of integral curves of a system
  of differential equations}, Bull. Acad. Sci. URSS. S\'er. Math. [Izvestia
  Akad. Nauk SSSR] \textbf{9} (1945), 233--240.

\bibitem[Tal13]{tal1}
F.~A. Tal, \emph{On non-contractible periodic orbits for surface
  homeomorphisms}, e-print: arXiv:1307.1664 (2013).

\bibitem[Whi33]{whi33}
H.~Whitney, \emph{Regular families of curves}, Ann. of Math. \textbf{34}
  (1933), no.~2, 244--270.

\bibitem[Whi41]{whi41}
\bysame, \emph{On regular families of curves}, Bull. Amer. Math. Soc.
  \textbf{47} (1941), 145--147.

\end{thebibliography}

\noindent \textsc{Instituto Tecnol\'ogico y de Estudios Superiores de Occidente, Perif\'erico Sur Manuel G\'omez Mor\'in 8585, C.P. 45604, Tlaquepaque, Jalisco, M\'exico.}\\
\ttfamily mailto: pablod@iteso.mx

\end{document}